\documentclass[12pt,twoside]{article}

\newcommand{\beq}{\begin{equation}}
\newcommand{\eeq}{\end{equation}}

\newtheorem{lem}{Lemma}[section]
\newtheorem{cor}{Corollary}[section]
\newtheorem{theor}{Theorem}[section]

\newtheorem{prop}{Proposition}[section]

\usepackage{amsmath}

\numberwithin{equation}{section}

\begin{document}

\title{Asymptotic Behavior of Subfunctions of
the Stationary Schr\"{o}dinger Operator}
\author{Boris Ya. Levin \\
{\small deceased} \vspace{.5cm} \\
Alexander I. \vspace{-.2cm} Kheyfits}

\date{{\small Bronx Community College/CUNY} \vspace{.2cm} \\
{\small alexander.kheyfits@bcc.cuny.edu }}

\maketitle

\vspace{.5cm}

\begin{abstract}
Subharmonic functions associated with the stationary
Schr\"{o}dinger operator are weak solutions of the inequality
$\Delta u - c(x)u \geq 0$ under appropriate assumptions on the
potential $c$. We derive for these functions analogs of several
classical results on analytic and subharmonic functions such as
the Phragm\'{e}n-Lindel\"{o}f theorem with the precise growth
rate, the Blaschke theorem on bounded analytic functions, the
Carleman formula, the Hayman-Azarin theorem on the asymptotic
behavior of special subharmonic functions in many-dimensional
cones. \vspace{4cm}
\end{abstract}

\footnoterule \vspace{.5cm}


\newpage

\tableofcontents

\newpage

\section{Introduction}

Some essential properties of holomorphic functions are valid for
F. Riesz subharmonic functions, that is, for locally-summable
solutions of the inequality $ \Delta u \geq 0 $, where $ \Delta =
\frac{\partial ^2}{\partial x^2 _1} + \cdots + \frac{\partial
^2}{\partial x^2 _n}$ is the Laplace operator - see, for example,
the monograph of W. Hayman and P. Kennedy  \cite{HK}. In many
issues the Laplacian can be replaced by a more general partial
differential operator. \\

We study from this point of view \emph{ generalized subharmonic
functions}, called here $c-$\emph{ subfunctions}, associated with
a stationary Schr\"{o}dinger operator
\begin{equation}
L_{c} = - \Delta + c(x) I
\end{equation}
with a potential $c(x)$, $I$ being an identical operator. These
are subsolutions of the operator (1.1), that is, weak solutions of
the inequality
\[\Delta u(x) - c(x)u(x) \geq 0.\]
An equivalent definition, similar to F. Riesz' definition of the
subharmonic functions, and exact assumptions on the potential
$c(x)$ will be stated in Section 2.

In this work we extend several results, concerning analytic and
subharmonic functions, onto the generalized subharmonic functions.
In Section 3 we generalize the Phragm\'{e}n-Lindel\"{o}f theorem
with the precise estimate of the limit growth rate. A general
Phragm\'{e}n-Lindel\"{o}f principle claims that, given a class
$\textrm{F}$ of functions in a domain $\Omega$ with a
distinguished boundary point $\xi \in \partial \Omega$, there
exists the limit rate, depending on the domain $\Omega$ and the
class $\textrm{F}$, such that if a function $f \in \textrm{F}$ is
upper bounded on $\partial \Omega \setminus \{\xi \}$, then either
$f$ is upper bounded throughout the entire domain $\Omega$ or else
$f(x)$ must grow at least with this limit rate when $\Omega \ni x
\rightarrow \xi$. In Theorem 3.1 we consider $c-$subfunctions in a
cone $K^D$ generated by a domain $D$ on the unit sphere and show
that the limit growth is given by an increasing solution of the
equation
\begin{equation}y''(r)+(n-1)r^{-1}y'(r)-\left(\lambda _0
r^{-2}+q(r)\right)y(r) = 0, \end{equation} where $\lambda _0 $ is
the smallest eigenvalue of the Laplace-Beltrami operator in $D$
and $q$ is any nonnegative radial minorant of the potential $c$.

With regard to this theorem, we want to mention the following
substantial circumstance. It turns out (under some mild regularity
conditions - see Section 3) that if $|x|^2 c(x) \rightarrow \infty
$ as $\rightarrow \infty $, then the lowest possible growth in a
cone $K^D$ of a subfunction $u$, which is upper bounded at the
boundary of the cone, is exponential rather than power as in the
case of the Laplacian. Moreover, this lowest growth does not
depend on the geometry of the spherical domain $D$ and coincides
with the precise growth rate in the Liouville theorem for
$c-$harmonic functions in the entire space. A version of the
Liouville theorem for $c-$harmonic functions is also proven in Section 3.\\

In Section 5 we prove a "sub"-analog of the Blaschke theorem on
bounded analytic functions. The theorem asserts that if $z_k, \;
k=1,2, \ldots $, are the zeros of a bounded analytic function in
the unit disk in the complex plane, then $\sum _k
(1-|z_k|)<\infty$. Our generalization of this theorem involves the
Riesz measure $d\mu$ of a subfunction $u$ in a cone $K^D$ and
requires the convergence of two integrals
\[\int _{K^D _1}V(r)\varphi _0(\theta )d\mu (r,\theta)<\infty \]
and
\[\int _{K^D \setminus K^D _1}W(r) \varphi _0(\theta )d\mu (r,\theta)<\infty ,\]
where $V$ and $W$ are increasing and decreasing, respectively,
solutions of the ordinary differential equation (1.2), $\varphi
_0(\theta )$ is a positive eigenfunction of the Laplace-Beltrami
operator in $D$, and $K^D _1$ is the truncated cone $K^D \cap
B(0,1)$. To prove this theorem, we extend the Carleman formula,
well-known in the case of analytic functions, to the
Schr\"{o}dinger operator $L_c$. This generalization is of its own
interest and can be carried over more general classes of
operators. The existence and asymptotic properties of special
solutions $V$ and $W$ are discussed in Appendix C. \\

In Section 6 we are concerned with the W. Hayman - V. Azarin
theorem \cite{Azar} on the asymptotic behavior of subharmonic
functions in $n-$dimensional cones. The original theorem of W.
Hayman \cite{Hay} states that if $u$ is a subharmonic function in
a half-plane, then under certain restrictions on its growth there
exists the limit $\lim _{r \rightarrow \infty } r^{-1} u(r  e^{i
\theta })$, provided that $z= r e^{i \theta}$ approaches infinity
avoiding some small exceptional set. In particular, this means
that $u$ cannot decay too fast when $z \rightarrow \infty$. We
generalize this statement
onto subfunctions in many-dimensional cones.\\

In the proofs we use the bilinear series representation and
estimates for Green's function of the operator $L_q$ in cones
derived in Section 4. These estimates generalize the known ones
for the Laplacian \cite{Azar}. In Sections 3 - 6 we deal with
subfunctions in cones. In Section 7 we restate our results for the
subfunctions in tube domains. Our results can be also extended
onto the case of a general second-order self-adjoint elliptic
operator
\[- \sum^{n}_{i, j = 1} \frac {\partial }{\partial x_{i}} (a_{ij}(x)
\frac {\partial u}{\partial x_{j}} ) + \sum^{n}_{i= 1} b_i (x)
\frac {\partial u}{\partial x_i} + c(x)u(x)\]
with sufficiently smooth coefficients $a_{ij}$ and $b_i$.  \\

Our proofs essentially use not only representation and estimates
but also some other properties of the Green function $G(x;y)$ of
the operator $L_c$ with the Dirichlet boundary conditions,
properties of the eigenvalues and eigenfunctions of the Laplace-
Beltrami operator $\Delta ^*$ in domains on the unit sphere
$S^{n-1} \subset \mathbf{R}^n$, and asymptotic properties of
radial solutions of the operator $L_q$ with a radial potential
$q=q(r)$, that is, solutions of (1.2). These results are mostly
known, but we often need them under less restrictions on the
smoothness of boundaries and potentials, than we were able to
locate in the literature. Due to this reason and for the sake of
completeness we have included proofs of some of these results. We
presented them in three appendices. Properties of Green's function
of the operator $L_c$ are studied in Appendix A. In particular, we
prove\footnote{The history of proof of the existence of Green's
function can be traced back to the article of E. Levy \cite{LeE},
who reduced this proof to solving an integral equation. The method
is described, for instance, in \cite{John} or \cite{Mir}, Chap.1.
However, Levy's proof contained a non-justified point. Apparently,
first complete proofs of the existence of Green's function for a
linear second-order elliptic operator with smooth coefficients
were given independently by P. Garabedian and M. Shiffman
\cite{GarSh} and Yu. Lyubich \cite{Lyu1, Lyu2}; see also P. Lax
\cite{Lax0}.} that for the class of potentials under consideration
Green's function of $L_c$ exists simultaneously with Green's
function of the Laplacian, namely, in any domain whose boundary is
not a polar\footnote{Terms \emph{ capacity, polar set, regular}
and \emph{ irregular} points are used here in the sense of
classical potential theory - see, for example, \cite{HK}.} set
(Cf. \cite{HK}, Theorem 5.24).\\

It is worth mentioning that these results have clear physical
interpretation - see Remark A.4 in Appendix A.\\

In Appendix B we prove the existence of the eigenfunctions of the
Laplace-Beltrami operator in any domain on $S^{n-1}$ whose
boundary is not a polar set. Appendix C contains necessary
properties of radial solutions of the operator $L_q$, that is,
solutions of the ordinary differential equation (1.2).

\vspace{1cm} The first draft of this article was written in
1985-1986, main results were published without proofs in
\cite{LK}. Many circumstances have since then delayed completing
the work. On August 24, 1993  Boris Yakovlevich Levin, a man of
great personality and a brilliant mathematician, passed away.
Despite the obvious futility of carrying out this task without
B. Ya., while finishing the paper, I have tried to follow his lessons. \\

\hfill  Alexander Kheyfits

\newpage

\section{Subfunctions of the operator $L_{c}$}

A generic point of the Euclidean $n-$dimensional space
$\mathbf{R}^n$ is denoted by $x = (x_{1}, x_{2}, \ldots ,x_{n})$.
The ball and the sphere of radius $r$ centered at the point $x$
are denoted by $B(x, r)$ and $S(x, r) = \partial B(x, r)$; also,
$B(r) = B(0,r)$ and $S(r) = S(0, r)$; $g=g(x;y)$ always stands for
Green's function of the Laplacian $- \Delta$. The boundary and the
closure of a domain $\Omega$ are denoted, respectively, by
$\partial \Omega $ and $\overline{\Omega}$, $a^+ = \max \{a; \; 0
\}$, $ \frac{\partial }{\partial n} $ always denotes the
derivative with respect to the inward unit normal. The end of a
proof or a statement, remark, etc., is denoted by a lozenge
$\diamondsuit$. Generic constants, which may be different from
point to point, are denoted by $b$, sometimes with subscripts.\\

Now we define a class of potentials $c(x)$ under consideration.
Hereafter, the potential $c$ is supposed to be nonnegative and
locally summable in a domain
$ \Omega \subset \mathbf{R}^n, \; n \geq 2 $, \\
\begin{equation}
c(x) \geq 0,
\end{equation}
\begin{equation}
c(x) \in L^{p}_{loc} ( \Omega )
\end{equation}
with an exponent $p>n/2$ if $n \geq 4$, and with $p=2$ if $n=2$ or
$3; \; \parallel c \parallel$ always means $\parallel c \parallel
_{L^p}$. In some cases (2.1) can be relaxed to the nonoscillatory
condition $4|x|^2 c(x) \geq - (n-2)^2$, but we do not pursue this
issue here. It should be noted that if $\Omega $ is unbounded, the
condition (2.2) imposes no restriction on the asymptotic behavior
of the potential $c(x)$ as $|x| \rightarrow \infty $.

The class of potentials satisfying the conditions (2.1)-(2.2) in a
domain $ \Omega \subset \mathbf{R}^n$ is denoted by ${\cal
C}(\Omega )$.  If $c \in {\cal C}(\Omega )$, then the differential
expression (1.1) can be extended in a standard way  to an
essentially self-adjoint operator $L_c$ on $L^2 (\Omega )$; $L_c$
always refers to this extended operator. The operator $L_c$ has
the positive Green function $G(x;y)$ with Dirichlet boundary
conditions;
its properties are discussed in Appendix A. \\

The classical F. Riesz' definition of  subharmonic functions has
been carried over the operators $L_{c}$ with continuous potentials
by Dinghas \cite{Ding} and Myrberg \cite{Myr}. \pagebreak \\
\noindent \textbf{Definition.} An upper semicontinuous function
\[u: \Omega \longrightarrow [-\infty , \infty ), \; u(x) \not\equiv -
\infty ,\] in an arbitrary domain $\Omega \subset \mathbf{R}^n$ is
called a \emph{ $c-$subfunction} or a {\em generalized subharmonic
function} (g.s.f.) associated with the operator $L_c$ if and only
if $u$ satisfies the generalized mean-value
inequality\footnote{The corresponding inequality in \cite{Ding}
differs from (2.3) with a multiplicative constant due to another
normalization of Green's function at its singularity.}
\begin{equation}
u(x) \leq {\cal M}^{G} (u,x,r)
\end{equation}
at each point $x \in \Omega$, where
\begin{equation}
{\cal M}^{G}(u,x,r) \equiv  \int _{S(x,r)} u(y) \frac {\partial
G_{r}(x; y)}{ \partial n(y) } \; d \sigma (y)
\end{equation}
and $G_r$ stands for the Green function of $L_c$ in the ball
$B(x,r)$.

The class of $c-$subfunctions in a domain $\Omega $ is denoted by
$\verb"SbH"(c, \Omega )$. If $ - u \in \verb"SbH"(c, \Omega)$,
then we call $u$ \emph{ a $c -$superfunction} and denote the class
of $c-$superfunctions by $\verb"SpH"(c, \Omega)$. If a function
$u$ is both sub- and superfunction, it is, clearly, continuous and
is called a $c-$harmonic function or a \emph{ generalized harmonic
function } (g.h.f.) associated with the operator $L_c$; the class
of g.h.f. is denoted by $\verb"H"(c, \Omega) = \verb"SbH"(c,
\Omega) \cap \verb"SpH"(c, \Omega)$. In terminology we follow
Beckenbach \cite{Beck} and Nirenberg \cite{Nir}. Other authors
have in the similar circumstances used various terms such as {\em
generalized convex functions}, \emph{ weakly} $L-$\emph{
subharmonic functions}, $c-${\em subharmonic functions}, {\em
metaharmonic} and {\em submetaharmonic functions}, {\em
subsolutions}, {\em subelliptic functions}, {\em panharmonic
functions}, etc., - see, for example, Courant \cite[p. 342]{Cour},
Duffin \cite{Duff}, Littman \cite{Litt1}, Topolyansky
\cite{Topol}, Vekua \cite{Vek}.

Some properties of generalized harmonic and subharmonic functions
with H\"{o}lder or summable potentials have been already studied -
in addition to the papers cited above, see, for instance, works
\cite{BraD, Bram, Her2, Kato, Khe73, KhePAMS, KhePOTA1, KheDIE,
KheSampl, KheCV, KhePOTA2, Litt2, Maed, Russ, Sim} and references therein. \\

We state a few properties of the subfunctions in the following
theorem, omitting proofs, which are similar to the case $c=0$.\\
\begin{theor}\hspace{-.09in}\textbf{.} Let $c \in {\cal C}(\Omega)$. Then \\

$1^{o}$ Green's function $G(x;y)$ of the operator $L_c$ is a
$c-$superfunction in $\Omega$ and a $c-$harmonic function in
$\Omega \setminus \{y\}$. \\

$2^{o}$ If $u \in \verb"SbH"(c, \Omega)$, then $u^{+} \in
\verb"SbH"(c, \Omega)$ and $bu \in \verb"SbH"(c, \Omega)$ for any
constant $b \geq 0$, however,
$bu \in \verb"SpH"(c, \Omega)$ if $b \leq 0$. \\

$3^{o}$ Let $c_k (x) \in {\cal C}(\Omega), \; k=1,2$, and $c_1
\leq c_2$ in $\Omega $. If $u \in \verb"SbH"(c_2, \Omega)$  and
$u(x) \geq 0$, then $u \in \verb"SbH"(c_1, \Omega)$; vice versa,
if $u \in \verb"SbH"(c_1, \Omega)$ and $u(x) \leq 0$, then $u \in
\verb"SbH"(c_2, \Omega)$. In particular, for any $c \in {\cal
C}(\Omega)$ every nonnegative $c-$subfunction is Riesz'
subharmonic function and every nonpositive Riesz' subharmonic
function is a $c-$subfunction for any $c\in {\cal C}(\Omega)$.\\

$4^{o}$ The maximum principle holds for subfunctions as follows
(Cf., for instance, \cite[p. 3, Exercise 1.1]{Landis}: \\

If a subfunction $u$ is upper bounded in a bounded domain $\Omega$
and
\[\limsup _{\Omega \ni x \rightarrow \xi }u(x) \leq M =
const \] at every boundary point $\xi \in \partial \Omega
\setminus E$, where $E \subset \partial \Omega $ is an exceptional
(maybe empty) polar set, then $u(x) \leq M^{+}$ for all $x \in
\Omega$. Moreover, if $\sup _{x \in \Omega } u(x) = M \geq 0$ and
$u (x_0) = M$  at even one point $x_0 \in \Omega$, then $u(x)
\equiv M, \; \forall x \in \Omega $. The conclusion fails
if we replace here $M^+$ with $M$.\\

$5^{o}$ The principle of $c-$harmonic majorant holds for the
subfunctions in the following form: \\

Let $u\in \verb"SbH"(c, \Omega)$ and $v \in \verb"H"(c, \Omega)
\cap C(\overline{\Omega})$ in a bounded domain $\Omega$. If an
inequality $\lim \sup _{\Omega \ni x \rightarrow \xi }u(x) \leq
v(\xi)$ holds at every boundary point $\xi \in \partial \Omega $,
then $u(x) \leq v(x)$ everywhere in $\Omega $. Moreover, an
equation $u(x_{0}) = v(x_{0})$, if only at one point $ x_{0} \in
\Omega $,
implies $u(x) \equiv v(x), \; \forall x \in \Omega $. \\

Vice versa, let $ u(x) $ be an upper semicontinuous function in $
\Omega $ such that for each ball $B(x_{0},r) \subset \Omega$ and
for any function $v \in \verb"H"(c, B(x_{0},r)) \cap C(\overline
{B(x_{0},r)})$, an inequality $\lim \sup _{B(x_{0},r) \ni y
\rightarrow x}u(y) \leq v(x)$ for all $x \in S(x_{0},r)$ implies
an inequality $u(x) \leq v(x)$ everywhere in $B(x_{0}, r)$. Then
$u \in \verb"SbH"(c, \Omega)$. \\

$6^{o}$ In the same way as the Riesz' subharmonic functions,
$c-$subfunctions are locally summable in $\Omega $.
Moreover\footnote{\cite{Sim}, the limiting case $\beta = 2-n/p$
was studied in \cite{KhePAMS}.}, if $u \in \verb"H"(c, \Omega)$
where $c \in L^p _{loc}(\Omega), n/2<p$, then $u$ is a
H\"{o}lder-continuous function with the index $\beta = \min
\{2-\frac{n}{p}; \; 1\}$ on each compact set $K \subset \Omega$;
in addition, if $p>n$, then $\nabla u$ is a H\"{o}lder-continuous
function with the index $\beta = 1-n/p$
on compact sets.\\

$7^{0}$ As in the classical case $c(x)=0$, the definition (2.3) is
equivalent to a differential inequality $L_{c}u(x) \leq 0$
(\cite{HK}, see also \cite{ItoS}). If $u \in C^{2}(\Omega)$, then
the latter inequality holds pointwisely almost everywhere in
$\Omega$. If $u$ is merely locally summable, then that inequality
is fulfilled in the sense of distributions,
therefore, $-L_{c}u$ is a positive generalized function, that is, a measure. \\

$8^{0}$ As a corollary, we obtain the following Riesz
representation (2.5) of subfunctions. Hereafter, the constant
$\theta _{2} = 2 \pi \; $, $\theta _{n} =(n-2) \sigma _{n-1}$ for
$n \geq 3$, and $\sigma _{n-1}$ is the surface area of the unit
sphere $S$ in $\mathbf{R}^n$. \\

For each $u \in \verb"SbH"(c, \Omega)$ there exists a unique
measure
\[ d \mu (x) = - \frac{1}{ \theta _{n}} L_{c}u \]
with the support in $\Omega $, such that for every subdomain $
\Omega _{0} \subset \subset \Omega $ there exists a $c-$harmonic
function $v \in \verb"H"(c, \Omega _0)$ satisfying the equation
\begin{equation}
u(x) = v(x) - \int _{\Omega _{0}}G(x;y) \; d \mu (y).
\end{equation}
almost everywhere in $\Omega _0$. If $\mu (\Omega)<\infty$ or $u$
has a $c-$harmonic majorant in $\Omega$, then $v\in \verb"H"(c, \Omega)$
and the representation (2.5) holds almost everywhere in $ \Omega $. \\

$9^{0}$ For every subfunction $u \in \verb"SbH"(c, \Omega)$ there
exist a sequence of domains $\Omega _n \subset \Omega$ and a
sequence of continuous subfunctions $u_n \in \verb"SbH"(c, \Omega
_n) \cap C(\overline {\Omega _n })$ such that
\[\lim _{n \rightarrow \infty } u_n (x) = u(x) \mbox{at each point $\; x \in
\Omega$} ,\]
\[\int _{\Omega _n} |u(x)-u_n (x)| dx \rightarrow 0, \; n
\rightarrow \infty ,\] and
\[\lim _{n \rightarrow \infty}{\cal M}^{G}(u_{n},x,r) =
{\cal M}^{G}(u,x,r).\] The domain $\Omega _n$ consists of all
interior points of $\Omega$ distant from the boundary $\partial
\Omega$ at least for $1/n$ and
${\cal M}^G (u,x,r)$ is given by (2.4).\\
\end{theor}

\noindent \textbf{Proof} of Part $9^{0}$ repeats the proof of a
similar statement in Littman \cite{Litt2} and uses Corollary
A.3 (Appendix A). \hfill \textbf{$\diamondsuit $}\\

\newpage

\section{Phragm\'{e}n-Lindel\"{o}f theorem for \\
subfunctions in an $n-$dimensional cone}

In this section we extend the following theorem of Phragm\'{e}n
and Lindel\"{o}f (\cite{PhL}; see, for example, \cite[p. 37-38]{L2}. \\

\emph{ Let $f(z)$ be a holomorphic function in an angle
$D=\{\alpha < \arg z < \beta \}$ such that
\begin{equation}
\liminf _{r\rightarrow \infty } r^{- \rho } \ln ^+ \max \{|f(z)|
\; \biggl| \; z\in B(r) \cap D \} = 0,
\end{equation}
where $0<\rho < \infty$. If $\limsup_{D\ni \zeta \rightarrow
z}|f(\zeta)| \leq A = const $ for all $z \in \partial D$ and
\begin{equation}
\rho (\beta - \alpha ) < \pi ,
\end{equation}
then $|f(z)| \leq A $ everywhere in the angle $D$.} \\

All similar statements have since been called Phragm\'{e}n-
Lindel\"{o}f theorems. The theorem has enjoyed numerous
generalizations and applications - see, for example, a beautiful
account of the earlier history of the subject in \cite[p.
15-22]{Gar}, works  \cite{BerCoh, Dahl, Ess, Evg, Fra, JiLa, Lax,
OlRad, PrWein}, and references therein. We state here only a
generalization due to Nevanlinna \cite{Nev} (Cf. Ahlfors \cite[Theorem 5]{Ahl1}):\\

\textit{The conclusion remains valid if} (3.1) \emph{ is replaced
by a weaker condition
\begin{equation}
\liminf_{r\rightarrow \infty}r^{-\rho} \int_{\alpha}^{\beta} \ln
^+ |f(re^{i \theta}| \sin \pi \frac{\theta -\alpha }{\beta -
\alpha } d \theta = 0 \end{equation}
and in} (3.2) $\rho (\beta - \alpha ) \leq \pi$.\\

Let $D$ be a domain on the unit sphere $S=S^{n-1} \subset
\mathbf{R}^n , \; n\geq 2$. We always assume that the boundary
$\partial D$ with respect to $S$ is not a polar set in the
classical potential theory meaning. Let
\[K^D = \left\{x= (r,\theta ) \in \mathbf{R}^n \biggl| \; 0<r< \infty,
\; \theta \in D \right\} \] be a cone generated by the domain $D$.
Truncated cones are denoted by $K^D _r = K^D \bigcap B(0,r)$, $K^D
_{r, \infty } = K^D \setminus \overline{K^D _r}$, $K^D _{r, R} =
K^D _R \setminus \overline{K^D _r}, \; 0<r<R< \infty
$. \\

Consider an equation $\Delta u(x) - q(|x|)u(x)=0$ with a radial
potential $q(r)$ in the cone $K^D$. Adjoining the zero boundary
conditions to the equation and separating variables in spherical
coordinates $r,\theta _1 , \ldots , \theta _{n-1}$, we get for the
radial component of the solution the well-known ordinary
differential equation (Cf. (1.2))
\begin{equation}
y''(r)+(n-1)r^{-1}y'(r)-\left(\lambda r^{-2}+q(r)\right)y(r)= 0,
\end{equation}
where $\lambda$ is a separation constant. Appendix C contains the
properties of solutions of (3.4) used in this paper. For the
angular component of the equation $\Delta u(x) - q(|x|)u(x)=0$ we
deduce a known eigenvalue problem for the Laplace-Beltrami
operator $\Delta ^*$,
\begin{equation} \left\{
\begin{array}{ll}
\Delta ^* \varphi (\theta )+\lambda \varphi (\theta )=0  &
    \mbox{if  $ \theta \in D$} \vspace{.1cm} \\
\varphi (\theta ) =0 & \mbox{if  $\theta \in \partial D \setminus E$} \vspace{.1cm} \\
\varphi \in L^2(D), & \mbox{}
\end{array}
\right.
\end{equation}
where $E$ is an exceptional polar set of irregular points on
$\partial D$. Necessary properties of the eigenvalues and
eigenfunctions of this problem are considered in detail in Appendix B. \\

Let $\lambda_0$ be the smallest eigenvalue of the problem (3.5)
and $\varphi_0$ be the corresponding eigenfunction; it is known
that $\lambda_0 >0$, for $\partial D$ is not a polar set. The
positive root $\mu ^+$ of the quadratic equation $\mu (\mu+n-2)=
\lambda _0$ is called \emph{ the characteristic constant}
of the domain $D\subset S$. \\

The following result of Deny-Lelong \cite{DeL} extends the
Phragm\'{e}n-Lindel\"{o}f theorem\footnote{Deny and Lelong stated
their result in a slightly different way.} onto the
subharmonic functions in $n-$dimensional cones. \\

\emph{ Let $u$ be a subharmonic function in a cone $K^D$, where
$D$ is a regular domain and
\[M(r, u)=\sup \left\{ u(x) \biggl| \; x \in K^D _r\right\}.\] If
\begin{equation}
\liminf _{r \rightarrow \infty } r^{- \mu ^+ } M(r, u) \leq 0
\end{equation}
and $\sup u(x) \bigl| _{x \in \partial K^D } \leq A $, then }
\[u(x) \leq A, \; \forall x \in K^D .\]

Now let $u \in \verb"SbH"(c, K^D)$ with $c \in {\cal C}(K^D)$. Let
$q(r)=q_c(r), \; r = |x|,$ be a nonnegative radial minorant of the
potential $c(x)$, that is, a measurable function $q(r), \; 0<r<
\infty$, such that
\[0 \leq q_c(r) \leq \inf \left\{ c(x) \biggl| \; x\in K^D _r , \; |x|=R \right\}. \]
Obviously, such functions exist, for example, $q(r)=0$ is one of
them. Due to (2.2), $q \in L_{loc}(0,\infty)$. An increasing
solution of the equation (3.4) with $\lambda = \lambda _0$ and
$q=q_c$ is denoted by $V(r)=V_q (r)$. We normalize it by $V(1)=1$.
Now we state the main result of this section. \\

\begin{theor}\hspace{-.09in}\textbf{.} Let $D\subset S$ be a domain such
that its boundary $\partial D$ with respect to $S$ is not a polar
set, $c \in {\cal C}(K^D)$, and $u \in \verb"SbH"(c, K^D)$. If
\begin{equation}
\limsup  _{K^D \ni x \rightarrow x_0 } u(x) \leq A = const, \; \forall x_0 \in \partial K^D,
\end{equation}
and
\begin{equation}
\liminf_{r\rightarrow \infty } \frac{1}{V_q (r)} \int _D u^+ (r,
\theta ) \varphi _0(\theta ) d\sigma (\theta ) = 0,
\end{equation}
then\footnote{The integral in (3.8) is called the Nevanlinna norm
of $u^+$ \cite{Hein}.}
\begin{equation}
u(x) < A^+ , \; \forall x \in K^D.
\end{equation}

The conclusion is the best possible in the sense that, as an
example below shows, $A^+$ in (3.9) cannot be replaced by $A$,
that is, if in (3.7) $A<0$, then, provided (3.8), one can claim in
(3.9) only $u(x)<0, \forall x \in K^D$. Moreover, the example of
$L_q$-harmonic function $u(r, \theta ) = V_q (r) \varphi _0
(\theta )$ shows that the condition (3.8) is the best possible in
the class of all the potentials $c$ dominating $q$, that is, such
that $c(x) \geq q(|x|)$.
\end{theor}

\noindent \textbf{Remark 3.1.} When $c=0$, Theorem 3.1 reduces to
the Deny-Lelong theorem. Since we impose no regularity conditions
on $\partial D$, the Deny-Lelong theorem also holds without any
such a restriction. \hfill  \textbf{$\diamondsuit $ }\\

\noindent \textbf{Remark 3.2.} If $c=0$, $n=2$, and $D$ is the
angle $\{ \alpha < \arg z < \beta \}$, then $\rho = \pi /(\beta -
\alpha )$, $\varphi _0 (\theta ) =\sin \left(\pi \frac{\theta -
\alpha
}{\beta - \alpha}\right)$, and (3.8) becomes (3.3). \hfill  \textbf{$\diamondsuit $ }\\

\noindent \textbf{Remark 3.3.} It should be noted \cite{DeL} that
if $\partial D$ is not a polar set, then for all truncated cones
$K^{\partial D} _R, \; 0<R<\infty$, their lateral surfaces
$K^{\partial D} _R$ are not polar sets in $\mathbf{R}^n$.
Moreover, due to (2.2), if $\Omega _1$ is a bounded proper
sub-domain of $\Omega $, then its regularity with respect to the
Dirichlet problem for the Laplacian is equivalent to its
regularity for the Dirichlet problem for the operator $L_c$ and
for $L_q$ with any nonnegative radial $q, \; 0 \leq q \leq c$
\cite{Maed}. \hfill  \textbf{$\diamondsuit $ }\\

\noindent\textbf{Proof of Theorem 3.1.} We shall prove the theorem
in several steps. We always assume that $A=0$, since the general
case can be immediately reduced to this one. First, suppose
additionally that all truncated cones $K^D _R, \; \forall R>0$,
are regular, that is, the Dirichlet problem for $\Delta , \; L_c$
and so for $L_q$ is solvable in each of these cones.

We use in the following an idea of Ahlfors \cite{Ahl1}. Due to
Theorem 2.1, Part $2^o$, $\; u^+ (x) \in \verb"SbH"(q, K^D)$. So,
by solving the Dirichlet problem for $L_q$ in $K^D _R, \; R>0$,
with the boundary data
\[u(x)= \left\{\begin{array}{ll}
           0 & \mbox{ for $x\in \partial K^D _R \; $ and $\; |x|<R$} \vspace{.2cm}\\
           u^+ (R,\theta ) & \mbox{ for $x=(R, \theta ), \; \theta \in D$},
                \end{array} \right. \]
we construct the lowest $q-$harmonic majorant $\widetilde{u}$ of
the function $u^+ (x)$; this majorant satisfies the same
boundary conditions \vspace{.1cm} \\
$\widetilde{u}(x) \bigl|_{\partial K^D \cap B(R)}=0 \; $ and $\;
\widetilde{u}(r, \theta ) \bigl|_{r=R}= u^+ (R, \theta ) \; $
for $\; \theta \in D$. \\

Since $L_q\widetilde{u}(r,\theta)=0$ in the sense of
distributions, we can write
\begin{equation}
\int _{K^D _R} \widetilde{u}(x) L_q \chi (x)dx=0,
\end{equation}
where $\chi$ is any finitary function in $K^D _R$. Obviously, we
can use $\chi(r,\theta)=\alpha(r)\varphi _0(\theta)$ with $\alpha
\in C^{\infty }_0 (0,R)$. Therefore, from (3.10) we get
\[0 = \int _{K^D _R} \widetilde{u}(x) L_q \chi (x)dx =\]
\[\int _{K^D _R} \widetilde{u}(r,\theta) \{\varphi_0 (\theta)
\alpha ''(r) + \varphi_0 (\theta) \frac{n-1}{r} \alpha'(r)- q(r)
\varphi_0 (\theta) \alpha (r) + \frac{1}{r^2} \alpha(r) \Delta^*
\varphi_0 (\theta) \}dx,\] where $\Delta ^*$ is the
Laplace-Beltrami operator. By the definition of $\varphi _0$,
\[\int _{K^D _R} \widetilde{u}(x) r^{-2}\alpha(r) \Delta ^* \varphi _0(\theta)dx =
- \lambda _0 \int _{K^D _R}\widetilde{u}(x) r^{-2}\alpha(r)
\varphi _0(\theta )dx.\] Combining these equations together, we
deduce
\[\int _{K^D _R} \varphi_0(\theta) \widetilde{u}(r,\theta) \left\{\alpha ''(r) +
(n-1)r^{-1}\alpha '(r) - \left(\lambda_0 r^{-2} + q(r) \right)
\alpha (r)\right\}dx = \]
\[\int_0 ^R y_0(r) \left\{\alpha ''(r) +
(n-1)r^{-1}\alpha '(r) - \left(\lambda_0 r^{-2} + q(r) \right)
\alpha (r) \right\} r^{n-1} dr=0,\] where
\[y_0(r)=\int _D\widetilde{u}(r,\theta)\varphi _0(\theta)d\sigma (\theta).\]
Since $\alpha$ is an arbitrary finitary function in $(0,R)$, we
conclude that $y_0 (r)$ satisfies the equation (3.4) with $\lambda
=\lambda _0>0$ and $q = q_c (r)$ for $0 < r < R$. It also
satisfies the boundary condition $y_0 (0) = 0$. So that, $y_0 (r)
= b V_q (r), \; b = const$. Substituting here $r = R$, we get $b=
y_0(R)/V_q(R)$ and
\[y_0 (r)=\frac{y_0(R)}{V_q(R)} V_q(r), \; 0< r\leq
R<\infty .\] Letting here $R\rightarrow\infty$ over a subsequence,
while $r$ is fixed and considering (3.8), we get $y_0 (r) =
0, \; \forall r \geq 0$. Thus, by the definition of
$\widetilde{u}$ we have
\[0\leq \int_D u^+ (r,\theta)\varphi_0(\theta )d\sigma(\theta)\leq
\int_D \widetilde{u}(r,\theta) \varphi _0 (\theta )d\sigma (\theta
) = y_0 (r) = 0, \; \forall r > 0.\] Since $\varphi _0(\theta)>0$
for $\theta \in D$, it follows that $u^+ (\theta)=0$ almost
everywhere in $K^D$. Therefore, $u(x) \leq 0$ almost everywhere in
the cone, and if $u$ is continuous, then $u(x) \leq 0$ everywhere
in the cone. If $u$ is merely an upper-semicontinuous subfunction,
the inequality $u(x) \leq 0$ follows from the mean-value property
(2.3). Let us remind that we suppose $A = 0$. Now the strict
inequality (3.9) follows from the
maximum principle. \\

Next we remove the regularity restriction on the boundary of the
cone $K^D$. Fix an $R>0$ and denote by $Q_R = \partial K^D \cap
B(R)$ the lateral surface of the truncated cone $K^D _R$. If
$n\geq 3$, we denote by $\Gamma (x;y)= \theta _n |x-y|^{2-n}$ a
fundamental solution of the operator $L_q$ in the entire space
$\mathbf{R}^n$, $x=(r,\theta) \in K^D$, $y=(\rho,\psi)$, and
assume $y \not\in \overline{K^D}$. If $n=2$, let $\Gamma (x;y)=
2\pi \ln \frac{1}{|x-y|} + h(x,y)$ be the Green function of $L_q$
in any domain strictly containing the sector $K^D$.\\

Since $\limsup_{x\rightarrow x_0} u^+(x)=0, \; x_0 \in Q_R$, by
the finite covering theorem for any $\eta >0$ we can find a
neighborhood of the lateral surface $Q_R$, such that
$u^+(r,\theta)<\eta, \; 0 \leq r \leq R$. Therefore, for each
$\varepsilon >0$ we can choose a neighborhood of $Q_R$, where the
function
\[v_{\varepsilon}(r,\theta) \stackrel{def}{=} u^+(r,\theta) - \varepsilon \Gamma
(r,\theta;\rho,\psi) < 0 \] and so that $v^+ _{\varepsilon}(x) =0$
in this neighborhood.

Now we exhaust the domain $D$ from within by an expanding sequence
of domains $\{D_l\}_{l=1} ^{\infty}$ with smooth boundaries. All
truncated cones $K^{D_l} _R , \; l= 1,2, \ldots $, are regular -
see Remark 3.3. We denote $Q^l _R = K^{\partial D_l} \cap B(R)$.
For $l \geq l(\varepsilon)$, $Q^l _R$ is situated in the domain
where $v^+ _{\varepsilon}(x) =0$. Moreover, after continuing by
zero outside $K^D$, $v^+ _{\varepsilon}$ becomes a $q-$subfunction
everywhere in $\mathbf{R}^n$.

On the other hand, for the upper semi-continuous function $v^+
_{\varepsilon}$ there exists a sequence of continuous functions $v
_{\varepsilon, \delta}$ monotonically decreasing to $v^+
_{\varepsilon}$ in $K^D$. If $l>l(\delta)$, $v_{\varepsilon,
\delta} (x)=0$ for $x\not\in K^{D_l}$. By solving the Dirichlet
problem for the operator $L_q$ in the truncated cone $K^{D_l} _R$
with boundary data $v_{\varepsilon, \delta}(x)$ and letting
$\delta \searrow 0$ we obtain, in the limit, a solution
$\widetilde{v}_{\varepsilon, l}$ of the Dirichlet problem for the
equation $\Delta \widetilde{v}_{\varepsilon, l}
-q(r)\widetilde{v}_{\varepsilon, l} = 0$, which vanishes at the
lateral boundary $Q^l _R$ and coincides with $v^+
_{\varepsilon}(R, \theta)$ on the cap $S(R)\cap K^{D_l} _R$.

It is clear that
\[\int _{D_l} \left(\Delta \widetilde{v}_{\varepsilon, l}
- q(r)\widetilde{v}_{\varepsilon, l} \right) \varphi _0 ^{(l)}
(\theta) d\sigma (\theta) =0, \] that is,
\[\int _{D_l} \left(\frac{\partial ^2 \widetilde{v}_{\varepsilon,
l}}{\partial r^2} + \frac{n-1}{r}\frac{\partial
\widetilde{v}_{\varepsilon, l}}{\partial r} - \frac{\lambda ^{(l)}
_0}{r^2} \widetilde{v}_{\varepsilon, l} - q(r)
\widetilde{v}_{\varepsilon, l} \right) \varphi _0 ^{(l)} (\theta)
d\sigma (\theta) =0, \] or, after extending $\varphi _0
^{(l)}(\theta) = 0, \theta \not\in D_l$,
\begin{equation}
\int _{D} \left(\frac{\partial ^2 \widetilde{v}_{\varepsilon,
l}}{\partial r^2} + \frac{n-1}{r}\frac{\partial
\widetilde{v}_{\varepsilon, l}}{\partial r} - \frac{\lambda ^{(l)}
_0}{r^2} \widetilde{v}_{\varepsilon, l} - q(r)
\widetilde{v}_{\varepsilon, l}\right) \varphi _0 ^{(l)} (\theta)
d\sigma (\theta) =0.
\end{equation}

We also extend $\widetilde{v}_{\varepsilon, l}=0$ outside $K^{D_l}
_R$, which makes it a $q-$subfunction in the entire
$\mathbf{R}^n$. When $l$ increases, the function
$\widetilde{v}_{\varepsilon, l}(x)$ also monotonically increases
and it does not exceed a solution of the Dirichlet problem for
$L_q$ in the ball $B(R)$ with boundary data $u^+ (R,\theta)$ for
$\theta \in D$ and $0$ for $\theta \in S(R) \setminus D$.
Therefore, there exists a finite limit
\[\widetilde{v}_{\varepsilon}(x) = \lim _{l \rightarrow \infty}
\widetilde{v}_{\varepsilon, l}(x),\] which is a $q-$harmonic
function in $K^D _R$. This limit $\widetilde{v}_{\varepsilon}(x)$
is equal to $v^+ _{\varepsilon}(R,\theta)$ as $\theta \in D$ and
vanishes outside a polar set at $Q_R$.

Letting $l \rightarrow \infty$ in the integral
\begin{equation}
y_l (r) = \int _{D} \widetilde{v}_{\varepsilon, l}(r, \theta)
\varphi _0 ^{(l)} (\theta) d\sigma (\theta),
\end{equation}
we get
\begin{equation}
y_l (r) \rightarrow y(r) = \int _{D} \widetilde{v}_{\varepsilon
}(r, \theta) \varphi _0 (\theta) d\sigma (\theta).
\end{equation}\\

We prove next that the function $y(r)$, defined in (3.13),
satisfies the equation
\begin{equation}
y'' + (n-1)r^{-1}y' - \left(\lambda _0 r^{-2} + q(r)\right)y = 0
\end{equation}
and a boundary condition $y(0^+)=0$. The latter follows
immediately from the equation $\widetilde{v}_{\varepsilon}(0,
\theta)=0$. To prove the former, we note that by virtue of (3.11)
and (3.12), the function $y_l (r)$ satisfies the
equation
\begin{equation}
y'' _l + (n-1)r^{-1}y' _l - \left(\lambda ^{(l)} _0 r^{-2} +
q(r)\right)y_l = 0, \; y_l (0^+) = 0.
\end{equation}
It is shown in Lemma C.1, Part $3^o$, that if $y_l (0^+) = 0, \;
y_l (1) = 1$, and $\lambda ^{(l)} _0 \searrow \lambda _0 \geq 0$,
then a solution $y_l$ of the equation (3.15) converges,
monotonically decreasing, to a solution\footnote{If $D$ is "almost
the entire sphere" $S$ and $q(r)=0$, then $V(r)= const$. That is
why we assume that $\partial D$ is not a polar set.} $V(r)$ of
(3.14). Thus,
\[\int _{D} \widetilde{v}_{\varepsilon}(r, \theta) \varphi _0 (\theta) d\sigma
(\theta)= bV(r).\]

To find the constant factor $b$, we can repeat the same reasoning
as in the smooth case and use the equation $\int _{D}
\widetilde{v}_{\varepsilon}(R, \theta) \varphi _0 (\theta) d\sigma
(\theta)= bV(R)$. Noticing also that the $q-$subfunction $v^+
_{\varepsilon}(r, \theta) \leq \widetilde{v}_{\varepsilon}(r,
\theta)$, we get an inequality
\[\frac{1}{V(r)} \int _{D} v^+ _{\varepsilon}(r, \theta) \varphi _0
(\theta) d\sigma (\theta) \leq \frac{1}{V(R)} \int _{D} v^+
_{\varepsilon}(R, \theta) \varphi _0 (\theta) d\sigma (\theta)
\leq \]
\[\frac{1}{V(R)}\int _{D} u^+ (R, \theta) \varphi _0 (\theta) d\sigma (\theta),
\; 0\leq r \leq R.\]
The condition (3.8) implies that
\[\int _{D} v^+ _{\varepsilon}(r, \theta) \varphi _0 (\theta) d\sigma
(\theta)= 0, \; 0 \leq r < \infty .\] Since $\varphi _0$ is
positive and $v^+ _{\varepsilon}$ is a nonnegative subfunction, it
follows that $v^+ _{\varepsilon}(r, \theta) =0$, that is,
\[u^+(r,\theta) \leq \varepsilon \; \Gamma (r,\theta;\rho,\psi),
\; \forall\varepsilon >0, \] and so $u(r,\theta) \leq 0$,
which completes the proof of (3.9) since we assume $A=0$.\\

We finally show that the conclusion of Theorem 3.1 fails if $A$ is
negative and we replace $A^+$ by $A$. It is enough here to
consider regular cones. First, we construct a $c-$harmonic
function $v_{-1}$ in $K^D$ such that $-1 < v_{-1} (x) < 0$ for $x
\in K^D$ and $v_{-1} \bigl| _{\partial K^D } = -1 $. Indeed, since
a cone $K^D _R $ is regular, the boundary value problem $L_c v(x)
= 0$ in $K^D _R$ with the boundary condition $v(x) = -1$ at
$\partial K^D _R$ has a solution $v^R _{-1}(x)$. The latter
function cannot have a nonnegative maximum inside the domain and
must attain its smallest negative value at the boundary,
therefore, $-1 < v^R _{-1} (x) < 0$ in $K^D _R$. If $R_1 < R_2$,
then $v^{R_1} _{-1} (x) = -1$ on the "cap" $D_{R_1} = K^D \cap
S(0, R_1)$, however, $v^{R_2} _{-1} (x) > -1$ there. Therefore,
$-1 < v^{R_1} _{-1} (x) < v^{R_2} _{-1} (x) < 0, \; x \in K^D
_{R_1}$, and there exists the limit $\lim _{R \rightarrow \infty }
v^R _{-1} (x) = v _{-1} (x)$, which is a $c-$harmonic function. In
addition, due to the maximum
principle, $-1<v_{-1}(x)<0, \; x\in K^D$.\\

It is obvious that $v_{-1}(x)=-1$ as $x\in \partial K^D$ and
$v_{-1}(x)$ is not a constant in $K^D$. Define a function
\begin{equation}
v_0(x)=\frac{A}{\alpha +1} \left(\alpha - v_{-1}(x) \right),
\end{equation}
where $\alpha = \sup _{x \in K^D}v_{-1}(x)$, so that $-1 < \alpha
\leq 0$. If $A<0$, then $v_0(x) \bigl| _{\partial K^D}=A<0$, but
$\sup _{x \in K^D} v_0 (x)=0$; therefore, we cannot replace $A^+$
in the statement  by $A$. The proof of Theorem 3.1 is complete.
\hfill \textbf{$\diamondsuit $}\\

\noindent \textbf{Remark 3.4. } We show now that actually in
(3.16) $\alpha = \sup _{x \in K^D} v_{-1}(x)=0$. To simplify
calculations, we assume that $\inf_{x \in K^D} c(x)=c_0 >0$. Fix a
ball $B(x_0,\rho) \subset K^D$. Since $v_{-1}(x)>-1$ on the sphere
$S(x_0,\rho)$, we have
\[v_{-1}(x)= \int _{S(x_0,\rho)} v_{-1}(y) \frac{\partial G(x;y)}{\partial
n(y)}d\sigma (y) > \]
\[-\int_{S(x_0,\rho)} \frac{\partial
G(x;y)}{\partial n(y)}d\sigma (y) \geq -\int_{S(x_0,\rho)}
\frac{\partial G_{c_0}(x;y)}{\partial n(y)}d\sigma (y),\] where
$G_{c_0}$ is the Green function of the operator $L_{c_0}$ in
$B(x_0,\rho)$ and the second inequality follows from Theorem A.2.
The latter integral represents a solution $u(x)$ of the Dirichlet
problem for $L_{c_0}$ in $B(x_0,\rho)$ such that $u \bigl|
_{S(x_0,\rho)} =-1$. The solution can be written explicitly,
$u(x)= - \frac{\Psi_{c_0}(r)}{\Psi_{c_0}(\rho)}$, where
$r=|x-x_0|$ and $\Psi _{c_0}$ is the (normalized) modified Bessel
functions $I_{n/2-1}$, namely,
\[\Psi_{c_0}(r) = b r^{1 - n/2} I_{n/2 - 1} (r \sqrt {c_0}), \;
\Psi_{c_0}(0)=1.\] Thus, $u(x_0)=-\frac{1}{\Psi _{c_0}(\rho)}$ and
we get $v_{-1}(x_0) \geq -\frac{1}{\Psi _{c_0}(\rho)}$. For a
fixed $x_0 \in D$ we can let $x\rightarrow \infty$ and therefore
$r\rightarrow \infty$, so that $\rho = \rho (r)$ can be made
arbitrarily large. Since $I_{n/2-1}(\rho \sqrt{c_0} ) \rightarrow
\infty $ as $\rho \rightarrow \infty$, we have $v_{-1}(r,\theta)
\rightarrow 0$ when $r \rightarrow \infty$ and $\theta \in Q
\subset D$, where $Q$ is a compact set. Therefore, $\sup _{x \in
K^D} v_{-1}(x)=0$. \hfill  \textbf{$\diamondsuit $}\\

\begin{cor}\hspace{-.09in}\textbf{.} Obviously, the conclusion
of Theorem 3.1 holds true if (3.8) is replaced by
\[\liminf _{r\rightarrow \infty } V^{-1} _q (r) M(r,
u^+ ) = 0.\]
\end{cor}
\hfill  \textbf{$\diamondsuit$}

In particular, this implies
\begin{cor}\hspace{-.09in}\textbf{.} If $u(x) \geq 0$ and
$\liminf _{r\rightarrow \infty } V^{-1} _q (r) M(r,u^+ ) = 0$,
then $u \equiv 0$ in $K^D$. \hfill  \textbf{$\diamondsuit $}
\end{cor}

\begin{cor}\hspace{-.09in}\textbf{.} Under the conditions of Theorem 3.1,
either $u(x) \leq A^+$ in $K^D$, or else
\[||u(r,\cdot )||_{L^2(D)} \geq b V_q (r).\]
\end{cor}
\noindent \textbf{Proof.} If $u(x_0)>A^+$ at a point $x_0 \in
K^D$, then
\[V_q(r) \leq b \int _D u(r,\theta)\varphi_0
(\theta) d\sigma(\theta), \; 0\leq r \leq \infty ,\] and it
suffices to apply the Bunyakovskii-Cauchy-Schwartz inequality to
the latter integral. \hfill  \textbf{$\diamondsuit $} \\

Applying now the principle of $c-$harmonic majorant (Theorem 2.1,
Part $5^0$) we get the following statement.
\begin{cor}\hspace{-.09in}\textbf{.} If $u \in \verb"SbH"(c, B(R))$,
$c(x) \geq c_0 = const \geq 0$ and $u(x) \biggl |_{S(R)} \leq -1$,
then
\[\max \left\{ u(x) \biggl | \; x \in B(R) \right\} \leq -
\frac{1}{\Psi_{c_0}(R)},\] where $\Psi_{c_0}$ is the same is in
Remark 3.4; this maximum is attained only on the solution $v(x) =
- \frac{\Psi_{c_0}(|x|)}{\Psi_{c_0}(R)}$ of the corresponding
Dirichlet problem. \hfill  \textbf{$\diamondsuit $} \\
\end{cor}
In these statements, the minorants $q(r)$ of the potential $c(x)$
play an essential role. The next statement involves its radial
majorant, that is, a measurable radial function
$Q(r)\geq \sup \left\{c(x) \biggl| \; x\in K^D_r \right\}$. \\
\begin{prop}\hspace{-.09in}\textbf{.} If a potential $c$ has a radial
majorant $Q(r)\in L^p _{loc} [0, \infty)$ with the same $p$ as in
(2.2), then for every real $A$ there exists a function $u_A (x)
\in \verb"SbH"(c, K^D)$ such that $u_A (x) \bigl| _{x \in
\partial K^D} = A, \; \; \sup _{x \in K^D} u_A (x) > A^+ $, and
\[\liminf _{r \rightarrow \infty }V^{-1} _q(r)M(r,u^+ _A)> 0.\]
\end{prop}

\noindent\textbf{Proof.} Let $V_Q(r)$ be an upper (increasing)
solution of the equation (3.4) with $\lambda = \lambda _0$ and
$Q(r)$ instead of $q(r)$. It follows from Lemma C.1 that \\
$\liminf _{r \rightarrow \infty } V_Q (r) / V_q (r) > 0$. If $u_0
(x) = V_Q (r) \varphi _0 (\theta )$,
then $u_0 \in \verb"SbH"(c, K^D)$ and we can set \\
\[ u_A (x) = \left\{ \begin{array}{ll}
        A + u_0 (x) & \mbox{if $ A \geq 0 $ } \\
        v_0 (x) + u_0(x) & \mbox{if $A < 0,$ }
                                \end{array}
                  \right. \]
where $v_0 (x)$ is given by (3.16). \hfill  \textbf{$\diamondsuit $ }\\

Theorem 3.1 can be stated as follows: \\

If a $c-$subfunction $u$ is positive at even one point in a cone
and satisfies in the cone the boundary condition (3.7), then it
grows in the cone at least as $V(r), \; r \rightarrow \infty $. If
we have more information on the asymptotic behavior of the
potential $c$, the conclusion of Theorem 3.1 can be made more
precise. The borderline case occurs if the potential behaves as
the inverse square, $|x|^{-2}$. Indeed, if the potential is a weak
perturbation of the Laplacian in the sense that
\[\limsup _{|x|\rightarrow \infty } |x|^2 c(x) < \infty ,\]
then many results are similar to those in the harmonic case
$c(x)=0$. On the other hand, if
\[\limsup _{|x|\rightarrow \infty } |x|^2 c(x) = \infty ,\]
the results may be essentially different. To distinguish these two
possibilities, we introduce the following two classes of
potentials; these definitions include also some mild regularity
conditions. \\

Let $0 \leq q(r) \in L_{loc} (0, \infty )$. Denote $s(r)=r^2
q(r)$. If there exists the finite limit
\[\lim _{r \rightarrow \infty } s(r) =k \in [0, \; \infty )\]
and
\[\int ^{\infty } \left(s(t) - k \right)^2 \frac{dt}{t} < \infty
,\] then we write $q \in (\mathbf{A})$ and refer to this case as
the case \textbf{(A)}. If a potential $c \in {\cal C}(\Omega)$ has
a minorant $q \in \mathbf{(A)}$, we write $c \in {\cal C}(\Omega ,
\mathbf{A})$. \\

On the other hand, if $s(r)=r^2 q(r)$ is monotonically increasing,
moreover,
\[\lim _{r \rightarrow \infty } s(r) = \infty ,\]
\[r^{-2} q^{-1/2} (r) \in L(1, \; \infty ),\]
and in addition $q$ has the second derivative such that
\[\int ^{\infty } \biggl| 4 \widetilde{q}(r) \widetilde{q}''(r) -
5 (\widetilde{q} '(r))^2  \biggl| (\widetilde{q}(r))^{-5/2} dr <
\infty , \] where $\widetilde{q}(r) = q(r) +\left(\frac{n^2 - 4n
+3}{4} + \lambda _0 \right)r^{-2}$, then we write $q \in
(\mathbf{B})$ and refer to this case as the case $\mathbf{(B)}$.
If the potential $c \in {\cal C}(\Omega)$ has a minorant $q \in
\mathbf{(B)}$, we write $c \in {\cal C}(\Omega , \mathbf{B})$. \\

With these notations, Lemmas C.2 and C.3 imply immediately

\begin{theor}\hspace{-.09in}\textbf{.} Under the conditions of Theorem 3.1,
if $c \in {\cal C}(\Omega , \mathbf{A})$, then the lowest possible
growth rate of an unbounded $c$-subfunction $u$ in $K^D$,
satisfying the boundary condition (3.7), is given by the
right-hand side of (C.9), that is,
\[V(r)= b \; r^{(2-n+ \chi _{k}) \; /\; 2} \exp \left\{\frac{1}{\chi _{k}}
\int ^r _1 (t^2 q(t) - k) \frac{dt}{t} \right\}
(1+\bar{\bar{o}}(1)), \; r \rightarrow \infty ,\] where $\chi _k
^2 = (n-2)^2 + 4(\lambda _0 + k)$.

In the case \textbf{(B)} this lowest rate does not depend on a
domain $D$ and is given by the \emph{JWKB-}asymptotic formula
(C.11)
\[V(r)= b \; r^{(1-n)/2}q^{-1/4}(r) \exp \left\{\int
^{r}_{1}q^{1/2}(t)dt \right\} \; (1+\overline{\overline{o}}(1)),
\; r \rightarrow + \infty .\]
\end{theor} \hfill  \textbf{$\diamondsuit $}\\

\noindent \textbf{Remark 3.5.} As was mentioned above, under the
condition (3.7) the lowest possible growth of an unbounded
$c$-subfunction in the cone $K^D$ in the case \textbf{(B)} does
not depend on the domain $D$ generating the cone. It turns out
that this growth coincides with the minimal possible growth of a
$c$-harmonic function in the entire space. Namely, the following
simple analog of the classical Liouville
theorem\footnote{Generalizations of the Liouville theorem on
$c-$harmonic functions with respect to a non-constant radial
potential $q(r)$ are studied in \cite{KhePOTA2}; see also
\cite[Theorem 3]{VeM2}.} on bounded analytic functions is valid
for the $c$-harmonic functions with respect to a constant
potential $c(x)= c_0$. Below $I_k$ is
the modified Bessel function of first kind. \hfill \textbf{$\diamondsuit $}\\
\begin{prop}\hspace{-.09in}\textbf{.} Let $u(x) \in H(c_0, \mathbf{R}^n )$,
where $c(x) \equiv c_0 = const \geq 0$. If
\[\liminf_{r\rightarrow
\infty}\frac{M(r,|u|)}{I_0(r\sqrt{c_0})}=0,\]
then $u(x) \equiv 0 $. \\
\end{prop}
\noindent\textbf{Proof.} We expand a continuous function $u$ in a
series against the spherical harmonics on the unit sphere,
\[u(r, \theta ) = \sum ^{\infty } _{k=0} I_k (r \sqrt{c_0})
\left( \sum _{l=0} ^{l_k }a_{l,k} P_{l,k} (\theta ) \right). \]
Due to the orthogonality of spherical harmonics,
\[a_{l,k} I_k (r \sqrt{c_0}) = \int _{S(0, 1)} u(r, \theta )
P_{l,k} (\theta ) d\sigma (\theta ), \] therefore, $|a_{l,k}| \leq
M(r, |u|) / I_k (r \sqrt{c_0})$. This and the assumption imply
that all coefficients $a_{l,k} = 0$. \hfill \textbf{$\diamondsuit$}\\

For example, let us consider a minorant $q(x)=c_0 =\inf _{x\in
K^D}c(x)$. When $c_0 =0$, we get an extension of the Deny-Lelong
theorem with the same limit growth (3.6). Now, if $c_0 >0$, then
the case \textbf{(B)} takes place. An increasing solution of (3.4)
in the cone $K^D$ is

\[V_1(r)=b_1r^{(2-n)/2}I_{\alpha}(r\sqrt{c_0}),\; b_1= const, \;
\alpha =(1/2)\sqrt{n(n-2)+4\lambda _0}.\] In the entire space
$\mathbf{R}^n $, an increasing solution is

\[V_2(r)=b_2r^{(2-n)/2}I_{\beta}(r\sqrt{c_0}),\; b_2= const,\; \beta =
(1/2)\sqrt{n(n-2)}.\] Due to the well-known asymptotic formulas
for the Bessel functions, if the appropriate constants $b_1$ and
$b_2$ are chosen, then $c_0-$harmonic functions
\[u_1 (r, \theta ) = b_1 r^{(2-n)/2} I_{\alpha } (r \sqrt{c_0})
\varphi _0 (\theta ) \in H(c_0, K^D),\] in the cone $K^D$ and
\[u_2 (r, \theta ) = b_2 r^{(2-n)/2} I_{\beta } (r \sqrt{c_0})
\in H(c_0, \mathbf{R}^n)\] in the entire space satisfy the
relation
\[M(r, u_1 ) \sim  u_2 (r) \sim r^{-1/2} \exp \left\{r \sqrt{c_0 } \right\},
\; r \rightarrow \infty . \] In the class of potentials with $\inf
c(x) = c_0 > 0$ this conclusion is the best possible and improves
the results of \cite{HerZh} related to the operator $L_c$. \\

\noindent \textbf{Remark 3.6.} Theorem 3.1 deals with the
convexity of the ratio
\begin{equation}
\frac{1}{V(r)} \int _D u^+ (r, \theta ) \varphi _0 (\theta )
d\sigma (\theta ).
\end{equation}
Related questions are considered by Maz'ya and Verzhbinski
\cite{VeM1}. Their Theorem 6.1 asserts the inequality
\begin{equation}
\frac{1}{y(r)} \left\{\int _{D_r} \widehat{u}^2 (r, \theta )
d\sigma (\theta ) \right\}^{1/2} \leq \frac{1}{y(R)} \left\{\int
_{D_R} \widehat{u}^2 (R, \theta ) d\sigma (\theta )\right\}^{1/2},
\end{equation}
where $r>R\searrow 0$ and $y(r)$, similarly to $V(r)$, satisfies
an ordinary differential equation close to (3.4). Clearly, (3.17)
and (3.18) involve different integral norms, but it is not very
essential. A more important distinction, from our point of view,
is that the function $\widehat{u}$ in (3.18) vanishes at the
boundary in some neighborhood of $0$ and is a subfunction outside
this neighborhood, while the function $u$ in our Theorem 3.1,
after the Kelvin transformation $u(r, \theta ) \longmapsto r^{2-n}
u(\frac{1}{r}, \theta )$, is a subfunction in a neighborhood of
$0$ and a solution elsewhere. \hfill \textbf{$\diamondsuit$}\\

\noindent \textbf{Remark 3.7.} Similar results (Cf. \cite{DeL})
can be proved on the behavior of subfunctions at the vertex of a
cone, as $r \rightarrow 0$, with an obvious substitution of the
decreasing solution $W(r)$ of (3.4) instead of $V(r)$ in all the
statements above. \hfill \textbf{$\diamondsuit $}\\

The classical Phragm\'{e}n-Lindel\"{o}f theorem can be made more
precise. First, let us notice that if $u \in \verb"SbH"(c, K^D)$
and $\sigma \geq 0$, then the difference $u(r, \theta ) - \sigma
V_q (r) \varphi _0 (\theta )$ is a $c$-subfunction. Applying
Theorem 3.1 to this difference, we deduce the following simple
result.
\begin{cor}\hspace{-.09in}\textbf{.} If a subfunction $u$ satisfies (3.7)
with $A=0$ and either
\[\liminf_{r\rightarrow \infty}V^{-1} _q (r) \int_D \left(u(r,\theta ) -
\sigma V_q (r) \varphi_0(\theta) \right)^+ \varphi _0 (\theta )
d\sigma (\theta ) = 0 \] or
\[\liminf_{r\rightarrow \infty} \int_D \left(\frac{u(r,\theta )}{V_q (r) \varphi_0(\theta)}
- \sigma \right) ^+ \varphi ^2 _0 (\theta ) d\sigma (\theta ) =
0\] with a $\sigma \geq 0$, then $u(r, \theta ) \leq \sigma V_q
(r)\varphi _0(\theta)$ everywhere in $K^D$. \hfill \textbf{$\diamondsuit$}\\
\end{cor}

The following result is known for any real $A$ \cite[p. 50]{L1}. \\

\emph{ Let $f(z)$ be an analytic function in an angle $\{| \arg z
| < \alpha \}$ and continuous up to the boundary. If $\ln |f(r
e^{\pm i \alpha })| \leq A, \; \forall r>0$, and
\[ \liminf _{r \rightarrow \infty } r^{-\rho } \ln ^+ M(r, f) = \sigma ,\]
where $2 \rho \alpha \leq \pi $, then
\[ \ln |f(r e^{i \theta })| \leq A + \sigma r^{\rho } \cos
\left(\frac{\pi \theta }{2 \alpha } \right), \; 0 \leq r < \infty
, \; |\theta | \leq \alpha .\]}

Let us notice that $\cos \left(\frac{\pi \theta }{2 \alpha }
\right)$ is a positive eigenfunction of the two-dimensional
Laplace-Beltrami operator in $|\theta |< \alpha $, vanishing at
$\pm \alpha$, but it is not normalized in $L^2$, rather
$\max_{|\theta | \leq \alpha} \cos \left(\frac{\pi \theta }{2
\alpha } \right) =1$. In the following analogous result for
subharmonic functions $\widehat{\varphi}$ stands for a positive
eigenfunction of the Laplace-Beltrami operator $\Delta^*$ in a
domain $D\subset S$, vanishing at $\partial D$, and normalized in
$L^{\infty}(D)$ as $\max _{\theta \in D} \widehat{\varphi}=1$. By
$\widehat{\theta} \in D$ we denote a point where the maximum is
attained at, that is, $\widehat{\varphi}(\widehat{\theta})=1$.

\begin{prop}\hspace{-.09in}\textbf{.} Let $u$ be a subharmonic function
in a cone $K^D$, such that $\partial D$ is not a polar set. If $u
\bigl|_{\partial K^D} \leq A = const $ and
\[\liminf_{r\rightarrow \infty}r^{-\mu^+}M(r,u^+)=\sigma <\infty ,\]
then $u(r, \theta ) \leq A + \sigma r^{\mu ^+} \widehat{\varphi}
(\theta )$ everywhere in $K^D$.
\end{prop}

\noindent\textbf{Proof.} Again, without loss of generality we set
$A=0$. Denote
\[D_{\delta} = \left\{\theta \in D \biggl| \; |\theta - \widehat{\theta}|
< \delta , \; \delta>0 \right\}\] and consider a subharmonic
function
\[u_{\epsilon } (r, \theta ) = u^+ (r, \theta ) - (\sigma + \epsilon
) r^{\mu ^+ } \widehat{\varphi}(\theta ), \; \epsilon > 0.\] Given
an $\varepsilon >0$, we can find a small $\delta = \delta
(\varepsilon )>0$ and a sequence $\{r_l\}^{\infty}_{l=1}$ such
that $u_{\varepsilon}(r_l, \theta ) \leq 0$ for all $l$ and $x \in
\partial \left(K_{r_l}^{D\setminus D_{\delta}}\right)$. The
maximum principle infers that $u_{\varepsilon} \leq b(\varepsilon
) =const$ at $\partial K^{ D_{\delta}}$. When a domain decreases,
the eigenvalues and the corresponding characteristic constant
increase, therefore, with the obvious notations, $\lambda _0
(D\setminus D_{\delta}) > \lambda _0 (D)$ and $\mu ^+ (D\setminus
D_{\delta}) > \mu ^+ (D) \equiv \mu ^+$. Thus we can apply the
Deny-Lelong theorem (see Remark 3.1) to the function $u_{\epsilon
}$ in the cone $K^{D\setminus D_{\delta}}, \; \delta = \delta
(\varepsilon )$, and conclude that $u_{\varepsilon} \leq
b_{\varepsilon}$ in $K^{D\setminus D_{\delta}}$.

On the other hand, decreasing $\delta $ if necessary, we can make
true the inequality $\mu ^+ (D _{\delta }) > \mu ^+ (D\setminus D_
{\delta })$, thus $\liminf _{r \rightarrow \infty } r^{-\mu ^+
(D_{\delta})}M(r, u^+ _{\epsilon }) = 0$. Applying the Deny-Lelong
theorem again, we get the inequality $u_{\epsilon}(x) \leq
b_{\epsilon }$ in the cone $K^{D_{\delta}}$, therefore this
inequality is valid in the entire cone $K^D$. Applying the
Deny-Lelong theorem once more, we obtain the inequality
$u_{\epsilon }(x) \leq 0$ in $K^D$, that is, $u(x) \leq (\sigma +
\epsilon )r^{- \mu ^+} \widehat{\varphi}(\theta )$, and it
suffices now to let $\epsilon \rightarrow 0$. \hfill \textbf{$\diamondsuit$}\\

However, for subfunctions such a conclusion in general does not
hold. To see that, it is enough to consider a constant potential
$c_0 = const > 0$ in the right half-plane $K= \left\{|\theta |<
\frac{\pi}{2}\right\}$ in $R^2$. The eigenvalues are $\lambda _k =
k^2, \; k =0,1,2,...$, and the eigenfunctions are
\[\varphi_k (\theta) = \left\{\begin{array}{ll}
        \sin(k\theta), & \mbox{if $k=2l,$ } \vspace{.2cm} \\
        \cos(k\theta), & \mbox{if $k=2l+1.$ }
                                \end{array}
                  \right. \]

The functions
\[u_k (r, \theta )= \sigma I_k \left(r \sqrt{c_0}\right) \varphi _k (\theta ),\]
where $I_k$ are the modified Bessel functions, are $c_0 -$harmonic
functions in $K$, vanish at its boundary, and due to the
well-known asymptotic formula for $I_k$, satisfy the following
equation for all $k = 0,1,2,...$,
\[\lim  _{r \rightarrow \infty } \max _{|\theta | \leq \pi /2}
\frac{u^+ _k (r, \theta )}{I_k (r\sqrt{c_0})} = \sigma .\]

Yet, the inequality $u_k (z) \leq u_0 (z)$ can not be valid
everywhere in the angle $K$ for all $k \geq 2$, since it reduces
to inequalities $\sin (2k \theta ) \leq \cos \theta $ or $\cos
((2k+1) \theta ) \leq \cos \theta $, which obviously fail for all
large enough $k$ uniformly in $|\theta | \leq \pi /2$.
\hfill \textbf{$\diamondsuit $}\\

Nonetheless, in the case \textbf{(A)} subfunctions again exhibit
the same behavior as subharmonic functions.

\begin{theor}\hspace{-.09in}\textbf{.} Let $u \in \verb"SbH"(c, K^D)$ with
$c \in {\cal C}(K^D, \mathbf{A})$ and $\partial D$ is not a polar
set. If the conditions (3.7) and $\liminf _{r \rightarrow \infty }
V^{-1}(r) M(r, u^+ ) = \sigma , \; 0 \leq \sigma < \infty$, are
valid, then everywhere in $K^D$
\[u(r, \theta ) \leq A^+ + \sigma V(r) \widehat{\varphi}(\theta ).\]
\end{theor}

\noindent \textbf{Proof.} It suffices to repeat the proof of
Proposition 3.3 with $r^{\mu ^+ }$ replaced by $V(r)$ and refer to
Theorem 3.1 instead of the Deny-Lelong theorem. \hfill \textbf{$\diamondsuit$}\\

In general case we claim the following. \\

\begin{theor}\hspace{-.09in}\textbf{.} Under the conditions of Theorem 3.1,
let the quantity
\begin{equation}
M_0 (r) \equiv M_0 (r, u^+) = \sup _{\theta \in D} \frac{u^+(r,
\theta )}{\widehat{\varphi}(\theta )} < \infty
\end{equation}
be upper bounded for $0 \leq r <\infty$. If
\[\liminf _{r \rightarrow \infty } M_0 (r) / V_q (r)= \widehat {\sigma }, \; 0
\leq \widehat{\sigma } < \infty ,\] then
\begin{equation}
u(r,\theta)\leq \widehat{\sigma }V_q(r)\widehat{\varphi}(\theta).
\end{equation}
\end{theor}

\noindent\textbf{Proof.} Let us notice that (3.19) implies
\[u^+ (r, \theta ) = 0, \; \theta \in \partial D \setminus E, \; 0 < r <
\infty .\] A $c-$subfunction $\widehat{u}(r, \theta ) = u(r,
\theta ) - \widehat{\sigma }V_q (r) \widehat{\varphi}(\theta )$
satisfies the equation
\[\liminf_{r\rightarrow \infty}M_0(r,\widehat{u})/V_q (r)=0.\]
Since $\max \widehat{\varphi}(\theta )=1$, one has $M(r) \leq M_0
(r)$, which together with (3.19) implies $\liminf _{r \rightarrow
\infty } M(r, \widehat{u}^+) / V_q (r) = 0$. Therefore,
$\widehat{u}^+ (r, \theta )=0$ for $\theta \in \partial D
\setminus E$, that is, $\widehat{u}(r, \theta ) \leq 0$ at the
boundary of $K^D$. Now (3.20) follows from Theorem 3.1. \hfill \textbf{$\diamondsuit$}\\

\noindent\textbf{Remark 3.8.} We mention one more distinction
between the cases \textbf{(A)} and \textbf{(B)}. It is known that
there exists a unique, up to a constant factor, positive harmonic
function in a cone, vanishing at the boundary - namely, the
so-called reduced function $r^{\mu ^+ }\varphi_0 (\theta)$.  For
the operator $L_c$ this is generally not true. In the case
\textbf{(B)}, that is, when the potential grows fast enough, there
exist infinitely many linearly independent $c-$harmonic functions,
$c= const$, which are positive in a cone and vanish at the
boundary. For example, we can set $u_\nu (r, \theta)=V(r)
\varphi_0(\theta) + b_{\nu } V_{\nu}(r) \varphi_{\nu} (\theta)$
with different indices $\nu $, if the constants $|b_{\nu }|$ are
sufficiently small. That follows from the asymptotic formulas of
Appendix C and the estimates of the
eigenfunctions - see, for example, (4.2) below. \\

However, in the case \textbf{(A)} the situation is again similar
to the harmonic case. For the second order divergence form
elliptic operators without lower order terms this property was
established in \cite{LaN}.

\begin{prop}\hspace{-.09in}\textbf{.} Let $D \subset S$ be an arbitrary
domain whose boundary is not a polar set and a radial potential
$q=q(r) \in (\mathbf{A})$.Then, up to a constant factor, there
exists only one positive solution of the equation $L_q =0$ in the
cone $K^D$ vanishing at the boundary $\partial (K^D)$ outside a
polar set. To be precise, the solution is
\[u(r, \theta ) = b V_0 (r) \varphi _0 (\theta ), \; b = const
\geq 0.\]
\end{prop}

\noindent\textbf{Proof.} Let $\{\lambda _{\nu} , \varphi _{\nu}\}$
be the eigenvalues and the corresponding eigenfunctions of the
Laplace-Beltrami operator in $D$, where eigenvalues repeat
according to their multiplicities. First, let $\partial D \in
C^2$. Then everywhere on $\partial (K^D)$ $\frac{\partial \varphi
_0 (\theta )}{\partial n(\theta )} \geq 0$, and $|\varphi _{\nu}
(\theta )|+|\nabla \varphi_{\nu}(\theta)|\leq N_{\nu},\; \nu \geq
1$. Therefore, $\frac{\varphi_{\nu}(\theta)}{\varphi_0 (\theta )}
\leq M_{\nu}<\infty$ for $\theta \in \overline{D}$ and $\nu \geq 1$. \\

Now we repeat the argument used in the proof of Proposition 3.2.
If $L_q u=0$ in $K^D$ and $u \bigl| _{\partial K^D} =0$, then for
each $r>0$
\[u(r, \theta ) = \sum ^{\infty } _{\nu =0} Y_{\nu}(r)\varphi
_{\nu}(\theta), \] where
\begin{equation}
Y_{\nu}(r)=\int_Du(r,\theta)\varphi_{\nu}(\theta )d\sigma (\theta)
\end{equation}
and the series converges in $L_2 (K^D)$.

Clearly, $0 = \int _D \left(\Delta u - q(r)u\right) \varphi _{\nu}
(\theta ) d\sigma (\theta )$ and (3.21) implies the equation
\[Y_{\nu}''(r)+\frac{n-1}{r}Y_{\nu}'(r)-q(r)Y_{\nu}(r)+\int_D
\Delta ^*u(r,\theta)\varphi_{\nu}(\theta ) d\sigma (\theta ) =0.\]
Using the self-adjointness of the Laplace-Beltrami operator
$\Delta ^*$, we conclude that $Y_{\nu}$ satisfies the equation
(3.4) with $\lambda = \lambda _{\nu}$, that is,
\begin{equation}
Y_{\nu}(r) = a_{\nu}V_{\nu}(r) + b_{\nu}W_{\nu}(r).
\end{equation}
We know that as $r\searrow 0 \; $, $W_{\nu}(r)\nearrow \infty$ and
$V_{\nu}(r)\searrow 0$. However, due to (3.21), $Y_{\nu}$ must be
bounded as $r\rightarrow 0$, so that in (3.22) all $b_{\nu}=0, \;
k=0,1,2, \ldots$, and $Y_{\nu}(r) = a_{\nu}V_{\nu}(r)$.

Now, let $u(r, \theta )$ be positive in $K^D$. Then
\[\biggl|\int_Du(r,\theta)\varphi_{\nu}(\theta)d\sigma(\theta)\biggl|
\leq \int_Du(r,\theta)\varphi_0(\theta)\frac{|\varphi
_{\nu}(\theta )|}{\varphi _0 (\theta )}d\sigma (\theta ), \] that
is, $|a_{\nu}V_{\nu}(r)|\leq M_{\nu}|a_0 | V_0 (r)$. However,
Lemma C.2 implies that \\
$V_{\nu}(r)/V_0(r)\rightarrow \infty $ as $r \rightarrow \infty $,
whence $a_{\nu}=0$ for $\nu \geq 1$ and
\[u(r, \theta ) = a_0 V_0 (r) \varphi _0 (\theta ), \; a_0 \geq
0.\]

If $\partial D$ is not smooth, we again exhaust $K^D$ from within
with a sequence of smooth cones $K^{D_j}, \; j = 1,2, \ldots $,
and use the uniform convergence of $\lambda _{\nu}^{[j]}, \;
\varphi _{\nu}^{[j]}, \; V_{\nu}^{[j]}$ and $W_{\nu} ^{[j]}$ on
compact sets (Lemma C.1, Part $3^o$). \hfill  \textbf{$\diamondsuit$}\\

\newpage

\section{Bilinear series and estimates of the Green function of the
operator $L_q$ with a radial potential $q$ in an $n$-dimensional
cone}

In this section we derive a bilinear series representation and
some estimates of the Green function $G(x;y)$ of the operator
$L_q$ with a radial potential $q = q(r), \; r=|x|$, in a smooth
$n-$dimensional cone $K^D$. As before, the eigenvalues and
eigenfunctions of the boundary value problem (3.5) are denoted by
$\lambda _{\nu }$ and $\varphi _{\nu }, \; \nu = 0, 1,2,...$, with
each eigenvalue repeated according to its multiplicity. We often
omit the subscript $\nu =0$. The eigenfunctions are normalized in
$L^2 (D)$, $\| \varphi _{\nu } \| _{L^2 (D)} = 1$, and are
continued on for all $x \in K^D$ as $\varphi _{\nu } (x) = \varphi
_{\nu } (x/|x|)$. The following properties of the eigenvalues and
eigenfunctions are known (see, for example, \cite[p. 335]{Mikh1}
and \cite{IlShis}):
\begin{equation}
b_1 (\nu +1)^{\frac{2}{n-1}} < \lambda _{\nu } < b_2 (\nu
+1)^{\frac{2}{n-1}},
\end{equation}
\begin{equation}
\max _{\theta \in \overline{D}} |\varphi _{\nu } (\theta )| = \tau
_{\nu } < b_3 \lambda _{\nu } ^{n/4}, \; \; \max _{\theta \in
\overline{D}} |\varphi ^{'} _{\nu } (\theta )| = \tau ^{'} _{\nu }
< b_4 \lambda _{\nu } ^{\frac{n+2}{4}},
\end{equation}
where $b_j$ are constants and $\varphi ' _{\nu}$ stands for any
first order partial derivative of $\varphi _{\nu}$.

Let $\mu ^+ _{\nu}\geq 0$ and $\mu ^- _{\nu} < 0$ be the roots of
the quadratic equation
\[\mu (\mu + n-2)=\lambda _{\nu} \]
and
\[\chi _{\nu}=\mu ^+_{\nu}-\mu ^- _{\nu}=\sqrt{(n-2)^2
+4\lambda _{\nu}}.\] Also, let $V_{\nu }(r)$ and $W_{\nu}(r)$
stand, respectively, for the increasing and decreasing, as
$r\rightarrow \infty$, solutions of the equation (3.4) with
$\lambda =\lambda _{\nu }$, normalized by
$V_{\nu}(1)=W_{\nu}(1)=1$, and
\[\chi ^{'} _{\nu}=w (W_{\nu}(r), V_{\nu}(r))\biggl|_{r=1}\]
be their Wronskian at $r=1$. Since by Part 1$^o$ of Lemma C.1,
$V'_{\nu}(1) \geq \mu ^+ _{\nu}$ and $W'_{\nu}(1) \leq \mu ^-
_{\nu}$, we have
\begin{equation}
\chi ^{'}_{\nu}= V'_{\nu}(1)W_{\nu}(1)-W'_{\nu}(1)V_{\nu}(1)=V'
_{\nu}(1)-W'_{\nu}(1)\geq \mu _{\nu}^+ -\mu _{\nu}^- =\chi _{\nu}.
\end{equation}\\

It is readily verified that the Green function $Q_{\nu }$ of the
boundary value problem
\[- y''(r)-(n-1)r^{-1}y'(r)+\left(\lambda _{\nu}r^{-2}+q(r)\right)y(r) =
 \gamma (r), \; 0<r< \infty,\]
\begin{equation}
0\leq q(r)\in L_{loc}(0, \infty ), \; r^{n-1}\gamma (r) \in L^2(0,
\infty ), \; y(0)=y(\infty ) = 0,
\end{equation}
is given by
\begin{equation}
Q_{\nu }(r, \rho ) = \left\{ \begin{array}{ll}
        \rho ^{n-1} (\chi ^{'} _{\nu })^{-1} V_{\nu }(r) W_{\nu }(\rho )
& \mbox{if $\; 0 < r \leq \rho < \infty $ } \vspace{.5cm} \\
        \rho ^{n-1} (\chi ^{'} _{\nu })^{-1} V_{\nu }(\rho ) W_{\nu }(r)
& \mbox{if $\; 0 < \rho \leq r < \infty $.}
                                \end{array}
                  \right.
\end{equation}

Indeed, the corresponding homogeneous equation
\[-y''(r)-(n-1)r^{-1}y'(r)+\left(\lambda _{\nu}r^{-2}+q(r)\right)y(r) =0\]
has the general solution $y(r) = b_1 V_{\nu}(r) + b_2 W_{\nu}(r)$
and the boundary conditions (4.4) imply $b_1 =b_2 =0$. It should
be also noted that if $n=2$ and
$D= \{|\theta |<\alpha \}, \; 0<\alpha <\pi$, then $\lambda _{\nu} > 1/4$.\\

The following statement on the comparison of the Green functions
is a one-dimensional analog of Theorem A.2.

\begin{lem}\hspace{-.09in}\textbf{.} If $p_j (r) \in L_{loc} (0, \infty ), \; j = 1, 2$,
and $0 \leq p_1 (r) \leq p_2(r)$ for \\
$0< r < \infty $, then the Green functions $Q_j$ of the operators
\[L_j (y) = - y'' (r) - (n-1) r^{-1} y' (r) + p_j (r) y(r) \]
satisfy the inequality
\begin{equation}
Q_2 (r, \rho ) \leq Q_1 (r, \rho ), \; 0  <r, \rho < \infty .
\end{equation}
\end{lem}

\noindent \textbf{Proof.} A function $u(r)=Q_1(r,\rho
)-Q_2(r,\rho)$ is continuous in $(0, \infty)$ together with $u'$
for each fixed $\rho >0$. If $r \neq \rho$, then $L_1u=(p_2
-p_1)Q_2$, whence
\[u(r)= \int _0 ^{\infty } Q_1 (r,t) \left(p_2(t)-p_1(t)\right)
Q_2 (t,\rho) dt \geq 0, \]
maybe in the sense of distributions, and the conclusion follows.
\hfill \textbf{$\diamondsuit$}\\

\begin{cor}\hspace{-.09in}\textbf{.} Setting here $p_1 (r)=\lambda _0 r^{-2}+q(r)$
and $p_2 (r)= \lambda _{\nu}r^{-2}+q(r), \; \nu \geq 1$, where
$0\leq q(r) \in L_{loc}(0, \infty )$, we get from (4.5) and (4.6)
the inequalities
\[\begin{array}{c}
\frac{1}{\chi ^{'} _{\nu}}V_{\nu}(r)W_{\nu}(\rho)\leq
\frac{1}{\chi ^{'} _0}V_0 (r)W_0 (\rho) \; \mbox{ if } \; r \leq \rho , \vspace{.5cm} \\
\frac{1}{\chi ^{'} _{\nu }}V_{\nu }(\rho ) W_{\nu }(r) \leq
\frac{1}{\chi ^{'} _0 }V_0 (\rho ) W_0 (r) \; \mbox{ if } \; \rho
\leq r.
\end{array} \] \hfill \textbf{$\diamondsuit$}
\end{cor}

\begin{cor}\hspace{-.09in}\textbf{.} Choosing in Lemma 4.1
$p_1 (r)=\lambda _{\nu}r^{-2}$ and $p_2 (r) =\lambda
_{\nu}r^{-2}+q(r)$, where $0\leq q(r)\in L_{loc}(0, \infty )$, we
deduce the inequalities
\begin{equation}
\begin{array}{c}
\frac{1}{\chi ^{'} _{\nu }} V_{\nu }(r) W_{\nu }(\rho ) \leq
\frac{1}{\chi _{\nu }} r^{\mu ^+ _{\nu }} \rho ^{\mu ^- _{\nu }}
\; \mbox{ if } \; r \leq \rho , \vspace{.5cm} \\
\frac{1}{\chi ^{'} _{\nu }}V_{\nu }(\rho ) W_{\nu }(r) \leq
\frac{1}{\chi _{\nu }} \rho ^{\mu ^+ _{\nu }} r^{\mu ^- _{\nu }}
\; \mbox{ if } \; \rho \leq r.
\end{array}
\end{equation}
\hfill \textbf{$\diamondsuit$}
\end{cor}

\begin{cor}\hspace{-.09in}\textbf{.} If $0 < \gamma \leq 1 \leq \delta < \infty $
and $\delta r \leq \gamma \rho $, then \\
\begin{equation}
\frac{1}{\chi ^{'} _{\nu }} V_{\nu }(r) W_{\nu }(\rho ) \leq
\frac{1}{\chi _{\nu }} \gamma ^{- \mu ^- _{\nu }} \delta ^{- \mu
^+ _{\nu }} V(\delta r) W(\gamma \rho ).
\end{equation}
\end{cor}

\noindent \textbf{Proof.} It suffices to note (Lemma C.1, Part
$1^0$) that $\frac{d}{dr} \left(r^{- \mu ^+ _{\nu }}V_{\nu
}(r)\right) \geq 0$, thus $r^{- \mu ^+ _{\nu }}V_{\nu }(r) \leq
(\delta r)^{- \mu ^+ _{\nu }}V_{\nu }(\delta r)$ and so $V_{\nu
}(r) \leq \delta ^{- \mu ^+ _{\nu }}V_{\nu }(\delta r)$. On the
other hand, $\frac{d}{d \rho } \left(\rho ^{- \mu ^- _{\nu
}}W_{\nu }(\rho ) \right) \leq 0$, which implies $W_{\nu }(\rho )
\leq \gamma ^{- \mu ^- _{\nu }}W_{\nu }(\gamma \rho )$.
Multiplying these two inequalities and using (4.3), we obtain
(4.8). \hfill \textbf{$\diamondsuit$}\\

The following bilinear series representation of Green's function
$g(x;y)$ of the Laplacian in a cone $K^D$ is well-known (see
Lelong-Ferrand \cite{LF}):
\begin{equation}
g(x;y)=\left\{\begin{array}{ll}
        \beta \sum _{\nu =0}^{\infty}|x|^{\mu ^+ _{\nu}}
        |y|^{\mu ^- _{\nu}}\frac{\varphi _{\nu}(\theta)\varphi _{\nu}(\psi)}{\chi _{\nu}}
            & \mbox{if $|x|<|y|$} \vspace{.5cm} \\
        \beta \sum _{\nu =0}^{\infty}|y|^{\mu ^+ _{\nu}}
        |x|^{\mu ^- _{\nu}}\frac{\varphi _{\nu}(\theta)\varphi _{\nu}(\psi)}{\chi _{\nu}}
            & \mbox{if $|y|<|x|$;}
                \end{array}
                \right.
\end{equation}
here $x=(r,\theta ), \; y=(\rho , \psi )$, and $\beta = \int _D
d\sigma (\theta)$. Due to (4.1)-(4.2), these series converge
absolutely and uniformly if, respectively, either $|x|<\gamma |y|$
or $|y|<\gamma |x|$ with $0<\gamma <1$. Now we derive an analog of
(4.9) for the Green function $G_q (x;y)$ of the operator $L_q$
with a radial potential $q(r)$. We suppress the subscript $q$ from
now on in this section.

\begin{theor}\hspace{-.09in}\textbf{.} For all $x = (r, \theta )$ and
$y = (\rho , \psi ), \; 0 < r, \rho < \infty, \; (\theta , \psi )
\in D \times D$, one has
\begin{equation}
G(r, \theta ; \rho , \psi ) = \left\{ \begin{array}{ll}
        \sum ^{\infty } _{\nu = 0} \frac{1}{\chi '_{\nu }} V_{\nu }(r) W_{\nu }(\rho )
         \varphi _{\nu }(\theta ) \varphi _{\nu }(\psi )
            & \mbox{if  $\; 0 < r < \rho < \infty $} \vspace{.5cm} \\
        \sum ^{\infty } _{\nu = 0} \frac{1}{\chi '_{\nu }} V_{\nu }(\rho ) W_{\nu }(r)
        \varphi _{\nu }(\theta ) \varphi _{\nu }(\psi )
            & \mbox{if  $\; 0 < \rho < r < \infty $.} \\
                                \end{array}
                  \right.
\end{equation}
On any compact set in $K^D$ these series converge absolutely and
uniformly if $r\leq \gamma \rho$ or respectively, if $\rho \leq
\gamma r, \; 0< \gamma <1$.
\end{theor}

\noindent \textbf{Remark 4.1.} The factor $\beta$ in (4.9) appears
since the normalization of $g$ at its singularity used in \cite{LF},
differs from that in Theorem A.1, Part $3^o$. \hfill \textbf{$\diamondsuit$}\\

\noindent \textbf{Proof of Theorem 4.1.} The convergence of the
series (4.10) follows immediately from (4.7), (4.3) and the
convergence of the series (4.9). Thus, we have to establish that
the right-hand side of (4.10) is equal to the Green function $G$.
The latter function represents a solution of the equation $L_q u =
f$, where $f\in L^2 (K^D)$, by the integral $u(x)=\int
_{K^D}G(x;y)f(y)dy \equiv (Gf)(x)$. Since $G(x;y)\in L^2 _{loc}
(K^D \times K^D)$ by Theorem A.1, $Gf$ is a completely continuous
integral operator and the system of eigenfunctions $\{\varphi
_{\nu }\}$ is complete in $L^2 (D)$. Thus, we can expand a
function $f(\cdot \; , \psi ) \in L^2 (K^D)$ in a series
\[f(\rho ,\psi)=\sum ^{\infty}_{\nu =0}\gamma _{\nu}
(\rho)\varphi _{\nu}(\psi),\] where
\[\gamma _{\nu}(\rho)= \int _D f(\rho , \psi)\varphi
_{\nu}(\psi) d\psi .\] The series converges in $L^2 (D)$, that is,
$\sum ^{\infty } _{\nu =0} |\gamma _{\nu } (\rho )|^2 < \infty $.
Considering this, to solve the equation $L_q u=f$ we first
look for "elementary" solutions of equations
\begin{equation}
L_q u(r, \theta ) = \gamma _{\nu } (r) \varphi _{\nu } (\theta ),
\; \nu =0,1, \ldots ,
\end{equation}
where now $\gamma _{\nu}$ is a finitary function in $(0, \infty)$.
Substituting $u(x)= y(r) \varphi _{\nu } (\theta )$ in (4.11), we
get for $y(r)$ a boundary-value problem
\begin{equation}
\left\{ \begin{array}{l}
        -y''(r)-(n-1)r^{-1}y'(r)+\left(\lambda _{\nu } r^{-2}
        + q(r)\right) y(r) = \gamma _{\nu }(r), \vspace{.2cm} \\
        y(+0) = 0, \; y \in L^2 _{loc} (0, \infty ),\\
             \end{array}
                  \right.
\end{equation}
whose Green's function is given by (4.5). Therefore, the unique
solution of (4.12) in $L^2 (0, \infty )$ is
\begin{equation}
y_{\nu }(r) = \int ^{\infty } _0 Q_{\nu } (r, \rho ) \gamma _{\nu } (\rho) d\rho .
\end{equation}

Now, for a solution of the equation $L_q u(x) = f(x) = \sum
^{\infty}_{\nu =0}\gamma _{\nu} (r)\varphi _{\nu}(\theta)$ we
obtain a representation $u(r, \theta ) = \sum ^{\infty } _{\nu =
0} y_{\nu } (r) \varphi _{\nu } (\theta )$, where the functions $
y_{\nu } (r)$ are given by (4.13). By virtue of the Parseval
formula we have the equations
\[\int _D |f(r, \theta )|^2 d\sigma (\theta ) = \sum ^{\infty } _{\nu = 0}
|\gamma _{\nu } (r)|^2 \]
for almost all $r>0$, and
\[\int _D \int ^{\infty } _0 |f(r, \theta )|^2 r^{n-1} dr d\sigma (\theta )
= \sum ^{\infty } _{\nu = 0} \int ^{\infty } _0 |\gamma _{\nu }
(r)|^2 r^{n-1} dr.\]

It follows that the series
\[u(r,\theta )=\sum ^{\infty} _{\nu =0}y_{\nu}(r)
\varphi _{\nu}(\theta)=\int _{K^D}G(r,\theta ;y)f(y)dy \]
also converges in $L^2 (K^D)$. \\

Consider now a series
\[G_1 (r, \theta ; \rho , \psi ) = \sum ^{\infty } _{\nu = 0}
\rho ^{1-n} Q_{\nu } (r, \rho ) \varphi _{\nu }(\theta ) \varphi
_{\nu }(\psi ).\] Due to (4.7), it is dominated by the similar
series (4.9) converging absolutely and uniformly for $r \leq
\gamma \rho$ or $\rho \leq \gamma r, \; \gamma < 1$. Therefore,
the series $G_1 (x;y)$ also converges absolutely and uniformly.
So, in the sense of $L^2$
\[\int _{K^D}G_1 (r,\theta ; y)f(y)dy=\sum ^{\infty} _{\nu = 0}
y_{\nu } (r) \varphi _{\nu}(\theta),\] implying that the
difference $G(x;y)-G_1 (x;y)$ is orthogonal to all finitary
functions $f(y)$, and finally, $G(x;y)=G_1 (x;y)$
almost everywhere in $K^D \times K^D$. \hfill \textbf{$\diamondsuit$}\\

Using (4.10), we now estimate\footnote{Cf. \cite{VeM1}} $G(x;y)$.

\begin{lem}\hspace{-.09in}\textbf{.} If
$0 <\gamma \leq 1 \leq \delta < \infty , \; \gamma \neq \delta$,
and $\delta r \leq \gamma \rho $, then \\
\begin{equation}
G(r, \theta ; \rho , \psi ) \leq V(\delta r) W(\gamma \rho
)(\gamma \delta )^{(n-2)/2} \sum ^{\infty } _{\nu = 0}
\left(\gamma / \delta \right)^{\chi _{\nu }/2} \frac{|\varphi
_{\nu } (\theta ) \varphi _{\nu } (\psi )|}{\chi _{\nu }}.
\end{equation}
\end{lem}
\noindent \textbf{Proof.} It suffices to apply the inequality
(4.8) to (4.10) and use the above mentioned remark on convergence.
\hfill  \textbf{$\diamondsuit $ }\\

From now on in this section we assume $\partial D$ to be smooth,
so that $\varphi _{\nu} \in C^2(\overline{D})$. To bound the
series in (4.14) we use the following theorem by Oleinik
\cite{Ole}, which we state only in the case of the operator $L_c$:\\

\emph{ Let $u(P)$ be a non-constant solution of the equation $L_c
u = f \geq 0$ in a domain $D$; thus, in particular, $u$ is a
$c-$superfunction. Let at each point $P\in \partial D$ there exist
a tangent plane and a ball $B\subset D$, such that $P\in
\partial B$. If $u$ is continuous in $\overline{D}$ and at a point
$P_1 \in \partial D$ takes on the lowest negative value, then
\[\liminf _{P_2 \rightarrow P_1} \frac{u(P_2) - u(P_1)}{|P_2 - P_1
|} >0 ,\] where $P_2 \in D$ and varies along a direction $h$
making an acute angle with the inward normal to the boundary
$\partial D$ at the point $P_1$.}\\

This statement implies that $\frac{\partial \varphi _0 (\theta
)}{\partial n} \bigl| _{\theta \in \partial D} \geq b>0$, so that
\[\lim _{\theta \rightarrow \theta _0 \in \partial D} \frac{\varphi _{\nu }(\theta )}
{\varphi _0 (\theta )} = \frac{\frac{\partial \varphi _{\nu
}(\theta _0)}{\partial n}}{\frac{\partial \varphi _0 (\theta _0
)}{\partial n}} \leq b(D) \nu ^{\frac{n+1}{2(n-1)}}.\] From here
and (4.14) we obtain

\begin{cor}\hspace{-.09in}\textbf{.} If $\partial D \in C^2$ and
$0 < \gamma \leq 1 \leq \delta < \infty , \;
\gamma \neq \delta$, then \\
\[G(r, \theta ; \rho , \psi ) \leq \left\{ \begin{array}{ll}
        b(D, \gamma, \delta) V(\delta r) W(\gamma \rho )
                \varphi _0 (\theta ) \varphi _0(\psi )
& \mbox{if  $\; 0 < \delta r \leq \gamma \rho < \infty $ } \vspace{.5cm} \\
        b(D, \gamma, \delta) V(\delta \rho ) W(\gamma r)
                \varphi _0(\theta ) \varphi _0(\psi )
& \mbox{if  $\; 0 < \delta \rho \leq \gamma r < \infty .$ }
                                \end{array}
                  \right. \]
\end{cor} \hfill \textbf{$\diamondsuit$}

These inequalities can be made more precise. For that we need two
more lemmas, including the positive homogeneity property of the
harmonic Green function $g$ in a cone $K^D \subset \textbf{R}^n$.

\begin{lem}\hspace{-.09in}\textbf{.} For any $k>0$
\[g(kx;ky)= k^{2-n} g(x;y).\]
\end{lem}
\noindent \textbf{Proof.} If $n\geq 3$, then $g(x;y)=\theta _n
|x-y|^{2-n} + $ a harmonic function in both $x$ and $y$. Thus, the
function $u(x,y)=g(kx;ky) - k^{2-n}g(x;y)$ is harmonic in $K^D$ in
both variables, vanishes at $\partial K^D$, and uniformly tends to
zero when one variable approaches infinity while another is fixed.
By the maximum principle \cite[p. 47]{HK}, $u(x,y)=0$. If $n=2$,
the proof is similar. \hfill \textbf{$\diamondsuit$}

\begin{lem}\hspace{-.09in}\textbf{.} If $\partial D \in C^2$, $r_0 >0$
and $\delta >1$ are fixed, then for $0<r \leq r_0 < \rho _0 =
\delta r_0 < \rho < \infty$ and for $\theta \in D$
\[g(r, \theta ; \rho ,\psi )\leq b r^{\mu ^+} \rho
^{\mu ^- } \varphi _0 (\theta ) \varphi _0 (\psi)\] and if $0<
\rho \leq r_0 < \rho _0 = \delta r_0 < r < \infty$, then
\[g(r, \theta ; \rho ,\psi )\leq b r^{\mu ^-} \rho ^{\mu ^+ }
\varphi _0 (\theta ) \varphi _0 (\psi ).\]
\end{lem}
\noindent \textbf{Proof.} We apply here a method used by
Rashkovskii \cite{Rash}. Since the inward normal derivative
$\frac{\partial \varphi _0 (\theta )}{\partial n}\bigl|_{\partial
D}$ is continuous and positive on the entire boundary $\partial
D$, the function
\[S(\theta ,\psi )=\frac{g(1, \theta ; \delta , \psi )}
{\varphi _0 (\theta ) \varphi _0 (\psi)}\] is continuous and,
therefore, bounded in $\overline{D} \times \overline{D}$. Thus,
$g(1, \theta ; \delta , \psi ) \leq b \varphi _0 (\theta )\varphi
_0 (\psi )$. This inequality and Lemma 4.3 imply, with any $r_0$,
\[g(r_0, \theta ; \delta r_0, \psi ) = r^{2-n}_0 g(1, \theta ; \delta, \psi )
\leq b r^{2-n} _0 \varphi _0 (\theta ) \varphi _0 (\psi ) \] and
after setting $\rho _0= \delta r_0$ we arrive at the inequality
\[g(r_0, \theta ; \rho _0,\psi )\leq b r_0 ^{\mu ^+} \rho _0 ^{\mu
^-}\varphi _0 (\theta )\varphi _0 (\psi ),\] since $\mu ^+ + \mu
^- =2-n$.

Fix an $r_0 >0$, $\rho _0= \delta r_0$ and consider a function
\[g(r, \theta ; \rho _0, \psi )- b \rho _0 ^{\mu
^-} \varphi _0 (\psi ) r^{\mu ^+} \varphi _0 (\theta ).\]
\pagebreak It is harmonic when $x=(r,\theta ) \in K^D _{r_0}$,
vanishes at its lateral surface and nonpositive when $r=r_0$.
Thus,
\begin{equation}
g(r, \theta; \rho _0, \psi )- b \rho _0 ^{\mu ^-} \varphi _0 (\psi
) r^{\mu ^+} \varphi _0 (\theta ) \leq 0
\end{equation}
for all $r \leq r_0$ and $\theta ,\psi \in D$. Now we unfreeze
$\rho$ and consider a function
\[g(r, \theta; \rho , \psi )- b \rho ^{\mu
^-} \varphi _0 (\psi ) r^{\mu ^+} \varphi _0 (\theta )\] with
fixed $r \leq r_0$ and $\theta \in D$. It is harmonic in $y=(\rho
, \psi )$ as $\rho > \rho _0$ and $\psi \in D$, vanishes at the
boundary of $\partial K^D \setminus \partial K^D _{\rho _0}$,
uniformly tends to zero as $\rho \rightarrow \infty$ and is
nonpositive at $\rho = \rho _0$ by virtue of (4.15). Thus,
\begin{equation}
g(r, \theta; \rho , \psi ) \leq b \rho ^{\mu ^-} \varphi _0 (\psi
) r^{\mu ^+} \varphi _0 (\theta )
\end{equation}
when $r\leq r_0$ and $\rho \geq \delta r_0$. Since $r_0 >0$ is
arbitrary, (4.16) is valid for all $\rho \geq \delta r_0; \;
\delta >1$. The second assertion of the lemma is due to the
symmetry of Green's function. \hfill \textbf{$\diamondsuit$}

\begin{theor}\hspace{-.09in}\textbf{.} If $\partial D\in C^2$ and a
constant $\delta > 1$, then
\begin{equation}
G(r, \theta ; \rho , \psi ) \leq b \; V(r) W(\rho ) \varphi _0
(\theta ) \varphi _0 (\psi ), \; 0 < r \leq r_0 \leq \delta r_0
\leq \rho < \infty ,
\end{equation}
\begin{equation}
G(r, \theta ; \rho , \psi ) \leq b \; V(\rho ) W(r) \varphi _0
(\theta ) \varphi _0 (\psi ), \; 0 < \rho \leq r_0 \leq \delta r_0
\leq r < \infty .
\end{equation}
\end{theor}

\noindent \textbf{Proof.} If $r_0$ and $\rho _0$ are fixed, $\rho
_0 =\delta r_0, \; \delta > 1$, then from Theorem A.2 and Lemma
4.4 we get
\[G(r_0 , \theta ; \rho _0 , \psi ) \leq g(r_0 , \theta ; \rho _0 , \psi
)\]
\[\leq b r_0 ^{\mu ^+} \rho _0 ^{\mu ^-} \varphi _0 (\theta )
\varphi _0 (\psi ) \leq b' V(r_0 ) W(\rho _0 )\varphi _0 (\theta )
\varphi _0 (\psi ),\] where $b' =b r_0 ^{\mu ^+} \rho _0 ^{\mu ^-}
\left(V(r_0 ) W(\rho _0)\right)^{-1}$. Therefore, the difference
\[G(r, \theta ; \rho _0 , \psi ) - b' V(r) W(\rho _0 )
\varphi _0 (\theta ) \varphi _0 (\psi ),\] which is an
$L_q-$harmonic function of the argument $x = (r, \theta )$ in $K^D
_{r_0}$, has nonpositive boundary values in this cone. From the
maximum principle we obtain the inequality
\[G(r, \theta ; \rho _0 , \psi ) \leq b' V(r) W(\rho _0 )
\varphi _0 (\theta ) \varphi _0 (\psi )\] for all $x \in K^D _{r_0}$.\\

To establish (4.17) with an arbitrary $\rho$ instead of a fixed
$\rho _0$, we apply the same argument to a $q-$harmonic function
\[G(r, \theta ; \rho , \psi ) - b V(r) W(\rho )\varphi _0 (\theta )
\varphi _0 (\psi )\] for $y = (\rho , \psi )\in K^D _{\rho _0,
\infty }, \; \rho _0 \leq \rho $. Finally, (4.18) follows from
(4.17) due to the symmetry of $G(x;y)$. \hfill \textbf{$\diamondsuit$}\\

\noindent \textbf{Remark 4.2.} The Oleinik theorem can be applied
directly to the quotient
\[\frac{G(x; \rho _0, \psi )}{\varphi _0 (\theta )\varphi _0 (\psi )}\]
if we observe that
\[G(r_0 , \theta ; \rho _0 , \psi ) \leq b(D, r_0 , \delta )
V(r_0) W(\rho _0 )\varphi _0 (\theta ) \varphi _0 (\psi ), \;
\delta > 1; \] our proof shows in addition that the constants $b$
do not depend on $r_0$. \hfill \textbf{$\diamondsuit$}

\begin{cor}\hspace{-.09in}\textbf{.} Combining (4.17) or (4.18) with (C.1)
we obtain, respectively, inequalities
\begin{equation}
G(r, \theta ; \rho , \psi ) \leq b r^{2-n} W(\rho ) W^{-1} (\delta
 r) \varphi _0 (\theta ) \varphi _0 (\psi ), \; 0 < r \leq r_0 \leq \delta
r_0 \leq \rho < \infty
\end{equation}
and
\[G(r, \theta ; \rho , \psi ) \leq b r^{2-n} V(r) V^{-1} (\rho / \delta )
\varphi _0 (\theta ) \varphi _0 (\psi ), \; 0 < \rho \leq r_0 \leq
\delta r_0 \leq r < \infty .\] \hfill \textbf{$\diamondsuit$}
\end{cor}

\begin{cor}\hspace{-.09in}\textbf{.} In the case \textbf{(A)}, using the
formulas (C.9)-(C.10) we have a generalization of a known result
for the Green function of the Laplacian \cite{Azar}:
\[G(r, \theta ; \rho , \psi ) \leq b \rho ^{2-n}
\left(\frac{r}{\rho } \right)^{\mu ^+ _k } \varphi _0 (\theta )
\varphi _0 (\psi ), \; k r \leq \rho ,\] where $\mu _k ^+ =
\frac{1}{2} \left(2-n + \sqrt{(n-2)^2 + 4(\lambda _0 + k)}
\right)$. In the case \textbf{(B)} (4.19) takes the form
\[G(r, \theta ; \rho , \psi ) \leq b \frac{r^{1-n}}
{\left(q(r) q(\delta r) \right)^{1/4}} \left(e^{- \int _r ^{\delta
r} \sqrt{q(t)} dt}\right) \varphi _0 (\theta ) \varphi _0 (\psi ),
\; \delta r \leq \rho .\] \hfill \textbf{$\diamondsuit$}
\end{cor}

\begin{cor}\hspace{-.09in}\textbf{.} If $\partial D \in C^2$ and
$0 < r \leq r_0 < \delta r_0 \leq \rho < \infty$, then
\begin{equation}
\frac{\partial G(r, \theta ; \rho , \psi )}{\partial n(y)} \leq b
V(r) W(\rho ) \varphi _0 (\theta ) \frac{1}{\rho } \;
\frac{\partial \varphi _0 (\psi )}{\partial n(\psi )}, \; \theta
\in D, \; \psi \in \partial D
\end{equation}
and
\begin{equation}
\frac{\partial G(r, \theta ; \rho , \psi )}{\partial n(y)} \leq b
V(\rho ) W(r) \varphi _0 (\theta ) \frac{1}{\rho} \;
\frac{\partial \varphi _0 (\psi )}{\partial n(\psi )}, \; \theta
\in D, \; \psi \in \partial D,
\end{equation}
if $0 < \rho \leq r_0 < \delta r_0 \leq r < \infty $.
\end{cor}

\noindent \textbf{Proof.} It is enough to note that a nonpositive
$q-$harmonic function
\[G(x;y)-bV(r)W(\rho)\varphi _0(\theta)\varphi _0(\psi )\]
with zero boundary values in $K^D _{r_0 ,\infty}$ has a
nonpositive inward normal derivative \cite[p. 14]{Mir}. We mention
also that $\frac{\partial \varphi _0 (x/r)}{\partial n(x)} =
\frac{1}{r} \; \frac{\partial \varphi _0
(\theta)}{\partial n(\theta)}$. \hfill \textbf{$\diamondsuit$}\\

A proof of the following statement is similar.

\begin{cor}\hspace{-.09in}\textbf{.} Let $G_R (x;y)$ be the Green function
of $L_q$ in a cone $K^D _R$, \\
$0 < R < \infty$. If $\partial D \in C^2$ and $0 < r \leq \gamma
R, \; \gamma < 1$, then
\[- \frac{\partial }{\partial \rho } G_R (r, \theta ;
\rho , \psi ) \biggl|_{\rho = R} \leq b(\gamma ) V(r) \left\{-
W'(R) \right\} \varphi _0 (\theta ) \varphi _0 (\psi ). \] \hfill
\textbf{$\diamondsuit$}
\end{cor}

\noindent \textbf{Remark 4.3.} Due to Theorem A.2, all upper
bounds of this section hold for the Green function of the operator
$L_c$ with any $c(x) \geq q(|x|)$. \hfill \textbf{$\diamondsuit$}\\

\noindent \textbf{Remark 4.4.} If $r=|x|$ is fixed, then $G(x;y)$
vanishes as $\underline{\underline{O}}\left(W(\rho ) \right)$ when
$\rho = |y| \rightarrow
\infty$. \hfill \textbf{$\diamondsuit$}\\

We will need lower bounds of Green's functions as well.

\begin{lem}\hspace{-.09in}\textbf{.} If $0<\gamma <1$ and $\rho _0>0$
are fixed, then
\[G(r, \theta ; \rho , \psi ) \geq \left\{
    \begin{array}{ll}
        b V(r) W(\rho ) \varphi _0 (\theta ) \varphi _0 (\psi )
            & \mbox{if  $\; r \leq \gamma \rho _0 < \rho _0 \leq \rho$ } \vspace{.5cm} \\
        b V(\rho ) W(r) \varphi _0 (\theta ) \varphi _0 (\psi )
            & \mbox{if  $\; \rho \leq \gamma \rho _0 < \rho _0 \leq r$.}
                                \end{array}
                  \right. \]
\end{lem}

\noindent\textbf{Proof.} We can repeat the argument of Lemma 4.4 -
Theorem 4.2 if we notice that the quantity $S(\theta , \psi )$ in
Lemma 4.4 and its analog for the Green function $G$ in Theorem 4.2
are bounded away from zero as well. \hfill \textbf{$\diamondsuit$}

\begin{cor}\hspace{-.09in}\textbf{.} Under the conditions of Lemma 4.5, if
$r \leq \gamma \rho _0 < \rho _0$, then
\[-\frac{\partial}{\partial \rho}G_R(r,\theta ;\rho ,\psi)\biggl|_{\rho = R}
\geq bV(r)\left\{- W'(R) \right\} \varphi _0 (\theta ) \varphi _0
(\psi ), \; b = b(K,\rho _0), \] where $G_R$ was defined in
Corollary 4.8.
\end{cor} \hfill \textbf{$\diamondsuit$}

\begin{cor}\hspace{-.09in}\textbf{.} If $y=(\rho , \psi ) \in \partial K^D$,
then
\[\frac{\partial G(r, \theta ; \rho , \psi )}{\partial n(y)}
\geq \left\{   \begin{array}{ll}
        b V(r) W(\rho ) \varphi _0 (\theta ) \frac{1}{\rho }
            \frac{\partial \varphi _0 (\psi )}{\partial n(\psi )}
                & \mbox{if  $\; |x| < |y|$} \vspace{.5cm} \\
        b V(\rho ) W(r ) \varphi _0 (\theta ) \frac{1}{\rho }
            \frac{\partial \varphi _0 (\psi )}{\partial n(\psi )}
                & \mbox{if  $\; |y| < |x|$.} \\
                                \end{array}
                  \right. \]
\end{cor} \hfill \textbf{$\diamondsuit$}\\

Using these inequalities, we can straightforwardly verify that the
series (4.10) can be termwise differentiated.

\begin{cor}\hspace{-.09in}\textbf{.} Under the conditions of Theorem 4.1,
for $y=(\rho , \psi ) \in \partial K^D$
\[\frac{\partial G(r, \theta ; \rho , \psi )}{\partial n(y)}
= \left\{ \begin{array}{ll}
        \sum ^{\infty } _{\nu = 0} \frac{1}{\chi' _{\nu}} \;
            V_{\nu}(r) W_{\nu }(\rho ) \varphi _{\nu }(\theta ) \;
                \frac{1}{\rho} \; \frac{\partial \varphi _{\nu}
                    (\psi)}{\partial n(\psi)} & \mbox{if $ \;
                        0 < r < \rho < \infty $} \vspace{.5cm} \\
        \sum ^{\infty } _{\nu = 0} \frac{1}{\chi' _{\nu}} \;
            V_{\nu}(\rho) W_{\nu }(r) \varphi _{\nu }(\psi) \;
                \frac{1}{r} \; \frac{\partial \varphi _{\nu}
                    (\theta)}{\partial n(\theta)} & \mbox{if $ \;
                        0 < \rho < r < \infty $.} \\
           \end{array}
  \right. \]
\end{cor} \hfill \textbf{$\diamondsuit$}\\

As a simple application of Theorem 3.1 and Corollary 4.7, we prove
an extension of a uniqueness theorem for analytic functions
\cite[p. 219-220]{Evg1}, which will be needed in Section 5. Here
again $V$ is a growing solution of (3.4) with $\lambda =\lambda
_0$.
\begin{prop}\hspace{-.09in}\textbf{.} Let $u\in \verb"SbH"(c,K^D)$ with
$c\in {\cal C}(K^D,\mathbf{A})$ and the boundary $\partial D\in
C^2$. If $u(x) \leq b V(|x|)$ for an arbitrary small positive
constant $b$, then $u(x) \equiv - \infty $.
\end{prop}

\noindent \textbf{Proof.} The convergence of the integral
\[I(x) = \int _{\partial K^D}V(|y|)\frac{\partial G(x;y)}{\partial n(y)}d\sigma (y),\]
where $G$ is the Green function of the operator $L_q$ in the cone
$K^D$, follows by a direct estimation based on Lemma C.2 and
estimates (4.20)-(4.21). Thus, $I(x)$ solves the Dirichlet problem
for $L_q$ in $K^D$ with the boundary data $V(|x|)$. Therefore, a
function $bI(x)$ with any $b>0$ is a $q-$harmonic majorant of $u$
in $K^D$.

Consider now a $q-$harmonic function $v(x)= -bI(x)+AV(|x|)\varphi
_0 (\theta)$ with a positive constant $A$ and exhaust the domain
$D$ from within by smooth domains $D_j \nearrow D$. With $b$ and
$A$ fixed, the inequality $A\varphi _0 (\theta) <b$ holds
everywhere in $D\setminus \overline{D_j}$ if $j$ is sufficiently
large. Therefore, $v(x) \leq 0$ at the boundary $\partial
K^{D_j}$. Since $\lambda _0 < \lambda _0 ^{[j]}$, where $\lambda
_0 ^{[j]}$ is the smallest eigenvalue of the Laplace-Beltrami
operator in $D_j$, we have $V (r) =
\bar{\bar{o}}\left(V^{[j]}(r)\right), \; r\rightarrow \infty$, and
so that
\[u(x) \leq bI(x) = v(x) - AV(r) \varphi _0 (\theta) \leq
- AV(r) \varphi _0 (\theta).\] It remains to let $A \rightarrow
\infty $. \hfill \textbf{$\diamondsuit$}

\newpage

\section{The Blaschke theorem}

This theorem asserts that given a sequence of points $z_k, \;
k=1,2,...$, in the unit disk $B(0,1) \subset \mathbf{R}^2$, in
order for a bounded analytic function with zeros at these and only
these points to exist it is necessary and sufficient that these
points satisfy the Blaschke condition $\sum _{k=1}
^{\infty}\left(1-|z_k| \right)< \infty $. In subharmonic setting,
the Riesz associated measure of a bounded subharmonic function in
the unit disk must satisfy the condition
\[\int _{B(0,1)} \left(1-|z| \right)d\mu (z)< \infty .\]
Vice versa, if a Borel measure $d\mu $ in $B(0,1)$ satisfies this
condition, then there exists an upper bounded subharmonic function
in the unit disk such that its associated Riesz measure coincides
with $d\mu$ - see, e.g., \cite[p. 174]{Hr}. By making use of
conformal mapping one can replace the unit disk here with other
plane domains.

The same questions make sense in $\mathbf{R}^n$ - thus, an upper
bounded subharmonic function in $\mathbf{R}^n, \; n \geq 3$, whose
Riesz associated measure coincides with a given Borel measure
$d\mu$, exists if and only if \cite[p. 128]{HK}
\begin{equation}
\int _1 ^{\infty} r^{1-n} n(r) dr < \infty , \mbox{ where }
n(r)=\mu (B(r))=\int _{B(r)}d\mu (x).
\end{equation}
In this section we prove an analog of the Blaschke theorem for
subfunctions of the operator $L_c$ in $n-$dimensional cones. To
derive a complete analog of the Blaschke theorem, we limit
ourselves here by subfunctions $u(x) \in \verb"SbH"(q, K^D)$ with
respect to a radial potential $q(r)$ - the results concerning the
general case are discussed in Remark 5.4 at the end of this
section. The argument is based on an extension of the following
Carleman formula.

For a holomorphic function $f(z)$ in a closed semi-annulus
\[\Omega = \left\{z \in R^2 \biggl| \; 0 < a \leq |z| \leq R, \; \Im z \geq 0 \right\}\]
with the zeros $z_k = r_k e^{i \theta _k }, \; k = 1,2,...,l$,
this formula reads \cite{Carl}:
\begin{equation}
\begin{array}{c}
\sum _{a<r_k <R} \left(\frac{1}{r_k} - \frac{r_k }{R^2} \right)
\sin \theta _k = \frac{1}{\pi R} \int _0 ^{\pi } \ln |f(R e^{i
\theta })| \sin \theta \; d \theta \vspace{.5cm} \\
+ \frac{1}{2 \pi} \int ^R _a \left(\frac{1}{x^2} - \frac{1}{R^2}
\right) \ln |f(x)f(-x)|dx + A_f (a, R),
\end{array}
\end{equation}
where $A_f$ is a constant. As for generalizations and applications
of (5.2) see, for example, \cite{Ahl2, GoOs, Gov, Khe73, L1,
Thom}. In particular, the Carleman formula was extended onto
subharmonic functions in smooth cones in $\mathbf{R}^n$
\cite{RaRo}. We establish two extensions of (5.2) for
subfunctions, one is convenient for studying the asymptotic
behavior of subfunctions at infinity, and another one for that at
the origin. Here again $V=V_0 (r)$ and $W=W_0 (r)$ are the
normalized solutions of the equation (3.4) and $\chi ' =
w(W,V)|_{r=1}$ is their Wronskian. As before, truncated cones are
denoted by
$K^D _R$ and $K^D _{a, R}, \; 0 < a < R \leq \infty$.\\

Consider a $q$-harmonic function
\[C^R (x) = \left(W(r) - \frac{W(R)}{V(R)}V(r) \right) \varphi _0 (\theta ),
\; x=(r,\theta )\in K^D. \] It is clear that $C^R(x)>0$ in the
cone $K^D _{1,R}$, $C^R (r, \theta )=0$ as $\theta \in \partial D$
and $0<r<R$, and $C^R (R, \theta )=0$. Let $d\sigma = d\sigma (r,
\theta )$ denote the surface area measure on the sphere $S^{n-1}$
or on the lateral surface $\partial K^D$ of the cone $K^D$, and
let $dx$ be the Lebesgue measure in $\mathbf{R}^n$. For a
subfunction $u$ we denote by $d\mu$ its Riesz measure. First we
consider continuous subfunctions in a closed cone with the smooth
boundary. As usual, $u$ is said to be a subfunction in a closed
set if it is a subfunction in an open domain containing this set.

\begin{lem}\hspace{-.09in}\textbf{.} If $\partial D$ is smooth and
$u(x) \in \verb"SbH"(q, \; \overline{K^D _{a,R}}), \; 0<a<R<\infty
$, then
\begin{equation}
\begin{array}{c}
\theta _n \int _{K^D _{a,R}} C^R (x) d\mu (x) = \frac{\chi
'}{V(R)} \int _D u(R, \theta ) \varphi _0 (\theta ) d\theta + \vspace{.5cm} \\
\int _{\partial D \times (a,R)} u(r, \theta ) \left(W(r) -
\frac{W(R)}{V(R)}V(r) \right) \frac{1}{r} \frac{\partial \varphi
_0 (\theta )}{\partial n} d\sigma (r, \theta ) + A_u (a, R),
\end{array} \end{equation}
where the constant $\theta _n$ was defined in Theorem 2.1, Part
$8^0$, and
\[A_u (a, R) = a^{n-1} \int _D \left(u(a,\theta )
\frac{\partial C^R (r, \theta )}{\partial r} \biggl|_{r=a} - \;
\frac{\partial u(r, \theta }{\partial r} \biggl|_{r=a} C^R (a,
\theta ) \right) d\theta = \underline{\underline{O}}(1) \] as $R
\rightarrow \infty$ while $a>0$ is fixed.
\end{lem}

\textbf{Proof.} If $u \in C^2(\overline{K^D _{a,R}})$, then to
derive (5.3) it suffices to apply the second Green identity in
$K^D _{a,R}$ to the functions $u$ and $C^R$ and notice that, due
to Part $8^o$ of Theorem A.1, for any continuous $f$
\[\int _{K^D _{a,R}}f(x)\left(-L_q u(x)\right)d(x)
=\theta _n\int_{K^D_{a,R}}f(x)d\mu (x).\]

If $u$ is merely upper semicontinuous in $\overline{K^D _{a,R}}$,
then we approximate it by a decreasing sequence of continuous
subfunctions $u_j$, whose associated measures $d\mu _j$ weakly$^*$
converge to $d\mu$ by the Helly selection theorem. All terms in
(5.3), but the last one, have limits due to the monotone
convergence $u_j(x) \searrow u(x)$, and so that the limit $\lim
_{j \rightarrow \infty} \int _D \left(\frac{\partial u_j (r,
\theta )}{\partial r} \bigl|_{r=a} \right) C^R (a,\theta )d\sigma
(\theta )$ also exists. It should be also noticed that like the
subharmonic case, the partial derivatives $\frac{\partial u_j (r,
\theta )}{\partial r}$ exist almost everywhere and are locally
summable - see the proof of Theorem A.1. \hfill \textbf{$\diamondsuit$}\\

The next lemma has the same proof; here
\[C_a(r,\theta )=\left(V(r)-\frac{V(a)}{W(a)}W(r) \right)
\varphi _0 (\theta ), \; 0<a<r< \infty . \]
\begin{lem}\hspace{-.09in}\textbf{.} If $\partial D$ is smooth and
$u(x) \in \verb"SbH"(q, \; \overline{K^D}), \; 0<a<R<\infty $,
then
\[\begin{array}{c}
\theta _n \int _{K^D _{a,R}} C_a (x) d\mu (x) = \frac{\chi
'}{W(a)}\int _D u(a, \theta ) \varphi _0 (\theta ) d\theta + \vspace{.5cm} \\
\int _{\partial D \times (a,R)} u(r, \theta ) \left(V(r) -
\frac{V(a)}{W(a)}W(r) \right) \frac{1}{r} \frac{\partial \varphi
_0 (\theta )}{\partial n} d\sigma (r, \theta ) + B_u (a, R),
\end{array}  \]
where
\[B_u (a, R) = R^{n-1} \int _D \left(u(R,\theta )
\frac{\partial C_a (r, \theta )}{\partial r} \biggl|_{r=R} - \;
\frac{\partial u(r, \theta }{\partial r} \biggl|_{r=R} C_a (R,
\theta ) \right) d\theta = \underline{\underline{O}}(1) \] as $a
\rightarrow 0$ and $R>0$ is fixed. \hfill \textbf{$\diamondsuit$}
\end{lem}

Now we state main results of this section.
\begin{theor}\hspace{-.09in}\textbf{.} Let $d\mu$ be a Borel measure in
a cone $\overline{K^D}$ such that $\partial D$ is not a polar set.
Let also $0 \leq q(r) \in L^p _{loc} (0, \infty ), \; p>n$.  An
upper bounded $L_q -$subfunction $u$ in $K^D$ with the Riesz
associated measure $d\mu$ exists if and only if
\begin{equation}
\int _{K^D _b} V(r) \varphi _0 (\theta ) d\mu (r, \theta) < \infty
\end{equation}
and
\begin{equation}
\int _{K^D _{b, \infty }} W(r) \varphi _0 (\theta ) d\mu (r,
\theta) < \infty
\end{equation}
for some, and so for any constant $b>0$.
\end{theor}

For harmonic functions in the half-space the necessity of
conditions (5.4)-(5.5) was established in \cite{LF1} (Cf. (4)
there). The necessity part of Theorem 5.1 follows from the next
result.

\begin{theor}\hspace{-.09in}\textbf{.} Let $D$ be an arbitrary domain on
$S(0, 1)$ whose boundary is not a polar set, and a Borel measure
$d\mu (x)$ is the Riesz associated measure of a function
$u(x) \in \verb"SbH"(q, K^D)$.\\

$1^o$ If for each $x_0 \in \partial K^D , \; |x_0 | > b$, there
holds a condition
\[\limsup _{K^D _{b, \infty} \ni x \rightarrow x_0} u(x) \leq A = const < \infty ,\]
then
\begin{equation}
\int _{K^D _{b,R}}C^R(x)d\mu (x) \leq \frac{\chi '}{\theta _n
V(R)} \int _D u(R, \theta ) \varphi _0 (\theta ) d\theta +
\underline{\underline{O}}(1), \; R \rightarrow  \infty .
\end{equation}

$2^o$ If for each $x_0 \in \partial K^D, \; |x_0 | < b$, there
holds a condition
\[\limsup _{K^D _b \ni x \rightarrow x_0} u(x) \leq A = const < \infty ,\]
then
\begin{equation}
\int _{K^D _{a, b}} C_a (x) d\mu (x) \leq \frac{\chi '}{\theta _n
W(a)} \int _D u(a, \theta ) \varphi _0 (\theta ) d \theta +
\underline{\underline{O}}(1), \; a \rightarrow 0^+ .
\end{equation} \end{theor}

\noindent \textbf{Proof.} We shall prove only Part $1^o$, since
the proof of Part $2^o$ is similar. Without loss of generality we
can assume that $b=1$ and also that $A<0$, since otherwise instead
of $u$ we can consider a $q$-subfunction $v(x)=u(x)-M$ with $M>A
\geq 0$. The function $v$ (or the $u$ itself if $A<0$) is negative
in some vicinity of the boundary $\partial K^D$. Therefore (Cf.
the proof of Theorem 3.1), there exists a sequence of domains
$D_j$ with smooth boundaries such that the boundaries of cones
$\partial
K^{D_j}$ are in the domain where $v(x) < 0$.\\

Now we apply the formula (5.3) with $a=1$ to ${v}(x)$ in $K^{D_j}
_{1,R}$:
\begin{equation}
\begin{array}{c}
\theta _n \int _{K^{D_j} _{1, R}} C^R _j (x) d\mu _v (x) =
\frac{\chi ' _j}{V^{[j]}(R)} \int _{D_j } v(R, \theta ) \varphi
^{[j]} _0 (\theta )
d\theta + \vspace{.5cm} \\
\int _{\partial K^{D_j} _{1, R}} \left(W^{[j]}(r) -
\frac{W^{[j]}(R)}{V^{[j]}(R)}V^{[j]}(r) \right) \frac{1}{r}
\frac{\partial \varphi ^{[j]} _0 (\theta )}{\partial n(\theta )}
v(r, \theta ) d\sigma (r,\theta )+ A^{[j]}(1, R),
\end{array}
\end{equation}
where $A^{[j]}(1, R) = \int _{D_j} \left(C^R _j (1, \theta
)\frac{\partial v(r, \theta }{\partial r} \biggl|_{r=1}  - \;
v(1,\theta ) \frac{\partial C^R _j (r, \theta )}{\partial r}
\biggl|_{r=1} \right) d\theta$.

Since $\partial D$ is smooth, $\frac{\partial \varphi ^{[j]} _0
(\theta )}{\partial n(\theta )} > 0$. Moreover, $v(r,\theta )<0$
for $\theta \in \partial D_j$ by the construction. Therefore, the
second addend on the right in (5.8) is negative, and (5.8) implies
the inequality
\begin{equation}
\theta _n \int _{K^{D_j}_{1,R}}C_j ^R (x)d\mu _{v}(x) \leq
\frac{\chi ' _j }{V^{[j]}(R)}\int _{D_j}v(R,\theta ) \varphi
^{[j]} _0 (\theta )d\theta + A^{[j]}(1,R).
\end{equation}
Extending here all functions by zero, we can rewrite (5.9)
with $D_j$ replaced at all occurrences by $D$. \\

Since by Theorem C.1, Part $3^o$, and Theorem B.1, Part $4^o$,
$V^{[j]}(r) \rightarrow V(r), \; W^{[j]}(r) \rightarrow W(r)$, and
$\varphi ^{[j]} _0 (\theta ) \rightarrow \varphi _0 (\theta )$
uniformly on compact sets as $j \rightarrow \infty $, we conclude
that it is possible to let $j \rightarrow \infty $ and so $\bigl|
A^{[j]} (1, R) \bigl| \leq A(1,R) = O(1)$ as $R \rightarrow \infty
$. Thus,
\[\int_QC_j^R(x)d\mu_{v}(x)\rightarrow \int_Q C^R(x)d\mu_v(x) \;
\mbox{ as $j \rightarrow \infty $,}\] for any compact set $Q
\subset K^D$. Now, since $\parallel \varphi ^{[j]} _0 (\theta ) -
\varphi _0 (\theta )\parallel _{L^2 (D)} \rightarrow 0$ as $j
\rightarrow \infty $ by Theorem B.1, Part $3^o$,  we get from
(5.9) an inequality
\[\theta _n \int _Q C^R (x) d\mu _v (x) \leq \frac{\chi '}{V(R)}
\int _D v(R, \theta ) \varphi _0 (\theta ) d\theta +
\underline{\underline{O}}(1).\] Exhausting the cone $K^D _{1, R}$
with compact sets leads finally to the inequality
\begin{equation}
\int _{K^D _{1, R}} C^R (x) d\mu _v(x) \leq \frac{\chi '}{\theta
_n V(R)} \int _D v(R, \theta ) \varphi _0 (\theta ) d\theta +
\underline{\underline{O}}(1), \; R \rightarrow \infty .
\end{equation}
To replace here $v(x)$ with $u(x)$, it is enough to notice that,
by virtue of  Theorem A.1, Part $7^o$, $d\mu _u(x) \leq d\mu
_v(x)$, and (5.6) straightforwardly follows from (5.10). \hfill
\textbf{$\diamondsuit$}

\begin{cor}\hspace{-.09in}\textbf{.} If
\begin{equation}
\liminf _{R \rightarrow \infty } \frac{1}{V(R)} \int _D u^+ (R,
\theta ) \varphi _0 (\theta ) d\theta < \infty ,
\end{equation}
in particular, if $u(x)$ is upper bounded in $K^D$, then (5.5)
holds.
\end{cor}

\noindent \textbf{Proof.} It follows from (5.6) and (5.11) that
\[\int _{K^D _{1, R}} \left(W(r) - \frac{W(R)}{V(R)}V(r) \right)
\varphi _0 (\theta ) d\mu (r, \theta ) \leq b < \infty .\]

Since $\mu ^+ >0$ and by Lemma C.1, Part $1^o$, $V(r)/V(R) \leq
(r/R)^{\mu ^+}$, for a given $\rho > 1$ we can choose $R > \rho$
such that $V(\rho )/V(R) < 1/2$ and so
\[(1/2) \int _{K^D _{1, \rho }} W(r) \varphi _0 (\theta )
d\mu (r, \theta ) \leq b.\]

The latter integral does not decrease when $\rho $ increases, thus
\\ $\int _{K^D _{1,\infty}}W(r)\varphi _0(\theta )d\mu
(r,\theta )\leq 2b <\infty$. \hfill \textbf{$\diamondsuit$}\\

The next result follows in the same way from (5.7).

\begin{cor}\hspace{-.09in}\textbf{.} If
\begin{equation}
\liminf _{a \rightarrow +0} \frac{1}{W(a)} \int _D u^+ (a, \theta
) \varphi _0 (\theta ) d\theta < \infty ,
\end{equation}
in particular, if $u(x)$ is upper bounded in $K^D _1$, then (5.4)
holds. \hfill \textbf{$\diamondsuit$}
\end{cor}

\noindent \textbf{Remark 5.2.} Vice versa, for nonnegative
subfunctions the statements (5.4) and (5.5) imply, respectively,
(5.11) and (5.12). Moreover, we can replace the lower limit with
the limit in these relations. \hfill \textbf{$\diamondsuit$}\\

\noindent \textbf{Remark 5.3.} The condition (5.1) follows from
(5.5) by partial integration, since if $D = S(1)$, then
$\varphi _0=const >0$ and $W(r)=r^{2-n}$. \hfill \textbf{$\diamondsuit$}\\

Theorem 5.2 has the following converse.

\begin{theor}\hspace{-.09in}\textbf{.} If $\partial D$ is not a polar set,
then for any Borel measure $d\mu (x)$ satisfying the conditions
(5.4)-(5.5) in a cone $K^D$ there exists an upper bounded function
$u \in \verb"SbH"(q, K^D)$ such that its associated measure is
$d\mu (x)$. The Green potential
\[u(x) = - \int _{K^D} G(x;y) d\mu (y), \] where $G$ is
Green's function of the Dirichlet problem for $L_q$ in $K^D$, is
such a function, and moreover, $\sup _{x\in K^D}u(x)=0$.
\end{theor}

\noindent \textbf{Proof.} First, we consider smooth cones with
boundaries $\partial D \in C^2$. Introduce measures
\[d\mu _1 (\rho, \psi) = \left\{ \begin{array}{ll}
        d\mu (\rho, \psi) & \mbox{if $\rho \geq 1, \; \psi \in D$} \vspace{.25cm} \\
        0 & \mbox{if $0 \leq \rho < 1, \; \psi \in D$}
                                \end{array}
                  \right. \]
and $d\mu _2 = d\mu - d\mu _1$. By making use of (5.5) and (4.17)
with $r_0 = 1/2$ and $r = 1/4$ we conclude that the integral
\begin{equation}
U_1(x) = \int _{K^D_{1, \; \infty}} G(x;y) d\mu _1(y)
\end{equation}
converges and represents a positive $q-$superfunction in $K^D$
which may take on value $+ \infty $ only on a polar set. The
integral
\begin{equation}
U_2(x) = \int _{K^D _1} G(x;y) d\mu _2(y)
\end{equation}
converges for $r = |x| > 1$. Therefore, it may be infinite only if
$x$ belongs to  a polar set, and $U_2$ is also a positive
$q-$superfunction in $K^D$. Thus, the function $u(x) = -\left\{U_1
(x)+ U_2 (x) \right\}$ is a nonpositive $q-$subfunction in $K^D$
with the given associated measure $d\mu $.

To remove the restriction $\partial D \in C^2$, we again, as in
the proof of Theorem 3.1, approximate the cone $K^D$ from within
with a sequence of smooth cones $K^{D_j}, \; \partial D_j \in C^2
, \; D_j \nearrow D$. Due to the variational principle, $\lambda
^{[j]} _0 \searrow \lambda _0$. Therefore, taking also into
account the normalization $V^{[j]}(1) = W^{[j]}(1) = 1$, the
function $V^{[j]}(r)$ vanishes as $r \rightarrow 0^+$ faster than
$V(r)$, and the function $W^{[j]}(r)$ vanishes as $r \rightarrow
\infty $ faster than $W(r)$. Considering again the uniform
convergence $\varphi ^{[j]} _0 (\theta ) \rightarrow \varphi _0
(\theta )$ on compact sets $Q \subset D$, we obtain from
(5.4)-(5.5)
\[\int _{K^{D_j}_R}V^{[j]}(r)\varphi ^{[j]}_0 (\theta )d\mu (r,\theta )<\infty ,
\; j=1,2,...\]
and
\[\int _{K^{D_j}_{R,\infty }}W^{[j]}(r)\varphi ^{[j]}_0 (\theta )d\mu (r,\theta )<\infty ,
\; j = 1, 2, ... \; .\]

Next, by making use of the formulas (5.13)-(5.14), we construct a
sequence of superfunctions
\[U_j(x) = \int _{K^D} G^{[j]}(x;y) d\mu (y), \]
where Green's functions $G^{[j]}$ of the cones $K^{D_j}$ are
extended onto $K^D \setminus K^{D_j}$ by zero. The sequence
$\left\{G^{[j]}\right\}^{\infty}_{j=1}$ monotonically increases as
$j\rightarrow \infty $. So, the sequence $\left\{U_j
(x)\right\}^{\infty}_{j=1}$ also monotonically increases and tends
either to a positive superfunction $U(x)$ with the associated
measure $d\mu (x)$ or to infinity.

In the latter case, however, the functions
\[P_j (r) \equiv \int _D U_j (r,\theta ) \varphi ^{[j]} _0 (\theta )
d\theta \rightarrow +\infty \] as well, which is impossible as we
show now. Indeed, by the Fubini theorem,
\[P_j (r) = \int _{K^D} \int _D G^{[j]}(x;y) \varphi ^{[j]} _0
(\theta ) d\theta d\mu (\rho , \psi ). \] The eigenfunctions,
$\varphi ^{[j]} _{\nu } (\theta )$, of the Laplace-Beltrami
operator are orthonormal in $L^2 (D_j)$ and extended by zero into
$D \setminus D_j$. Therefore, Theorem 4.1 implies the equations
\begin{equation}
\int _D G^{[j]}(r, \theta ; \rho , \psi ) \varphi ^{[j]} _0
(\theta ) d\theta =
    \left\{ \begin{array}{ll}
        V^{[j]}(\rho ) W^{[j]}(r)\frac{\varphi ^{[j]} _0 (\psi )}{\chi ' _{0, j}}
            & \mbox{if $\rho < r$} \vspace{.5cm} \\
        V^{[j]}(r) W^{[j]}(\rho )\frac{\varphi ^{[j]} _0 (\psi )}{\chi ' _{0, j}}
            & \mbox{if $r < \rho $. }
                                \end{array}
    \right.
\end{equation}
We split $P_j$ in two terms, $P_j (r) = P_{j,1}(r) + P_{j,2}(r)$,
where
\[P_{j,1}(r) = \int _{K^D _r} \int _D G^{[j]}(r, \theta ; \rho , \psi )
\varphi ^{[j]} _0 (\theta ) d\theta d\mu (\rho , \psi ) \]
and
\[P_{j,2}(r) = \int _{K^D _{r, \infty }} \int _D G^{[j]}(r, \theta ; \rho , \psi )
\varphi ^{[j]} _0 (\theta ) d\theta d\mu (\rho , \psi ). \]

The equations (5.15) give us
\[P_{j,1}(r) = \int _{K^D _r} V^{[j]}(\rho )
\varphi ^{[j]} _0 (\psi ) \left(\chi ' _{0,j} \right) ^{-1} d\mu
(\rho , \psi )\] and
\[P_{j,2}(r) = \int _{K^D _{r, \infty }}W^{[j]}(\rho ) \varphi ^{[j]} _0 (\psi )
\left(\chi ' _{0,j} \right) ^{-1} d\mu (\rho , \psi ). \]

It follows from the above-mentioned properties of $\chi ' _j , \;
V^{[j]}, \; W^{[j]}$, and $\varphi ^{[j]} _0$, that
\[P_{j,1}(r) \leq b \int _{K^D _r} V(\rho )
\varphi _0 (\psi ) d\mu (\rho , \psi )\]
and
\[P_{j,2}(r) \leq b \int _{K^D _{r, \infty }} W(\rho )
\varphi _0 (\psi ) d\mu (\rho , \psi ),\] where the constants $b$
do not depend on $j$. Taking into account (5.4)-(5.5), we deduce
that $P_j (r) = P_{j,1}(r) + P_{j,2}(r) \leq b < \infty$ uniformly
with respect to $j$. So, the $q-$superfunction $U(x) = \lim _{j
\rightarrow \infty } U_j (x) < \infty $ and we can set $u(x)=
-U(x)$, that concludes the proof of Theorem 5.3. \hfill \textbf{$\diamondsuit$}\\

\noindent \textbf{Remark 5.4.} We stated Theorem 5.3 for
subfunctions with respect to radial potentials. Nonetheless,  both
its parts can be extended to general potentials, however, in
"opposite directions". Indeed, the subfunction $u$, we have just
constructed, is negative, $u \leq 0$ and $-L_q u \geq 0$, thus
this subfunction $u \in \verb"SbH"(c, K^D)$ for any potential
$c(x) \geq q(|x|)$. So, we constructed an upper bounded
subfunction with the given Riesz measure with respect to an
arbitrary, nor necessarily radial potential.

On the other hand, Theorem 5.2 implies that if $u\in \verb"SbH"(c,
K^D)$ and has an associated measure $d\mu$, then
\[\int _{K^D _b} \widetilde{V} (r) \varphi _0(\theta ) d\mu (r, \theta
) < \infty \]
and
\[\int _{K^D _{b, \infty}} \widetilde{W} (r) \varphi _0(\theta ) d\mu (r, \theta
) < \infty ,\] where $\widetilde{V}, \widetilde{W}$ are solutions
of the equation (3.4) with $\lambda =\lambda _0$ and any potential
$Q(r) \in L^p _{loc}, \; p>n$, such that $c(x) \leq Q(|x|)$.
\hfill \textbf{$\diamondsuit$}

\newpage

\section{Generalization of the Hayman - Azarin \\
theorem}

Consider a subharmonic function in a cone $K^D \subset
\mathbf{R}^n, \; n \geq 2,$ satisfying the Lindel\"{o}f boundary
condition
\[\limsup _{K^D \ni x \rightarrow x_0} u(x) \leq 0, \;
\forall x_0 \in \partial K^D. \] The Deny-Lelong theorem asserts
(see Remark 3.1) that if the boundary is not a polar set and
\begin{equation}
\liminf_{r \rightarrow \infty } r^{- \mu ^+} M(r,u) \leq 0,
\end{equation}
then $u \leq 0$ in $K^D$. A natural question arises what can be
concluded if (6.1) is replaced by
\begin{equation}
\limsup _{r \rightarrow \infty } r^{- \mu ^+} M(r,u) < \infty .
\end{equation}
Extending preceding results by Ahlfors and Heins \cite{AhlH},
Hayman \cite{Hay}, and Ushakova \cite{Ush}, Azarin \cite{Azar}
proved that if a subharmonic function $u(x)$ in a cone with the
sufficiently smooth boundary satisfies (6.1)-(6.2), then there
exists the limit
\begin{equation}
\lim _{r \rightarrow \infty } r^{- \mu ^+} u(r,\theta ) = \kappa
\varphi _0 (\theta ), \; \kappa = const,
\end{equation}
uniformly in $\theta \in \overline{D}$ provided that a point $x =
(r, \theta )$ tends to infinity outside a set of balls $\left\{B_j
(x_j, r_j ) \right\}^{\infty } _{j = 1}$ of \emph{ finite view}.
The latter means that $\sum ^{\infty
}_{j=1}\left(r_j/|x_j|\right)^{n-1}<\infty$. \\

The equation (6.3) implies, in particular, that $u$ can not decay
too fast outside of some small exceptional set. We prove here an
analog of this assertion for the subfunctions of the operator
$L_q$ with a radial potential $q(r)$. We consider only analogs of
the volume potentials, that is, we study the superfunctions
\begin{equation}
U(x) = \int _{K^D} G(x;y) d\mu (y),
\end{equation}
where $G(x;y)$ is Green's function of the Dirichlet problem for
the operator $L_q$ in $K^D$ and $d\mu$ is a given Borel measure in
$K^D$. To facilitate convergence of the integral (6.4), we assume
that there exists a positive radial minorant $q$ of the potential
$c, \; 0 \leq q(|x|)\leq c(x), \; q \in L^p _{loc} (0, \infty ),
\; p>n$, such that the measure $d\mu$ satisfies the conditions
(5.4)-(5.5). We use the approach of Azarin \cite{Azar} and the
method of normal points of Hayman \cite{Hay} as it was applied by
Govorov \cite[p. 7]{Gov}. To state our result, we need two
definitions.

\newtheorem{deff}{Definition}[section]
\begin{deff}\hspace{-.09in}\textbf{.} A set $\Omega \subset \mathbf{R}^n$
is said to be a \emph{ set of $q-$finite view} if and only if it
can be covered by the union of balls $\left\{B_j (x_j, r_j )
\right\}^{\infty } _{j = 1}$ such that
\[\sum ^{\infty } _{j = 1} \frac{r_j}{|x_j|} V\left(\frac{|x_j | }{r_j}\right)
W\left(\frac{|x_j |}{r_j}\right) < \infty .\] \hfill
\textbf{$\diamondsuit$}
\end{deff}

\noindent \textbf{Remark 6.1.} It follows from (C.1) that a set of
finite view always has a $q-$finite view, but not conversely, that
is, the sets of $q-$finite view are in general more massive than
those of finite view. For instance, let $|x_j | = j, \;
j=1,2,\ldots$, and $r_j = j^{(n-2)/(n-1)}$, thus, the series $\sum
^{\infty}_{j=1}(r_j / |x_j |)^{n-1} = \sum ^{\infty } _{j=1} 1/j $
diverges. However, if $q(r) = c_0 = const >0$, then (C.11)-(C.12)
imply that
\[\sum ^{\infty} _{j=1} \frac{r_j}{|x_j |} V\left(\frac{|x_j | }{r_j}\right)
W\left(\frac{|x_j |}{r_j}\right) \sim \sum ^{\infty } _{j=1}
j^{-n/(n-1)} < \infty .\]  \hfill \textbf{$\diamondsuit$}

\begin{deff}\hspace{-.09in}\textbf{.} A point $x$ is called
$(\varepsilon , q)-$\emph{ normal}, $\varepsilon >0$, with respect
to a Borel measure $d \mu $ if and only if for each $\lambda , \;
0 < \lambda < (1/2)|x|$, there holds the inequality
\begin{equation}
\int _{B(x, \lambda )} d\mu (y) < \varepsilon \frac{\lambda }{|x|} V
\left(\frac{|x|}{\lambda }\right) W\left(\frac{|x|}{\lambda }\right).
\end{equation}
\hfill \textbf{$\diamondsuit$} \end{deff}

Now we state the main result of this section.

\begin{theor}\hspace{-.09in}\textbf{.} If a Borel measure $d\mu $ satisfies
the conditions (5.4)-(5.5) in a cone with $C^2$-boundary, then the
function (6.4) admits an asymptotic estimate
\begin{equation}
U(r, \theta ) \leq b(\delta) \; \varepsilon \;
\frac{r^{2-n}}{W(\delta r)}, \; r\rightarrow \infty ,
\end{equation}
for any $\varepsilon >0$ and any $\delta >1$ outside a set $\Gamma
$ of $\; q-$finite view, where $b(\delta) >0$ is a constant. In
particular, (6.6) means that the $q-$subfunction $u(x) = -U(x)$
cannot decay too fast outside of an exceptional set. \hfill
\textbf{$\diamondsuit$}
\end{theor}

We use a measure (Cf. \cite{Azar})
\[dm(\rho , \psi ) = \left\{ \begin{array}{ll}
        W(\rho )\varphi _0 (\psi )d\mu (\rho ,\psi )
                &\mbox{if $\; \rho \geq 1$ } \vspace{.5cm}\\
        V(\rho )\varphi _0 (\psi )d\mu (\rho ,\psi )
                &\mbox{if $\; 0<\rho <1$. }\\
                             \end{array}
                     \right. \]
Observe that, due to (5.4)-(5.5),  $\; m(K^D)=\int
_{K^D}dm(y)<\infty$.

First, we estimate the kernel
\[K(r,\theta ; \rho ,\psi ) \stackrel{\rm def}{=}\left\{\begin{array}{ll}
   \frac{G(r,\theta ;\rho ,\psi )}{W(\rho )\varphi _0 (\psi )} &
         \mbox{if $\; \rho \geq 1$} \vspace{.5cm}\\
   \frac{G(r,\theta ;\rho ,\psi )}{V(\rho )\varphi _0 (\psi )} &
         \mbox{if $\; 0 <\rho < 1$.}
                                   \end{array}
                            \right. \]

\begin{lem}\hspace{-.09in}\textbf{.} If  $x = (r, \theta )$ is an
$(\varepsilon , q)-$normal point, then for every $\delta >1$ there
exists a constant $b(\delta ) > 0$, which does not depend on
$\varepsilon $, such that for all $x$ with $r= |x| > 2$
\[\int _{K^D _{{\delta}^{-1}r,\; \delta r}}K(r,\theta ; \rho ,\psi)
dm(\rho , \psi ) < b(\delta) \; \varepsilon \; V(r). \]
\end{lem}

\noindent\textbf{Proof.} Fix $x$ and consider truncated cones
$\Omega _{\delta } = K^D _{{\delta}^{-1} r, \; \delta r}$. Let
$\widetilde {\Omega }$ be a domain with the smooth boundary of
bounded curvature, such that
\[\Omega _{\delta +\eta}\subset \widetilde {\Omega}\subset
\Omega_{\delta +2\eta}\] with a sufficiently small $\eta >0$, and
$\widetilde {G}$ be its Green's function extended outside
$\widetilde {\Omega }$ by 0. The difference $\Phi
(y)=G(x;y)-\widetilde {G}(x;y)$ is a $q-$harmonic function in
$\Omega _{\delta }$ which vanishes at $\partial K^D$ and, due to
Theorem 4.2, satisfies the inequality
\[\Phi (y) \leq b_1 (\delta )V(r)W(\rho )\varphi _0 (\theta )\varphi _0(\psi )\]
at all other boundary points of $\Omega _{\delta }$. By the
maximum principle, the same inequality holds everywhere in $\Omega
_{\delta }$, that is,
\[G(x;y) \leq b_1 (\delta )V(r)W(\rho )\varphi _0 (\theta )\varphi
_0(\psi )+ \widetilde {G}(x;y).\] If $n\geq 3$ (for $n=2$
calculations are similar), we have from this estimate
\[G(x;y) \leq b_1 (\delta )V(r)W(\rho )\varphi _0 (\theta )\varphi _0(\psi ) +
\theta _n |x-y| ^{2-n},\] and so
\[K(r, \theta ; \rho , \psi ) \leq b_1 (\delta )V(r)\varphi _0 (\theta )+
\theta _n \frac{|x-y| ^{2-n}}{W(\rho )\varphi _0(\psi)}.\]

On the other hand, $\widetilde {G}(x;y) \leq g(x;y)$, where as
before, $g$ is the harmonic Green function, admitting the
following known bound \cite[(3.2)]{Azar},
\[g(x;y) \leq \rho ^{2-n}\varphi _0 (\theta )\varphi _0(\psi )
\left( b_2 (\delta )+ b_3 (\delta )r \rho ^{n-1} |x-y|^{-n}
\right), \; y \in \Omega _{\delta}.\] Thus,
\[G(x;y) \leq V(r)W(\rho )\varphi _0 (\theta )\varphi _0(\psi )
\left(b_1 (\delta )+ \frac{b_2 (\delta )\rho ^{2-n}}{V(r)W(\rho )}
+ \frac{b_3 (\delta )r \rho |x-y|^{-n}}{V(r)W(\rho )} \right), \]
and we get yet another bound of the kernel $K(x;y)$:
\[K(r, \theta ; \rho , \psi ) \leq V(r)\varphi _0 (\theta )
\left(b_1 (\delta )+ \frac{b_2 (\delta)\rho ^{2-n} + b_3 (\delta
)r \rho |x-y|^{-n}}{V(r)W(\rho )} \right). \]

Integrating these inequalities over $\Omega _{\delta}$, we obtain
the estimate
\begin{equation} \begin{array}{cc}
\int_{\Omega_{\delta}}K(r,\theta ;\rho ,\psi)dm(\rho ,\psi) \vspace{.3cm}\\
\leq b_1 (\delta )V(r)\varphi _0 (\theta )\int _{\Omega
_{\delta}}dm(y) + \int _{\Omega _{\delta }} \widetilde {K}(r,
\theta ; \rho , \psi )dm(y), \; y=(\rho , \psi ),
\end{array}  \end{equation}
where
\[\widetilde {K}(x;y)= \min \left\{\frac{\theta _n |x-y| ^{2-n}}{W(\rho) \varphi _0(\psi )};
\; \frac{\varphi _0 (\theta ) \left(b_2(\delta)\rho ^{2-n} +
b_3(\delta )r \rho |x-y|^{-n}\right)}{W(\rho)}\right\}.\]

Since $m(K^D)<\infty $, there holds $m(\Omega _{\delta}) =
\bar{\bar{o}}(1), \; r\rightarrow \infty $, and so the first term
on the right in (6.7) is $\bar{\bar{o}}(V(r)), \; r\rightarrow
\infty $. To estimate the second integral on the right-hand side
of (6.7), we split $\Omega _{\delta}$ into two parts:
\[\Omega _{\delta} ' = \left\{y=(\rho ,\psi ) \in \Omega _{\delta }
\; \biggl| \; \; \sin \widehat {(\theta ,\psi )} < \frac{1}{4}
\right\}\] and $\Omega _{\delta }'' = \Omega _{\delta } \setminus
\Omega _{\delta } '$. If $y \in \Omega _{\delta }''$, then $|x-y|
\geq \frac{1}{4}\rho $ and $\varphi _0 (\psi ) \rightarrow 0, \;
\psi \rightarrow \partial D$, therefore, $\widetilde {K}(r, \theta
; \rho , \psi ) \leq b(\delta ) \frac{|x-y| ^{2-n}}{W(\rho)
\varphi _0(\psi )}$. Thus, (6.5) with $\lambda = r b(\delta )$
gives
\[\int _{\Omega _{\delta}''} \widetilde {K}(r, \theta ; y)dm(y)
\leq b(\delta ) \; \varepsilon \; r^{2-n}, \; r\rightarrow \infty
.\]

To bound the integral over $\Omega _{\delta }'$, we denote $a=
a(x)= dist \left\{x; \; \partial K^D \right\}$ and decompose
$\Omega _{\delta }'$ into spherical layers
\[\Upsilon _m = \left\{y \in \Omega _{\delta }' \; \biggl| \;\;
2^{m-1}a \leq |x-y|< 2^m a, \; -\infty <m<\infty \right\}.\]
Because $x$ is an $(\varepsilon ,q)-$normal point, we have $\mu
(x) = m(x) = 0$ and
\[\int _{\Omega _{\delta}'}\widetilde {K}(x;y)dm(y) = \sum ^{\infty } _{m= -\infty }
\int _{\Upsilon _m }\widetilde {K}(x;y)dm(y).\] Notice also that
by the definition of $\Upsilon _m$, $|x-y| ^{n-2} \geq
2^{(m-1)(n-2)} a^{n-2}$ and (Cf. \cite{Azar}) $r \varphi _0 (\psi
) \geq \mbox{const} \cdot a$. Moreover, we always have
\[\widetilde {K}(x; \rho , \psi ) \leq \frac{\theta _n |x-y|
^{2-n}} {W(\rho) \varphi _0(\psi )}.\]

First, let be $m<0$. If $\rho \geq \rho _0$, then obviously,
$2^{2m}a \leq \rho $, so, (6.5) with $\lambda =2^m a$ implies
\[\int _{\Upsilon _m } \widetilde {K}(x;y)dm(y)= \theta _n
\int _{\Upsilon _m } |x-y| ^{2-n} dm(y) \leq \]
\[\leq \theta _n a^{2-n} 2^{(m-1)(n-2)} \int _{\Upsilon _m }
dm(y)\]
\[\leq \theta _n a^{2-n} 2^{(m-1)(n-2)} \; \varepsilon \;
\frac{2^m a}{r} V\left(\frac{r}{2^m a}\right) W\left(\frac{r}{2^m
a}\right).\]

It is easy to draw from here and (C.1) that
\[\sum ^{-1} _{m= -\infty }\int _{\Upsilon _m } \widetilde {K}(x;y)dm(y)
\leq b(\delta)\varepsilon r^{2-n} \sum ^{\infty } _{k=1} 2^{-k} =
b(\delta )\varepsilon r^{2-n} .\]

Now, let be $m \geq 0$. In this case
\[\widetilde {K}(x; \rho ,\psi ) \leq
\frac{b(\delta ) \rho ^2 }{|x-y|^n W(\rho )\varphi _0 (\psi )}
\leq \frac{2^n b(\delta ) \rho ^2 }{2^{nm} a^m W(\rho )\varphi _0
(\psi )},\] and
\[\int _{\Upsilon _m } \widetilde {K}(x;y)dm(y)
\leq \frac{b(\delta ) r^2 }{2^{nm} a^m } \mu (\Upsilon _m )\leq
b(\delta ) \; \varepsilon \; a \; r^{1-n},\] where $b(\delta )$
does not depend on $m$ and as before, we set $\lambda =2^m a$ in
(6.5). From the definition of $\Omega _{\delta }$ we see that it
can contain no more than $b(\delta )r/a$ of the layers $\Upsilon
_m$ with $m \geq 0$, therefore, $\sum ^{\infty } _{m=0} \int
_{\Upsilon _m } \widetilde {K}(x;y)dm(y)\leq b(\delta ) \;
\varepsilon \; r^{2-n}$.

To finish the proof, it is enough now to combine all the drawn
bounds, take into consideration (C.1), and let $\eta \rightarrow
0$. \hfill \textbf{$\diamondsuit$}\\
\begin{lem}\hspace{-.09in}\textbf{.} The set $\Delta (\varepsilon , q, \mu )$
of points, which are not $(\varepsilon , q)-$normal, has
$q-$finite view.
\end{lem}
\textbf{Proof} mimics the proof of Lemma 7 in \cite{Azar}. \hfill \textbf{$\diamondsuit$}\\

\noindent\textbf{Proof of Theorem 6.1.} We choose an $\eta >0$ and
a $\rho _{\eta } >0$ such that $\int _{K^D _{\rho _{\eta} , \infty
}} W(\rho )\varphi _0 (\psi )d\mu (\rho , \psi ) < \eta$. If $r_0
> \rho _{\eta}$, then by Theorem 4.2,
\[\int _{K^D _{\rho _{\eta}}} G(r, \theta ; y)d\mu (y)
\leq b(r_0)W(r) \varphi _0 (\theta ) \int _{K^D _{\rho _{\eta}}}
V(\rho ) \varphi _0 (\psi )d\mu (\rho , \psi ),\] that is, the
integral over $K^D _{\rho _{\eta}}$ tends to zero as $r\rightarrow
\infty $. Next, we break up the remaining integral over $K^D
_{\rho _{\eta}, \; \infty }$  into three parts,

\[U_1(x) = \int _{K^D _{\rho _{\eta}, \; r/\delta }} G(x;y)d\mu (y)\]
\[U_2(x) = \int _{\Omega _{\delta }} G(x;y)d\mu (y)\]
and
\[U_3(x) = \int _{K^D _{\delta r, \; \infty }} G(x;y)d\mu (y).\]

For estimating $U_1$ we use (4.19). Thus,
\[U_1(x)\leq b(\delta )\frac{r^{2-n}}{W(\delta r)}\int _{K^D _{\rho _{\eta },
\; r/\delta }}W(\rho ) \varphi _0 (\psi )d\mu (\rho , \psi ) \leq
\frac{b(\delta ) \; \eta \; r^{2-n}}{W(\delta r)}.\] An estimation
of $U_2$ is given by Lemma 6.1. It should be noted that due to
(C.1), $V(r) < V(\delta r) \leq b(\delta ) \frac{r^{2-n}}{W(\delta
r)}$, that is,
\[U_2(x) \leq \frac{\varepsilon \; b(\delta ) \; r^{2-n}}{W(\delta r)}.\]

Finally, to estimate $U_3$, we again use (4.19) with $\rho >
\delta r$:
\[U_3(x) \leq \frac{b(\delta ) \; r^{2-n}}{W(\delta r)} \int _{K^D _{\delta r, \; \infty }}
W(\rho ) \varphi _0 (\psi )d\mu (\rho , \psi ) \leq \frac{\eta \;
b(\delta)\; r^{2-n}}{W(\delta r)}.\] \hfill
\textbf{$\diamondsuit$}

\noindent \textbf{Remark 6.2.} The function $r^{2-n}/W(r)$ is
monotonically increasing as the product of two monotone functions,
$\frac{r^{2-n}}{W(r)} = r^{\mu ^+} \left\{r^{\mu ^-}
W(r) \right\} ^{-1}$. \hfill \textbf{$\diamondsuit$}\\

\begin{cor}\hspace{-.09in}\textbf{.} In the case \textbf{(A)},
$W(r) \approx r^{\mu ^- _k}$, where $\mu ^- _k$ is the negative
root of the equation $\mu (\mu +n-2) = \lambda _0 +k$, and the
estimate (6.6) is equivalent to the inequality $U(x)\leq b \;
\varepsilon |x|^{\mu ^+ _k}$, where $b$ depends only on $D$.
\hfill \textbf{$\diamondsuit$}
\end{cor}

\begin{cor}\hspace{-.09in}\textbf{.} In the case \textbf{(B)}, (6.6) is
equivalent to the inequality
\[U(x)\leq b(\delta) \; \varepsilon \left(s(\delta r)\right)^{1/2}
V(\delta r).\]
\end{cor}
\textbf{Proof.} It suffices to apply asymptotic formulas of Lemma
C.3. \hfill \textbf{$\diamondsuit$}\\
\begin{cor}\hspace{-.09in}\textbf{.} In the case \textbf{(B)}, if
$r \rightarrow \infty $ provided that $x\not \in \Gamma $ and $r
\not \in E$, where $E$ is a set of finite logarithmic
length\footnote{That is, $\int _E t^{-1}dt < \infty $.}, then
asymptotically
\begin{equation}
U(x) \leq b(\delta ) \; \varepsilon \; V(\delta |x|).
\end{equation}
\end{cor}
\noindent\textbf{Proof.} From the asymptotic formulas
(C.11)-(C.12) for $V(r)$ and $W(r)$ we deduce, with $k > \delta
>1$,
\[\frac{r^{2-n}}{W(\delta r)} = \frac{r^{2-n}}{W(\delta r)V(k \delta r)}V(k \delta r)\]
\[\leq b \left(s(r)s(k \delta r) \right)^{1/4} \exp \left\{ - \int _r ^{k
\delta r}\sqrt{s(t)} \; \frac{dt}{t} \right\} V(k \delta r).\] Let
us denote $T(r) = \sqrt{s(r)}$. Since $s(r)$ is a monotone
function, then
\[\frac{r^{2-n}}{W(\delta r)}\leq A T(\delta r)
\exp \left\{-T(r)\ln \delta \right\}V(k \delta r).\] By making use
of the Borel inequality \cite[p. 121]{GoOs} with $\delta =1+1/\ln
T(r)$ we obtain the estimate $T(\delta r)<
\left\{T(r)\right\}^{1+\varepsilon },\; r \not \in E$, and from
here, in turn,
\[T(\delta r) e^{- T(r)\ln \delta }<\left\{T(r)\right\}^{1+\varepsilon }
\left\{\left(1+\frac{1}{\ln T(r)} \right)^{- \ln T(r)} \right\}
^{\frac{T(r)}{\ln T(r)}} < \]
\[< \left\{T(r)\right\}^{1+\varepsilon } 2^{-\frac{T(r)}{\ln T(r)}}\rightarrow 0,
\; r\rightarrow \infty . \] Thus, the asymptotic inequality
$\frac{r^{2-n}}{W(\delta r)} \leq b(\delta ) V(k \delta r)$ holds
for all $k > \delta >1$ and $r\not\in l^{-1}E$. After changing
notations, this inequality and (6.6) imply (6.8). \hfill \textbf{$\diamondsuit$}\\

\newpage

\section{Subfunctions in tube domains}

Our methods can be applied to other unbounded regions. In this
section we consider tube domains
\[T^D = \mathbf{R}^n \times D = \left\{(x,y) \; \biggl| \; \; x\in
\mathbf{R}^n, \; n \geq 2, \; y \in D \right\},\] where $D$ is an
arbitrary bounded domain in $\mathbf{R}^p , \; p \geq 1$. If a
potential $c \in {\cal C}(T^D)$, then an operator $L_c$, generated
by the differential expression
\[L_c u = - \Delta _x u - \Delta _y u + c(x,y) u, \]
can be defined as in the first paragraph of Appendix A. In this
section $\lambda _0$ and $\varphi _0 (y)$ denote the lowest
eigenvalue and the corresponding eigenfunction of the boundary
value problem
\[\left\{ \begin{array}{ll}
    \Delta _y u + \lambda u = 0, & y \in D \vspace{.3cm}\\
    u(y) = 0, & y \in \partial D \setminus E, \\
                                \end{array}
                  \right. \]
where $E \subset \partial D$ is an exceptional polar set and $V$
and $W$ are, respectively, the growing and decaying solutions of
the equation
\begin{equation}
y'' +(n-1)r^{-1}y' - \left(\lambda _0 +q(r)\right)y=0, \; 0<r<
\infty ,
\end{equation}
with a positive, locally summable function $q,\; 0\leq q(r) \leq
c(x,y), \; \forall y\in D$ and $\forall x\in S(r), \; r=|x|$,
where $S(r)$ is a sphere of radius $r$, centered at zero, in
$\mathbf{R}^n$. They are normalized by $V(1)=W(1)=1$. We state
without proofs, which are similar to the preceding ones, analogs
of our main results from sections 3-5.

\begin{theor}\hspace{-.09in}\textbf{.} Let $c(x,y)\in {\cal C}(T^D)$ and
$u(x,y)\in \verb"SbH"(c, T^D)$. If
\[\limsup _{T^D\ni (x_0,y)\rightarrow (x_0,y_0)} u(x_0,y) \leq A\]
for each point $(x_0 , y_0) \in \partial T^D = \mathbf{R}^n \times
\partial D$ and
\begin{equation}
\liminf _{R\rightarrow \infty } \frac{1}{V(R)} \int _D {\cal
M}(y,R,u^+) \varphi _0 (y) dy = 0,
\end{equation}
where ${\cal M}(y,R,u) = \frac{1}{|S(R)|} \int _{S(R)} u(x,y)
d\sigma _x (R, \theta )$, then $u(x,y)\leq A^+$ everywhere in
$T^D$. Here, $|S(R)|$ is the surface area of the sphere $S(R)$.

If we replace $A^+$ with $A$, then the conclusion fails. Moreover,
the condition (7.2) is the best possible in the class of all
potentials $c$ dominating the function $q$.
\end{theor}
\noindent\textbf{Sketch of the proof.} We solve the Dirichlet
problem for the operator $L_q$ in $T^D _R = B(R) \times D, \;
R>0$, with the boundary conditions to be $0$ on $B(R) \times
\partial D$ and $u^+ ((R, \theta ), y)$ on $S(R) \times D$.
Denoting its solution by $\widetilde{u}$, it can be proved as in
Section 3, that the function
\[w_R (x) = \int _D \widetilde{u} (x,y) \varphi _0 (y) dy\]
satisfies the equation
\begin{equation}
\Delta _x w(x) - (\lambda _0 + q(r))w(x)= 0.
\end{equation}
So, its average
\[\gamma _R (r, \widetilde{u}) = \int _D {\cal M}(y,r, \widetilde{u}) \varphi _0 (y)
dy\] satisfies the same equation (7.3) but does not depend on the
spherical component $\theta$ of $x = (r, \theta )$. The rest of
the proof is the same as in Theorem 3.1. \hfill \textbf{$\diamondsuit$}\\

\noindent \textbf{Remark 7.1.} In particular, the conclusion holds
if we replace (7.2) with the following condition involving the
maximum modulus of $u^+$,
\[\liminf _{r\rightarrow \infty }\frac{M(r,u^+)}{V(r)} = 0.\]
\hfill \textbf{$\diamondsuit$}
\begin{cor}\hspace{-.09in}\textbf{.} If $A=0$ and $\varpi = \liminf
_{r\rightarrow \infty } \frac{1}{V(r)} \int _D {\cal M}(y,r,u^+)
\varphi _0 (y) dy $, then $\int _D {\cal M}(y,r,u^+) \varphi _0
(y) dy \leq \varpi \; V(r), \; 0<r< \infty $. \hfill
\textbf{$\diamondsuit$}
\end{cor}
\begin{cor}\hspace{-.09in}\textbf{.} If $A=0$ and
\[\liminf _{r\rightarrow \infty} \int _D \left(\frac{u(x,y)}{V(|x|)\varphi _0 (y)}
- \sigma\right)^+ \varphi ^2 _0 (y) dy =0\] with a $\sigma \geq
0$, then $u(x,y)\leq \sigma V(|x|) \varphi _0 (y)$ in $T^D$.
\hfill \textbf{$\diamondsuit$}
\end{cor}

Comparing equations (7.1) and (3.4) we see that since $q(r) \geq
0$ and $\lambda _0 > 0$, solutions of (7.1) always have asymptotic
behavior similar to the solutions of (3.4) in the case
\textbf{(B)} (Cf. Lemma C.3), whenever we assume in addition that
$\int ^{\infty} \left( \widetilde{q}(t)\right)^{1/2} dt = \infty
$, where $\widetilde{q}(t) = \lambda _0 +q(t) + (1/4)(n-1)(n-3)
r^{-2}$. Thus, if this $\widetilde{q}$ satisfies also
\begin{equation}
\int ^{\infty} \biggl| \; 4 \widetilde{q}(r) \widetilde{q}''(r) -
5 \left(\widetilde{q}'(r)\right)^2 \biggl| \;
(\widetilde{q}(r))^{-5/2} dr < \infty ,
\end{equation}
then the growing solution of (7.1) satisfies
\begin{equation}
V(r)=b r^{(1-n)/2} \left(\lambda _0 +q(r)\right)^{-1/4} \exp
\left\{\int ^r _1 \sqrt{\lambda _0 +q(t)}dt \right\}
(1+\bar{\bar{o}}(1)),\; r\rightarrow \infty .
\end{equation}

Moreover, we see that instead of the "inverse square" growth rate,
which is the borderline rate in the case of cones, in the case of
tube domains we are to consider the boundedness condition. So
that, if in addition to (7.4) there holds
\[q(r)=\underline{\underline{O}}(1), \; r\rightarrow\infty ,\]
then the lowest possible growth of subfunctions in $T^D$, which
are upper bounded at the boundary $\partial T^D$, is the
exponential growth\footnote{That is the case for analytic
functions in a strip \cite[p.72]{L1}.} that depends on $D$.
However, if
\[q(r) \nearrow \infty, \; r\rightarrow\infty ,\]
then this lowest growth is faster than any exponential.\\

We mention a few special cases. If $p=1$ and $D= \{0<y<h \}$ is an
interval, then $T^D$ is an infinite layer in $\mathbf{R}^{n+1}$.
Here, $\lambda _0 = (\pi /h)^2$, $\varphi _0 (y) = \sin (\pi y /
h)$, and the asymptotic behavior of $V(r)$ is given by (7.5) with
$\lambda _0 = (\pi /h)^2$. In particular, if $q(r)= c_0 =const$,
the condition (7.2) becomes
\[\liminf_{R\rightarrow\infty}\frac{R^{(n-1)/2}}{\exp\left\{R\sqrt{(\pi /
h)^2 +c_0 }\right\}} \int _1 ^h {\cal M}(y,R,u^+) \sin (\frac{\pi
y}{h}) dy =0.\] This case was considered in \cite{Landis56} for
operators of any order with constant coefficients. In particular,
if in addition $n=1$ and the potential $q(r)= c_0 =const \geq 0$,
then the lowest possible growth in the strip
\[\left\{(x,y) \in \mathbf{R}^2 \biggl| \; - \infty <x< \infty, \;
0<y<h \right\}\] coincides with (7.5) and is given by the
$c_0-$subfunction

\[\cosh \left(\sqrt{c_0 + \pi ^2 /h^2} \; x \right) \sin (\pi y /h).\]

As another example, let $p=2$ and $D$ be a disk of radius $\varrho
$ in the complex plane, that is, $T^D$ is a bilateral
cylinder\footnote{Similar results can be proved for unilateral
cylinders.} in $\mathbf{R}^{n+2}$. In this case $\lambda _0 =
\gamma _0 / \varrho $ and $\varphi _0 (y) = J_0 (\frac{\gamma _0
}{\varrho } r)$, where $J_0$ is the Bessel function of first kind
and $\gamma _0$ is its first positive root. Again, $V(r)$ is a
solution of (7.1) with $\lambda _0 = \gamma _0 / \varrho $. Thus,
if $q(r)= c_0 =const \geq 0$, then the lowest growth is
$r^{(1-n)/2} \exp \{r \; \sqrt{c_0 + \lambda _0}\}$.\\

In \cite{ArEvg, Yarm} and references therein the authors studied
asymptotic behavior of the classical harmonic functions in a
cylinder, when $p=2,\; n=1$, and showed that if $u\bigl| _{y\in
\partial D} \leq b$ and $\frac{\partial u}{\partial n} \bigl|
_{y\in \partial D} \leq b$, then the lowest possible growth of
such a function is $\exp (\exp (a x))$ with some $a>0$. It should
be mentioned that for subfunctions the condition $\frac{\partial
u}{\partial n} \bigl| _{y\in \partial D} \leq b$ does not play an
essential role. Indeed, if the boundary $\partial D$ is
sufficiently smooth, then the $q-$subfunction $u(x,y)=
\left(V(|x|) \varphi _0 (y)\right)^{1+\varepsilon}, \; \varepsilon
>0$, satisfies $u\bigl| _{\partial D} = \frac{\partial u}{\partial
n} \bigl| _{\partial D} =0$, however, the growth of $V(r)$ is again
given by (7.5).\\

Finally, we state a generalization of the Blaschke theorem.
Carleman's formula for the smooth tube domain $T^D _{1,R}=
\left\{x\in \mathbf{R}^n :\; 1<|x|<R\right\}\times D$ is similar
to (5.3):
\[\int _{T^D _{1,R}} \left( W(|x|)- \frac{W(R)}{V(R)}
V(|x|)\right) \varphi _0 (y) d\mu (x,y) \]

\[+ \int _{\{1<|x|<R\}\times \partial D} \left(W(|x|)- \frac{W(R)}{V(R)}
V(|x|)\right) \frac{\partial \varphi _0 (y)}{|x|\partial n(y)}
u(x,y) d\sigma (x,y) \]

\[= \frac{\chi}{V(R)}\int _{S(R)\times D} u(R,y)\varphi _0 (y) dy +
\underline{\underline{O}}(1), \; R\rightarrow \infty ,\] and the
second formula alike.
\begin{theor}\hspace{-.09in}\textbf{.} Let $d\mu (x,y)$ be a Borel measure
in $T^D$. In order for an upper bounded subfunction in $T^D$ with
the Riesz associated measure $d\mu $ to exist, it is necessary and
sufficient that
\[\int _{T^D _R} V(|x|)\varphi _0(y) d\mu (x,y)< \infty \]
and
\[\int _{T^D \setminus T^D _R} W(|x|)\varphi _0(y) d\mu (x,y)< \infty .\]
This function is given by
\[u(x,y)=-\int _{T^D} G(x,y; x',y')d\mu (x',y').\]
\hfill \textbf{$\diamondsuit$}
\end{theor}

\newpage
\appendix
\section*{Appendices}
\addcontentsline{toc}{section}{\hspace{.65cm}Appendices}
\section{Green's function of the operator $L_c$}

Let $\Omega$ be an arbitrary domain in $\mathbf{R}^n, \; n \geq
2$. The differential expression (1.1), defined initially  on
smooth finitely supported functions in $ \Omega $, can be extended
in a standard way to an essentially self-adjoint operator $ L_{c}
$ on $ L^{2}( \Omega ) $. A function $ G(x;y) $ defined for $
(x,y) \in \Omega \times \Omega, \; x \not = y $, is called the
Green function of the extended operator $L_{c}$ if $G$ gives a
representation $u(x)= \int_{\Omega} G(x;y)f(y) \; dy $ of a weak
solution of the equation
\begin{equation}
L_{c}u(x) = f(x), \; x \in \Omega, \; f \in L^{2}(\Omega) ,
\end{equation}
satisfying the boundary condition $ u(x)=0 $ at every point $ x
\in \partial \Omega \setminus E$, where $E$ is an exceptional,
maybe empty, boundary set of zero capacity. In other words,
Green's function $G$ generates an inverse operator to (A.1). The
Green function satisfies the equation $L_{c}G(x;y) =\delta (x-y)$,
$\delta$ being the Dirac $\delta-$function, in the domain $\Omega$
 in the sense of distributions and takes on zero boundary values at
$\partial \Omega \setminus E$. Hereafter, terms \emph{ capacity},
\emph{ polar set}, \emph{ regular} or \emph{ irregular point},
etc., are used in the sense of the classical potential theory -
see, for example, \cite{HK}.

The existence of Green's function for the operator (A.1) with
measurable coefficients was established in \cite{Her2} under
stronger assumptions\footnote{For so-called $L_c-$adapted domains
\cite[p. 309]{Her2}.} on the domain $D$. It is known (see, for
example, \cite{Sim} or \cite{Bram}) that the Green function $G$ of
the operator $L_c$ exists for more general potentials than in the
class\footnote{This class of potentials was defined in Section 2,
after (2.1)-(2.2).} ${\cal C}$, but the assumptions (2.1)-(2.2)
allow us to prove in Theorem A.1 that the Green function of $L_c,
\; c\in {\cal C}$, exists for any domain whose boundary is not a
polar set. Moreover, as we have already mentioned, no restriction
on the growth of the potential at the boundary or at infinity is
assumed. In addition, in this theorem we also derive certain
analytic properties of $G$. \\

These results are mostly known\footnote{For a concise survey of
relevant results see, for example, \cite{Anc}.}, we just present
them in a form suitable for our goals. We use the approach of Levy
reducing (A.1) to an integral equation. Due to (2.1), the operator
$L_c$ is positive, that is, $(u, \; L_c u) \geq 0$. Therefore, the
point $\; -1$ does not belong to its spectrum and the
corresponding integral equation is uniquely solvable.

\begin{theor}\hspace{-.09in}\textbf{.} Let $\Omega \subset \mathbf{R}^n, \; n \geq 2 $,
be an arbitrary domain, whose boundary $\partial \Omega $ is not a
polar set, and $c \in {\cal C}(\Omega)$. Then \\

$1^{o}$ There exists, and unique, Green's function $G(x;y)$ of the
operator $L_{c}$ with the Dirichlet boundary conditions in $\Omega$. \\

$2^{o}$ This function is symmetric, that is, $G(x;y) = G(y;x)$. \\

$3^{o}$ At its singularity, as $r=|x-y|\rightarrow 0$, $G$ is
normalized as the harmonic Green's function $g$:
\[ G(x;y) = \left\{ \begin{array}{ll}
        \theta _2 \ln \frac{1}{r}+ \bar{\bar{o}}(\ln \frac{1}{r})
                & \mbox{if $ n=2$}\vspace{.5cm} \\
        \theta _n r^{2-n} + \bar{\bar{o}}(r^{2-n}) & \mbox{if $n \geq 3.$ }
                                \end{array}
                  \right. \]
The constant $\theta _{2} = 2 \pi $ and $\theta _{n} = (n-2)
\sigma _{n-1}$ for $n \geq 3$, where $\sigma _{n-1}$ is the
surface area
of the unit sphere in $\mathbf{R}^n.$ \\

$4^{o}$ If $ (x, y) \in (( \partial \Omega \setminus E ) \times
\Omega ) \cup (\Omega \times (\partial \Omega \setminus E )) $,
where $ E \subset \partial \Omega $ is the set of irregular
boundary points,\footnote{$E$ is a polar set \cite[p. 242]{HK}.}
then $G(x;y) = 0$. \\

$5^{o}$ $G(x;y) > 0$ everywhere in $\Omega \times \Omega $. \\

$6^{o}$ For any fixed point $y_{0} \in \Omega $, Green's function
$G(x;y_{0})$ is continuous for $x \in \Omega \setminus \{ y_{0}
\}$. Moreover, due to (2.2) $G(x;y_{0})$ is a
H\"{o}lder-continuous function with the index $\beta =\min \{1,\;
2-n/p\}$ on each compact set $K\subset\Omega\setminus \{y_{0} \}$.
If $p\geq n$, then $G$ is Lipschitz-continuous (that is,
H\"{o}lder-continuous with the index $\beta = 1$) on any
$K \subset \Omega \setminus \{y_{0} \}$. \\

$7^{o}$ $G(x;y)$ has the generalized first partial derivatives
belonging to the class $L^{s} _{loc} (\Omega ), \; \forall s \in
[1, \frac{n}{n-1})$. However, since we assume $c \in {\cal
C}(\Omega)$, the difference $G_1 (x;y) \stackrel{def}{=} G(x;y) -
\theta _n |x-y|^{2-n}$ has better differential properties. Namely,
if $n \geq 3$, then the first derivatives of $G_1$
belong\footnote{The inclusion $G(x;y)\in H_0^{1,q}(\Omega)$ with
$q<\frac{n}{n-1}$ was proved by Stampaccia \cite{Sta} under the
condition $c(x)\geq c_0 = const >0$ and by Maeda \cite{Maed}
without this restriction. Let us remark that
$\frac{n}{n-1}<\frac{np}{n+p(n-3)}$ if $n/2<p$, unless $n=2$ and
$p>2$. Hereafter, $H^{1,q}$ and $H_0^{1,q}$ are the usual Sobolev
spaces.} to $L^{s}( \Omega _1), \; \forall s \in \left[1,
\frac{np}{n+p(n-3)}\right)$, and its generalized second
derivatives belong to $L^{s}(\Omega _1)$, $\forall s \in \left[1,
\frac{np}{n+p(n-2)} \right)$, where $\Omega _1$ is a compactly
imbedded subdomain of $\Omega$. Moreover, if $c \in L^{p}
_{loc}(\Omega )$ with $p > n$, then these first derivatives are
continuous for $x \neq y$, and if $c \in L^{p} _{loc}(\Omega )$
with $n/2 < p \leq n$, they belong to $L^{s} (\Omega _1)$,
$\forall s \in [1, \; \frac{np}{n - p})$; let us note that
$\frac{np}{n+p(n-3)} < \frac{np}{n-p}$. Here $\Omega _1$ is any
domain such that $\Omega _1 \subset \subset \Omega \setminus
\{y\}$ or $\Omega _1 \subset \subset \Omega \setminus \{x\}$. \\
\par If $n=2$ (we assume in this case $p=2$), the first derivatives of $G_1$
belong to $L^{2 - \epsilon}( \Omega _1)$ with any $\epsilon > 0$. \\

$8^{o}$ If $E \subset \partial \Omega$ is a sufficiently smooth
part of the boundary (for example, it suffices for $E$ to be a
$C^1$ surface), then the inward normal derivative
\[\frac{\partial G(x_0;y)}{\partial n(y)}=\lim _{x \rightarrow x_0}
\frac{G(x;y)}{|x- x_0|} \geq 0 \] exists almost everywhere on $E$
and is locally-summable and nonnegative. Here, $x \in \Omega , \;
x_0 \in \partial \Omega $, and $x \rightarrow x_0$ along the
inward normal to $\partial \Omega $ at $x_0$. If the boundary
$\partial \Omega $ consists of a finite number of $C^1$-smooth
parts and $n/2 < p \leq n$, then $\partial G / \partial n \in
L^{s}(E_1)$, $\forall s \in [1, \frac{(n-1)p}{n-p})$, where $E_1$
is any compact part of the boundary $\partial \Omega $. In
addition, if a continuous normal vector exists everywhere on
$\partial \Omega $ and $p > n$, then $\partial G / \partial n$ is
continuous on $\partial \Omega $. In general case, if there exists
a normal vector at a boundary point $x_0 \in \partial \Omega $,
then
\[\frac{\partial G(x_0;y)}{\partial n(y)} \stackrel{\rm def}{=}
\liminf _{\Omega \ni x \stackrel{n}{\rightarrow} x_0}
\frac{G(x;y)}{| x- x_0 | } \geq 0. \]
\end{theor}
\noindent\textbf{Proof.} It suffices to prove the theorem for a
bounded domain $\Omega $. Indeed, if $\Omega $ is not bounded, we
exhaust it from within with an increasing sequence of smooth
bounded domains, whose Green's functions make up an increasing
sequence bounded from above in $\Omega $ by the generalized
harmonic Green's function $g$. The harmonic Green function $g$ in
$\Omega$ exists \cite[p. 250]{HK} due to our assumption that
$\partial\Omega$ is not a polar set and the inequality $G\leq g$,
which is known for smooth domains - see (A.14) below. So, in the
following proof we assume $\Omega $ to be a bounded domain and
moreover, in proving local statements we assume $c \in L^p(\Omega
)$. First, we estimate in the following lemma the iterated kernels
\begin{equation} g^{[k+1]}(x;y)=\int _{\Omega}g^{[k]}(x;t)c(t)g(t;y)dt,
\; k = 1, 2, ..., \; g^{[1]}(x;y)= g(x;y). \end{equation}

\begin{lem}\hspace{-.09in}\textbf{.} Let $\; \Omega$ be a bounded domain in
$\mathbf{R}^n , \; n \geq 2$. Denote $\gamma = \frac{2p-n}{p}$,
$\gamma > 0$ for $p>n/2$. If $c \in {\cal C} (\Omega) \cap L^p
(\Omega )$, then
\begin{equation}
g^{[k+1]}(x;y)\leq \left\{ \begin{array}{ll}
   b \left(b_1 \|c\| \right)^k |x - y|^{2-n+k\gamma } &
          \mbox{if  $\; n-2-k\gamma > 0$ } \vspace{.5cm} \\
     b \left(b_1 \|c\| \right)^k  & \mbox{if  $\; n-2-k\gamma < 0$.}
                                \end{array}
                  \right.
\end{equation}
\end{lem}

\noindent \textbf{Proof.} To estimate $g^{[2]}$, we make use of
the next basic property of the kernels with a weak singularity
\cite[p. 59]{Mikh}:\\

\emph{ If $0<\alpha <n,\; 0<\beta <n$, and $\alpha +\beta >n$,
then}
\begin{equation}
\int _{\Omega } \mid x-t \mid ^{- \alpha } \mid t-y \mid ^{- \beta
} dt \leq b |x-y|^{n- \alpha - \beta } .
\end{equation}
If $n \geq 3$, then $g(x;y)\leq b |x-y|^{2-n}$. Applying the
H\"{o}lder inequality to $g^{[2]}$, we get (hereafter $1/p +1/q
=1$)
\[g^{[2]}(x;y)\leq \|c\|_{L^p(\Omega)}\left\{\int_{\Omega}
\left(g(x;t)g(t;y)\right)^q dt\right\}^{1/q}\]
\[\leq b \|c\|_{L^p (\Omega )}\left\{\int_{\Omega }\left(|x-t||t-y|\right)^{q(2-n)}dt
\right\}^{1/q}.\] Now (A.4) with $\alpha = \beta = q(n-2)$ gives
\[g^{[2]}(x;y) \leq b \|c\|_{L^p (\Omega )} |x-y| ^{2-n+ \gamma}.\]
Iterating this estimate, we prove (A.3) for $n \geq 3$. If $n=2$,
we have to replace here $|x-y|^{2-n}$ with $\ln \frac{1}{|x-y|}$
and repeat the same calculations; in this case all the kernels
$g^{[k]}(x;y)$ with $k \geq 2$ are bounded.    \hfill
$\diamondsuit$

\begin{cor}\hspace{-.09in}\textbf{.} If $m$ is large enough, then all
iterated kernels $g^{[m]}$ are bounded in $\Omega \times \Omega $.
They are also continuous in $\Omega \times \Omega$ by
\cite[Th.1.3.1]{Mikh1}. Moreover, $g^{[m]}\in H_0^{1,s}(\Omega)$
with $s<\frac{np}{n-p}$ if $\frac{n}{2}<p\leq n$, and $g^{[m]}\in
C^{(1)}(\Omega)$ whenever $p>n$ (ibid., Theorem 11.5.2.) \hfill
$\diamondsuit$
\end{cor}

\noindent \textbf{Proof of Theorem A.1 (continued).} We reduce the
equation (A.1) to an integral equation, whose resolvent kernel
will be Green's function of  $L_c$. \\

The assertions of Theorem A.1 are valid for the generalized
harmonic Green's function $g$ \cite[Section 5.7.2]{HK}, therefore
in the proof we assume the potential to be strictly positive on a
set of positive measure in $\Omega $.

The equation (A.1) can be written as  $-\Delta u(x)+c(x)u(x)=f(x)$
and further, after multiplying by $g(x;y)$ and integrating over
$\Omega$, as
\begin{equation}
u(x)+\int _{\Omega } g(x;y) c(y) u(y) dy = F(x), \; x \in \Omega ,
\end{equation}
where $F(x) = (gf)(x) = \int _{\Omega } g(x;y) f(y) dy.$

We consider the integral operator of the (symmetrizable) equation
(A.5),
\[(g_c u)(x) = \int _{\Omega } g(x;y) c(y) u(y) dy \]
in a real Hilbert space $L^2 _c (\Omega )$ with the inner product
\[ < f_1 , f_2 >\; = \int _{\Omega } f_1 (x) f_2 (x) c(x) dx, \]
that is, $ \parallel f \parallel ^2 _c = \int _{\Omega } f^2 (x)
c(x) dx.$ The operator $g_c $ is defined on smooth functions with
compact supports in $\Omega $ and is symmetric due to the symmetry
of $g(x;y)$.

Applying the Laplacian $- \Delta $ to the equation $g_c u = 0$, we
arrive at the equation $u(x) c(x) = 0$ for almost all $x \in
\Omega $. The latter means that the condition $g_c u = 0$ implies
an equation $u(x) = \Theta $, where $\Theta $ is the zero element
of the space $L^2 _c (\Omega )$. Thus, the symmetric operator
$g_c$ has a trivial kernel in $L^2 _c (\Omega )$, that is, it is
selfadjoint and, therefore, its spectrum is real. Next we shall
prove that  this spectrum is actually positive.

It is known \cite[p. 364]{Sob}, that the positiveness of the
spectrum is equivalent to the positive definiteness of a quadratic
form
\begin{equation}
<g_c u, u> =  \int _{\Omega } \int _{\Omega } g(x;y) c(x) c(y)
u(x) u(y) dx dy.
\end{equation}
To prove the nonnegativity of (A.6), we use the inequality
\begin{equation}
(- \Delta v, v) = \int _{\Omega } \mid \nabla v \mid ^2 dx > 0,
\end{equation}
which is known to be valid for a non-constant function $v$ in a
domain $\Omega $ with the smooth boundary, provided that this $v$
vanishes at $\partial \Omega $. In (A.7), $(\cdot , \cdot )$ is
the standard inner product in $L^2 (\Omega )$. To prove the
nonnegativity of the form (A.6) in a smooth domain $\Omega $, it
now suffices to apply (A.7) to $v(x)=v_n (x)=\int_{\Omega}g(x;t)
\psi _n (t) dt $, where $\psi _n (t) $ is a finitary function in $
\Omega $, and pass to the limit provided that the sequence $ \{
\psi _n (t) \} ^{ \infty } _{n=1}$ weakly converges to the
function $u(x) c(x)$, where $u$ is an arbitrary function in $L^2
_c (\Omega )$.

Next, to prove the nonnegativity of the form (A.6) in an arbitrary
domain $\Omega $, it suffices to pass to the limit with respect to
Green's functions $g_l (x;y)$ of smooth domains $\Omega _l , \; l
=1, 2,...$, exhausting $\Omega $ from within.

We have proved that the operator $g_c$ has a positive spectrum.
So, $\lambda = -1$ is a regular value of $g_c$, an operator $I +
g_c = I - (-1)g_c $ is invertible and a nonhomogeneous equation
(A.5) is uniquely solvable for every right-hand side $F(x) =
(gf)(x) \in L^2 _c (\Omega )$. However, this latter inclusion
holds for each $f(x) \in L^2 (\Omega )$, which follows immediately
from the definition of the class ${\cal C}(\Omega)$ and from the
Young inequality \cite[p. 271]{Stein}. Therefore, applying the
operator  $ - \Delta $ to the equation (A.5), we straightforwardly
see for ourselves that a solution of (A.5) also satisfies the
equation (A.1). Moreover, this solution of (A.5) can be written as
\begin{equation}
u(x) = (I+g_c )^{-1} F(x) = (I+g_c )^{-1} (gf)(x).
\end{equation}

We will prove next that the operator in (A.8) is an integral one
and will find its kernel - the desired Green function $G(x;y)$.
First, let us note that (if $n=2$, calculations are similar)
\[\int \int_{\Omega \times \Omega}g^2(x;y)c(x)c(y)dxdy \leq
b \int \int _{\Omega \times \Omega} \frac{c(x)c(y)}{|x -
y|^{2n-4}}dxdy \]
\[\leq b \left \{ \int \int _{\Omega \times \Omega} \frac{dxdy}
{|x - y|^{(2n-4)q}} \right \}
^{1/q} \left \{\int \int _{\Omega \times \Omega} [c(x)c(y)]^p dxdy
\right \} ^{1/p}\]
\[\leq b \parallel c \parallel ^2 < \infty ,\]
since $(2n-4)q <2n$, therefore, (A.5) is a Hilbert-Schmidt
equation \cite[p. 1083, Exercise 44]{DS}.

Now, due to Lemma A.1 and this remark, the entire Hilbert-Schmidt
theory of integral equations in $L^2$ extends on integral
equations in $L^2 _c$. In particular, we mean here the results
such as the representation of solutions via the resolvent kernel,
properties of the equations with a weak singularity, the
representation of a symmetric kernel as a bilinear series against
the eigenfunctions of the operator and properties of the latter,
etc., - see, for example, \cite[Lectures 24-25]{Sob}, where the
case of a bounded weight $c(x)$ is treated, and also a remark in
\cite[p. 65-66]{Mikh1}.

It follows that $g_c$ is a compact operator with the positive
spectrum and (A.5) is a Fredholm equation with a symmetric kernel
in $L^2 _c (\Omega )$. Therefore, its solution (A.8) can be
written as $u(x) = F(x) - \int _{\Omega } R(x, y) F(y) c(y) dy, $
where $R(x, y)$ is a resolvent kernel of the operator $g_c $
corresponding to the regular value $\lambda = -1$. Since $ F(x) =
(gf)(x)$, this solution can be rewritten as $u(x) = \int _{\Omega
} G(x;y) f(y) dy, $ where
\begin{equation}
G(x;y) = g(x;y) - \int _{\Omega } R(x, t) c(t) g(t;y) dt
\end{equation}
is Green's function of the operator $L_c$ we sought for. Its
uniqueness follows by the construction, the symmetry is due to the
symmetry of $L_c$ - it is enough to substitute $\int G(x;y)u(y)dy$
for $u$ and $\int G(x;y)v(y)dy$ for $v$ into the equation $(L_c u,
v) = (u, L_c v)$. So, parts $1^{o}$ and $2^{o}$ of Theorem A.1 are proved. \\

We proceed on to the proof of Theorem A.1. Let $\lambda _k, \;
k=1,2,... $, be positive characteristic values of the selfadjoint
operator $g_c$ (where multiple values repeat), thus, the series
$\sum ^{\infty}_{k=1} \lambda _k ^{-2}$ converges, and let $\phi
_k (x)$ be the corresponding eigenfunctions orthonormal in $L^2 _c
(\Omega)$, so the expansion
\[g(x;y) = \sum_{k=1} ^{\infty} \frac{1}{\lambda _k } \phi _k (x) \phi _k
(y)\] converges in the $L^2 _c (\Omega )-$norm. The functions
$\phi _k (x) $ are also eigenfunctions of all iterations of $g_c
$, namely,
\begin{equation}
\phi _k (x) = \lambda _k ^m  \int _{\Omega } g^{[m]} (x;y) c(y)
\phi _k (y) dy, \; m = 1, 2,... ,
\end{equation}
and all series $\sum ^{\infty}_{k=1} \lambda _k ^{-2m}, \; m=1,2,
\ldots ,$ also converge. Moreover, the bilinear series for the
resolvent at a regular value $\lambda =-1$ is
\[R(x,y) = R(x,y,-1) = \sum _{k=1} ^{\infty} \frac{1}{\lambda _k +1}
\phi _k (x) \phi _k (y),\] also converging in $L^2 _c (\Omega)$.
Since $\phi _k (x) \in L^2 _c (\Omega)$ and $c \in L^p$, the
H\"{o}lder inequality implies that $c(x)\phi _k (x) \in L^r
(\Omega)$, where $r= \frac{2p}{p+1} >1$. Therefore, $\phi _k \in
L^{\infty}(\Omega)$ and now (A.10) necessitates that, together
with $g^{[m]}$, all the eigenfunctions $\phi _k \in
H^{1,s}(\Omega), \; s>n$.

Fix an index $m_ 0$ such that for $m \geq m_0$ $g^{[m]}(x;y) \in
H^{1,s} (\Omega), \; s>n$. Assuming $c \in L^p (\Omega)$ and
$\parallel \phi _k \parallel _{L^2 _c (\Omega } = 1$, (A.10) with
$m=m_0$ implies an inequality
\[\mid \phi _k (x) \mid \leq b \lambda  _k ^{m_0} \int _ {\Omega }
\mid \phi _k (y) \mid \sqrt {c(y)} \sqrt {c(y)}dy\]
\[\leq b \lambda  _k ^{m_0} \parallel \phi _k \parallel
\left\{\int _ {\Omega } c(y) dy \right\}^{1/2} = b \lambda  _k
^{m_0} \left\{\int _ {\Omega } c(y) dy \right\}^{1/2}.\] Applying
the H\"{o}lder inequality to the last integral, we get a bound
\[\mid \phi _k(x) \mid \leq b \lambda  _k^{m_0} \sqrt{\parallel c \parallel}
(mes \Omega )^{1/2q}, \; p^{-1} +q^{-1} =1.\]

Taking into account the relationship between the resolvents of the
original kernel and the iterated ones, we arrive at the equation
\[G(x;y)=g(x;y) - g^{[2]}(x;y) +...+(-1)^{m-1} g^{[m]}(x;y)\]
\begin{equation}
+ \sum _{k=1}^{\infty} \frac{1}{\lambda _k ^m (\lambda _k + 1)}
\phi _k (x) \phi _k (y).
\end{equation}
When $m \geq 2$, the latter series converges uniformly, so $G(x;
y)$ is continuous in $\Omega \times \Omega$ if $x \neq y$.
Moreover, it follows from all the above that the boundary behavior
of Green's function $G(x;y)$ and its singularity as $x-y
\rightarrow 0$ are determined by the first term of the right-hand
side of (A.11), that is, by the harmonic Green function $g(x;y)$
of the Laplacian, so all statements of parts $3^o$ and $4^o$ of
Theorem A.1 follow. \\

To prove part $5^o$, we fix a point $y_0 \in \Omega$. Due to part
$3^o$, there exists a ball $B(y_0, \delta)\subset \Omega$, such
that $G(x; y_0)>0$ in this ball. Suppose that $G(x_0; y_0)<0$ at
some point $x_0 \in \Omega $. Since $G$ is a continuous function
with the zero boundary values at $\partial \Omega$, the inequality
$G(x; y_0)<0$ holds in some domain $\Omega ' \subset \Omega
\setminus B(y_0 , \delta )$ and $G(x; y_0 )$ attains its negative
minimum value at a point $x' \in {\Omega '}$. By the definition,
Green's function satisfies the equation $L_c G(x; y_0) = \delta (x
- y_0)$, therefore, $\Delta G(x; y_0) = c(x) G(x; y_0) \leq 0$ in
$\Omega '$, maybe in the sense of distributions. Thus, $G(x; y_0)$
is a superharmonic function in $\Omega '$, which can not attain
its minimum value in a domain. This contradiction proves the
positiveness of $G$. \\

In part $6^o$, the conclusion on the H\"{o}lderian continuity of
$G(x;y)$ follows from the representation (A.11) and the Sobolev
lemma \cite[p. 70-72]{Sob1}: \\
\noindent
\par {\em Let a domain $\Omega \subset B(R), \; f \in L^p(\Omega)$, and
\[U(x,y)=\int _{B(R)}\frac{(r+r_1)^{\lambda -1}}{r^{\lambda }r_1 ^{\lambda }}f(t)dt,\]
where $0<\lambda <n, \; r = |x-t|, \; r_1 = |y-t|$, and $f$ is
extended by zero outside $\Omega$. Then $U(x,y) \leq b
\parallel f\parallel |x-y|^{\beta -1}$, where
$\beta = \min \left\{1, n(1-1/p)- \lambda \right\}$}.\\
\par It suffices to substitute here $\lambda =n-2$ and $f(t)=c(t) |t-y|^{2-n}$. \\

To prove part $7^o$ of Theorem A.1, first we consider the case $n
\geq 3$. Green's function $g(x;y)$ is infinitely differentiable if
$x \neq y$ and iterations can not worsen its smoothness,
therefore, differential properties of $G(x;y)$ depend upon the
second term, $g^{[2]}(x;y)$, of the right-hand side of (A.11).
Moreover, since $g(x;y) =\theta _n r^{2-n} + \; \mbox{a bounded
harmonic function}$, the result is determined by the principal
part of $g^{[2]}(x;y)$, that is, by the integral
\begin{equation}
I(x,y) = \int _{\Omega } \mid x-t \mid ^{2-n} c(t) \mid t-y \mid
^{2-n}dt.
\end{equation}

Again, fix a point $y_0 \in \Omega$ and remark that, due to the
definition of the class ${\cal C}(\Omega)$, a function $\gamma (t)
\equiv c(t) \mid t-y_0 \mid ^{2-n} \in L^{\lambda} (\Omega ) $ for
any $\lambda  \in [1, \frac{np}{n+p(n-2)} )$. Moreover, $\gamma
\in L^p (K)$ on each  compact set $K \subset \Omega \setminus
\{y_0 \}$, where $p$ is from the definition of ${\cal C}(\Omega)$.
Now, all statements of Part $7^o$ ensue from the following known
properties of the potential-type integrals \cite[p. 210-211]{Mikh1}: \\

{\em If $\gamma \in L^{\lambda }(\Omega)$, then}
\[w(x) \equiv \int _{\Omega } \mid x - t \mid ^{2-n} \gamma (t) dt \in
    \left\{\begin{array}{lll}
        H^{1,s} (\mathbf{R}^n ), \forall s \in [1, \frac{n \lambda }{n - \lambda })
            & \mbox{if $\; 1 \leq \lambda \leq n$} \vspace{.1cm} \\
        C(\mathbf{R}^n )  & \mbox{if  $\; \lambda  > n/2$} \vspace{.1cm} \\
        C^1 (\mathbf{R}^n)  & \mbox{if $\; \lambda > n,$}
           \end{array}
    \right. \]
{moreover, there exist generalized second derivatives}
$\frac{\partial ^2 w(x)}{\partial x_i \partial x_j } \in L^{\lambda }(\Omega )$.\\

If $n=2$, then $I(x,y)$ in (A.12) is to be replaced by $\int
_{\Omega} \ln \frac{1}{|x-t|} c(t) \ln \frac{1}{|t-y|} dt$.
Rewriting this integral as
\[\int _{\Omega } \frac{r^{1 + \delta }\ln (1/r)}{r^{\lambda }} \rho (t) dt, \]
where $r=|x-t|, \; 0 < \delta < 1, \; \lambda =1- \delta $ and
$\rho (t)=r^{-2\delta } c(t) \ln \frac{1}{|t-y|}$, and applying to
the function $\rho $ the triple H\"{o}lder inequality with the
exponents
\[\frac{2p\delta  +1 -\delta }{p\delta }, \; \frac{2p\delta  +1 -\delta }{1 - \delta }, \;
\frac{2p\delta  +1 -\delta }{p\delta },\] we see that $\rho  \in
L^s(\Omega)$ with $s=\frac{p(1- \delta)}{2p\delta +1 -\delta}$ for
each $\delta , \; 0<\delta < \frac{p-1}{3p-1}$. Now we can
immediately apply the following statement \cite[Theorem 2.8.1]{Mikh1}: \\

\emph{ Let $u(x)= \int _{\Omega } \frac{A(x,t)}{r^{\lambda }} \rho
(t) dt $ where $\Omega \subset \mathbf{R}^n$ is a bounded domain,
$\lambda < n-1, \; \rho \in L^s (\Omega ), \; 1 \leq s \leq
\infty$, and functions $A(x,t)$ and $\partial A(x,t) /
\partial x_j $ are continuous in $\overline {\Omega } \times
\overline {\Omega }$. Then there exist generalized partial
derivatives
\[\frac{\partial u(x)}{\partial x_j}=\int _{\Omega } \frac{\partial }{\partial x_j}
\left[\frac{A(x,t)}{r^{\lambda }} \right] \rho (t) dt \in L^t
(\Omega ), \; 1 \leq j \leq n, \] for $t < \frac{ns}{n-(n-\lambda
-1)p}$ if $\frac{(\lambda +1)s}{s-1} \geq n$.}

We substitute here $n=2$, $\lambda =1- \delta $, and $s$ as above. \\

Finally, we proceed on to part $8^0$. Nonnegativity of the inner
normal derivative at points, where it exists, follows directly
from the positiveness of $G$. If $p>n$, both the existence and the
continuity of $\frac{\partial G(x;y)}{\partial n}$ follow from the
known facts about the differentiability of the integrals with weak
singularity - see the proof of Part $7^0$ above. Therefore, in the
following we assume $n/2 < p \leq n$. Again, we represent $G$ by
the sum (A.11). It suffices to study its principal term - the
integral (A.12). Let us fix two points, $x_0 \in E_1, \; y_0 \in
\Omega $, and a domain $\Omega _0 \subset \Omega $, such that $y_0
\not\in \overline{\Omega _0}$ and $E_1 \subset \partial \Omega _0
$. Obviously, it suffices to study an integral
\[I_{00}(x) = \int _{\Omega _0 } \mid x-t \mid ^{2-n} c(t) \mid t-y_0 \mid ^{2-n} dt , \]
since a volume integral $\int _{\Omega \setminus \Omega _0 } \mid
x-t \mid ^{2-n} c(t) \mid t-y_0 \mid ^{2-n}dt $ is a harmonic
function in $x$ at a vicinity of the point $x_0 $. We extend
$c(t)$ with zero values into a domain $\Omega _1 $ lying in the
complement to $\Omega $ and such that $E_1 \subset \partial \Omega
_1$. Then $ c(t) \mid t-y_0 \mid ^{2-n} \in L^p (\overline{ \Omega
_0 \cup \Omega _1})$, and as we showed in Part $7^0$, there exist
the distributional partial derivatives
\[\frac{\partial I_{0}(x)}{\partial x_j}=
\int _{\Omega _0}b_j(x,t)|x-t|^{1-n}c(t)|t-y_0|^{2-n}dt ,\] where
$b_j(x,t),\; 1 \leq j \leq n$, are bounded functions. Now the
statement of Part $8^0$ follows from the following theorem \cite[p. 24]{Mikh1}:\\

{\em Let a function $\rho (t) \in L^p(\Omega)$, and there exist an
integer $l$ and a real number $\lambda$ such that
$n-(n-\lambda)p<l \leq n$ and $0<n-\frac{n}{p} \leq \lambda<n$;
let $\omega _l$ be a cross-section of the set $\Omega$ with an
$l-$dimensional $C^1$-smooth manifold. Then a function $w(t) =
\int _{\Omega _0 } \mid x-t \mid ^{- \lambda } \rho (t) dt $ is
defined on every such cross-section almost everywhere in the sense
of $l-$dimensional Lebesgue measure and $u(t) \in L^s (\omega _l)$
for each $s \in [1, \; \frac{lp}{n - (n - \lambda )/p})$.} \\

We substitute here $l=\lambda =n-1$; the condition $\lambda \geq n
- n/p$, which is equivalent to $\frac{p(n-1)}{p-1} \geq n$,
obviously holds true for now $p \leq n$. \\

The proof of Theorem A.1 is complete. \hfill \textbf{$\diamondsuit$}\\

\samepage{ \noindent\textbf{Remark A.1.} Expanding the resolvent
kernel into the Neumann series, we get from (A.9) a series
\begin{equation}
G(x;y) = \sum ^{\infty}_{m=1} (-1)^{m-1} g^{[m]} (x;y),
\end{equation}
converging in $L^2 _c(\Omega )$. Due to Lemma A.1, if $\|c\|$ is
sufficiently small, then the series (A.13) converges absolutely
and uniformly. In particular, this is the case if $G(x;y)$ is
Green's function of a small ball compactly imbedded in the domain
$\Omega $, where $c \in {\cal C}(\Omega)$. Moreover, if $\|c\|$ is
sufficiently small and the boundary $\partial \Omega$ is smooth,
then the series (A.13) can be differentiated termwise, therefore
\[\frac{\partial G(x;y) }{\partial n(y)} = \sum ^{\infty}_{m=1}
(-1)^{m-1} \frac{\partial g^{[m]} (x;y)}{\partial n(y)}.\] \hfill
\textbf{$\diamondsuit$}}

For completeness we include here a short proof of the following
known result; see, for example, a concise survey in \cite{Anc} and
references therein.

\begin{theor}\hspace{-.09in}\textbf{.} Let a domain $\Omega$ and functions
$c_{k}(x), \; k = 1,2$, satisfy the conditions of Theorem A.1, and
$G_{k}$ be Green's functions of the operators $L_{c_k}$. If $c_1
\leq c_2$ almost everywhere in $\Omega $, then $G_2 (x;y) \leq G_1
(x;y)$ everywhere in $\Omega \times \Omega $. Moreover, if $c_1
(x) < c_2 (x)$ on a set of positive measure, then $G_2 (x;y) < G_1
(x;y)$ everywhere in $\Omega \times \Omega , \; x \neq y$.

Vice versa, if there exists at least one point $(x_0 ,y_0 ) \in
\Omega \times \Omega, \; x_0 \neq y_0 $, where $G_1 (x_0; y_0 ) =
G_2 (x_0; y_0)$, then $c_1 = c_2$ almost everywhere in $\Omega $.

In particular, for each potential $c(x) \in {\cal C} (\Omega )$
\begin{equation}
G(x;y) \leq g(x;y)
\end{equation}
where, again, $g$ is Green's function of the Laplacian. \hfill
\textbf{$\diamondsuit$}
\end{theor}

\noindent \textbf{Remark A.2.} Opposite inequalities are discussed
in
\cite{Anc} and \cite{CrFaZh}. \hfill \textbf{$\diamondsuit$}\\

\noindent \textbf{Proof of Theorem A.2.} Applying the inverse
operator $\int _{\Omega } G_1 (x;y)f(y)dy$ of the operator
$L_{c_1}$ to the equation
\begin{equation}
L_{c_{2}} u = f(x),
\end{equation}
we deduce an equation
\begin{equation}
u(x)= \int _{\Omega } G_2 (x;y)\left(c_2 (y) - c_1 (y)\right)
u(y)dy + \int _{\Omega } G_1 (x;y)f(y) dy.
\end{equation}
On the other hand, from (A.15)
\begin{equation}
u(x) = \int _{\Omega } G_2 (x;y)f(y) dy.
\end{equation}
Fix a point $y_0 \in \Omega $. Choosing for $f$ in (A.16)-(A.17) a
$\delta $-sequence $f_n (x) \geq 0$ such that
\[\int _{\Omega}f_n (x) \chi (x)dx \rightarrow \chi (y_0),\; n\rightarrow \infty ,\]
we obtain from (A.17) that $u(x) \rightarrow G_2 (x; y_0) \geq 0$
and then from (A.16) the relation between the Green functions
\begin{equation}
G_2(x; y_0) = \int _{\Omega } G_1 (x;y) \left(c_1 (y) - c_2 (y)
\right) G_2 (y; y_0)dy + G_1 (x; y_0).
\end{equation}
The theorem follows. \hfill  \textbf{$\diamondsuit$}\\

\noindent \textbf{Remark A.3.} The equation (A.18) has appeared in
the literature in various instances. See, for example, the
equation (4.7.3) in \cite{CrFaZh}, where $c_1 =0$, $c_2$ is a Kato
class potential, and the authors consider a general divergence
form elliptic operator of second order in a Lipschitz domain.
\hfill \textbf{$\diamondsuit$}

\begin{cor}\hspace{-.09in}\textbf{.} Under the conditions of Theorem A.2,
\[\partial G_{2}(x;y) / \partial n\leq \partial G_{1} (x;y) / \partial n \]
everywhere where these derivatives exist. \hfill
\textbf{$\diamondsuit$}
\end{cor}

The next result is needed in the proof of Theorem 2.1.

\begin{cor}\hspace{-.09in}\textbf{.} If $\| c_1-c_2\|_{L^p (\Omega)}<
\varepsilon, \; \varepsilon >0, \; p>n/2$, then
\[\| G_2(\cdot ;\; y_0)-G_1(\cdot ;\; y_0\|_{L^r (\Omega)}<b\; \varepsilon \]
for any $r, \; 1\leq r < \frac{n}{n-2}$, uniformly in $y_0$.
\hfill \textbf{$\diamondsuit$}\\
\end{cor}

\noindent\textbf{Proof.} It suffices to estimate the integral in
(A.18) by the Minkowski inequality with any $r, \; 1\leq
r<n/(n-2)$, observe that $\| G_k (\cdot ; y)\|_{L^r _{(\Omega )}}
< \infty$ uniformly in $y\in \Omega$, and then apply the
H\"{o}lder inequality with the indices $p$ and $q,\;
p^{-1}+q^{-1}=1,\; p>n/2$. \hfill \textbf{$\diamondsuit$}\\

\noindent \textbf{Remark A.4.} Theorem A.2 has a clear physical
meaning. Let us consider a bounded domain $\Omega $ in
$\mathbf{R}^2$ (a homogeneous isotropic membrane) imbedded in
$\mathbf{R}^3$. Suppose that a stationary heat distribution is
given in $\Omega $ and the edge $\partial \Omega $ is maintained
at zero temperature. Then the temperature $u(x)$ at $x \in \Omega
$ satisfies the equation $\Delta u -c(x)u(x)=0$, where
$c(x)=\lambda ^{-1} c_0 (x)$, $\lambda $ being the heat
conductivity coefficient and $c_0 (x)$ the coefficient of heat
emission from the surface of the membrane. If there is a heat
source of unit capacity at a point $y \in \Omega $, then the
corresponding temperature distribution $u(x), \; x \in \Omega $,
on the membrane is nothing but Green's function $G(x;y)$ of the
operator $L_c$ in $\Omega $. It is now obvious that, if we
increase $c_0 (x)$, that is, the output of heat into the
surrounding space, while keeping all other parameters unchanged,
then the temperature on the membrane must decay, so $u(x) = G(x;y)$ decreases. \\

The same argument holds true for nonhomogeneous and anisotropic
membrane. In this case, we obtain an inequality between Green's
functions of a divergence form operators
\[L_k u(x)=- \sum _{j=1}^n \frac{\partial}{\partial x_j}\left\{\lambda _j (x)
\frac{\partial u}{\partial x_j} \right\} + c_k (x)u(x)=0, \;
k=1,2,\]
with $c_1 (x) \leq c_2 (x)$. \hfill \textbf{$\diamondsuit$}\\

\noindent \textbf{Remark A.5.} It is worth noting that the known
equation
\[\int _{\partial \Omega} \frac{\partial G(x;y)}{\partial n(y)}
d\sigma(y) + \int _{\Omega} c(y) G(x;y) dy =1, \; \forall x\in
\Omega ,\] has a physical meaning as well. The first
integral here is equal to the amount of heat outgoing through the
boundary of the membrane in a unit of time, given a unit heat
source at a point $x$, while the second integral represents the
amount of heat outgoing into the surrounding space from the
surface of membrane under the same conditions. \hfill
\textbf{$\diamondsuit$}\\

\newpage

\section{Eigenfunctions of the Laplace-Beltrami \\
operator on the unit sphere in $\mathbf{R}^n$}

Laplace-Beltrami operator $\Delta ^*$ is the spherical part of the
Laplacian $\Delta $, that is, its restriction onto functions
defined on the unit sphere $S\equiv S^{n-1}\subset \mathbf{R}^n$.
In this section we study eigenfunctions of the eigenvalue problem

\begin{equation} \left\{
\begin{array}{ll}
\Delta ^* \varphi (\theta )+\lambda \varphi (\theta )=0
        & \mbox{if  $\; \theta \in D$} \vspace{.1cm} \\
\varphi (\theta ) =0  & \mbox{if
        $\; \theta \in \partial D \setminus E$} \vspace{.1cm} \\
\varphi \in L^2(D)  & \mbox{}
\end{array}
\right.
\end{equation}
in a domain $D \subset S$, where $E$ is an exceptional set (maybe
empty) of irregular points on $\partial D$. \\

The coefficients of $\Delta ^*$ in standard spherical coordinates
have singularities of kind $\frac{1}{\sin ^2\theta}$ \cite[p. 393]
{Sob}, therefore, we cannot straightforwardly apply the known
theorems on the existence of eigenfunctions of the Dirichlet
problem  for elliptic operators, even if the boundary $\partial D$
is smooth. However, a representation of the Laplacian in a
coordinate-free form implies that those singularities are not
essential, they disappear under a rotation of coordinate axes.

Indeed (see, for example, \cite[p. 409]{Sob}, let $f(\theta)$ be a
function defined in a domain $D \subset S$. We extend it as $f(r,
\theta ) = f(\theta )$ for $0<r<2$. Now, $\Delta f(r, \theta ) =
r^{-2} \Delta ^* f(\theta )$ independently upon the spherical
coordinates chosen, and the left-hand side of the latter equation
does not depend on a choice of spherical coordinates, therefore,
the right-hand side does not depend either. So, the expression
$\Delta ^* f(\theta )$ always returns the same function
independently on a system of spherical coordinates. Thus, we can
locally select an appropriate coordinate system and the operator
$\Delta ^*$ has no singularities. Therefore, it has Green's
function and the eigenvalue problem (B.1) has eigenfunctions in
each domain with a smooth, for example, $C^2$-boundary (Oleinik
\cite{Ole}). The main result of this section states
that even this smoothness can be essentially relaxed.\\

We need two following remarks in Section 3. First, we show that
the gradient vector $\nabla f$ can also be defined independently
upon a coordinate system. Indeed, let $\gamma $ be a smooth curve
on $S$, $\theta \in \gamma $, and $ds$ be the arc differential on
$\gamma $. Then $\frac{df}{ds}$ does not depend on a coordinate
system by definition. Now, let us consider on $S$ two orthogonal
coordinate systems, $(t_1 ,..., t_{n-1})$ and $(\tau _1 ,..., \tau
_{n-1})$. We have
\[\frac{df}{ds}=\sum^{n-1} _{j=1} \frac{\partial f}{\partial t_j}
\frac{dt_j}{ds}=\sum ^{n-1} _{j=1} \frac{\partial f}{\partial t_j}
\cos \alpha _j = \left(\nabla _t f, \vec {e} \right),\] where
$\vec {e}$ is a unit tangent vector to $\gamma $, and similarly,
$\frac{df}{ds} = \left(\nabla _{\tau }f, \vec {e} \right)$.
Therefore, $\left(\nabla _t f - \nabla _{\tau }f, \vec {e} \right)
= 0$ for each tangent vector $\vec {e}$, and since $\nabla f$
itself lies in the tangent plane, $\nabla _t f = \nabla_{\tau
}f$.\\

Second, while calculating an integral $\int _D \left((\nabla u)^2
+ c(\theta )u^2 \right) d\sigma (\theta ), \; D \subset S$, we can
divide $D$ into small pieces and select an appropriate system of
spherical coordinates locally in each piece of $D$.\\

The following theorem ascertains the existence of the Green
function and the eigenfunctions of the Laplace-Beltrami operator
$- \Delta ^* $ in each domain on the sphere such that its boundary
is not a polar set. In the proof we exhaust $D$ from within with a
sequence of expanding domains $D_j \subset \subset D_{j+1}
\nearrow D, \; j=1,2,...,$ having $C^2$-smooth (and therefore,
regular) boundaries $\partial D_j$. Denote by $B^{(j)}(\theta ;
\psi )$ Green's function of $- \Delta ^* $ in $D_j$ with zero
boundary values on $\partial D_j$. The eigenvalues of the problem
(B.1) in the domain $D_j$, repeating accordingly to their
multiplicities, and the corresponding eigenfunctions are denoted,
respectively, by $\lambda ^{(j)} _{\nu}$ and $\varphi ^{(j)} _{\nu
}(\theta ), \; \nu =0,1, \ldots \; $. It is known that $\lambda
_0$ is simple and $\varphi ^{(j)} _0$ does not change its sign, so
we assume $\varphi ^{(j)} _0 >0$. The Green functions $B^{(j)}
(\theta ; \psi)$ and the eigenfunctions $\varphi ^{(j)}
_{\nu}(\theta )$ are extended by zero for $\theta \in S\setminus
D_j$. We normalize the eigenfunctions in $L^2(D)$, $\parallel
\varphi ^{(j)} _{\nu } \parallel _{L^2 (D)} = 1$.

\begin{theor}\hspace{-.09in}\textbf{.} Let $D \subset S$ be a domain whose
boundary $\partial D$ with respect to $S$ is not an
$(n-1)-$dimensional polar set. Then

$1^0$ There exists Green's function $B(\theta ; \psi )$ of the
operator $- \Delta ^*$ in $D$ with zero boundary values outside an
exceptional (maybe empty) polar set $E \subset \partial D$ of
irregular boundary points.

$2^0$ The eigenvalue problem (B.1) in $D$ has the eigenvalues
$\lambda _{\nu}, \; \nu =0,1,...$, with the corresponding
eigenfunctions $\varphi _{\nu }$; the $\lambda _0$ is simple and
$\varphi _0 >0$.

$3^0$ The eigenfunctions $\varphi ^{(j)} _0(\theta ) $ converge to
$\varphi _0 (\theta ) $ in $L^2-$norm,
\[\|\varphi ^{(j)}_0(\theta)-\varphi_0(\theta)\|_{L^2(D)}\rightarrow
0 \mbox{ as } j \rightarrow \infty .\] $4^0$ Moreover, $\varphi
^{(j)} _0 (\theta ) \rightarrow \varphi _0 (\theta )$ as $j
\rightarrow \infty $, uniformly on each compact set $Q \subset D$.
\end{theor}

\noindent\textbf{Remark B.1.} The same statement can be proved for
all other eigenfunctions $\varphi _{\nu}, \; \nu = 1,2,...\; $.
\hfill \textbf{$\diamondsuit$}\\

\noindent \textbf{Proof of Theorem B.1.} Green's function
$B^{(j)}(\theta ; \psi )$ of the operator $-\Delta ^*$ in
$D^{(j)}$ is strictly positive inside the smooth domain $D_j$ and
continuously takes on zero boundary values on $\partial D_j, \;
j=1,2,...\; $. Due to the monotonicity of the sequence $D_j, \;
j=1,2,...$, the sequence $B^{(j)}(\theta ; \psi )$ also steadily
increases at each point $(\theta , \psi ) \in D \times D, \;
\theta \neq \psi $, and is bounded from above by Green's function
of the same operator $- \Delta ^* $ in any regular domain
$\widehat{D} \supset D$. Therefore, there exists a function
$B(\theta ; \psi ) =  \lim _{j \rightarrow \infty } B^{(j)}
(\theta ; \psi )$. Clearly, $B(\theta ; \psi ) > 0$ everywhere in
$D$ and since the set of irregular points on $\partial D$ is
polar, $B(\theta ; \psi )$ takes on zero boundary values on
$\partial D$ outside an exceptional polar set $E$.

Now, by the definition of Green's function $B^{(j)}(\theta ; \psi
)$, for each finitely supported in $D_j$ smooth function $\chi
(\theta )$ there holds an equation
\[\chi (\theta )=\int _{D_j} B^{(j)}(\theta ; \psi) \left(-\Delta ^*
\chi (\psi )\right) d\sigma (\psi)= \int _D B^{(j)} (\theta ; \psi
)\left(-\Delta ^* \chi (\psi )\right) d\sigma (\psi ),\] for we
extend $B^{(j)} = 0$ outside $D_j$. On the right we can let $j
\rightarrow \infty $, since $B^{(j)}$ is dominated by the summable
function $B$. Thus,
\[ \chi (\theta ) = \int _D B (\theta ; \psi )\left(-\Delta ^* \chi (\psi )\right)
d\sigma (\psi ) \] for all finitary functions $\chi$ in $D$ and so
that over the entire domain of the operator $- \Delta ^* $. We
have shown that $B(\theta ; \psi)$ is the (generalized) Green's
function of $- \Delta ^* $ in the domain $D$ with zero boundary
values on $\partial D \setminus E$. \\

To construct the eigenfunctions of $-\Delta ^*$ in the domain $D$,
we consider integral operators
\[B^{(j)}f(\theta) \stackrel{\rm def}{=}
\int_{D_j}B^{(j)}(\theta ; \psi)f(\psi)d\sigma(\psi)\] and
\begin{equation}
Bf(\theta ) \stackrel{\rm def}{=} \int _D B(\theta ; \psi ) f(\psi
) d\sigma (\psi ).
\end{equation}
The former operator is defined on $L^2 (D_j)$ and the latter on
$L^2 (D)$.

Since Green's functions of the domains $D_j$ and $\widehat{D}$,
extended by zero, belong to  $L^2 (\widehat{D} \times
\widehat{D})$, the kernels $B^{(j)} (\theta ; \psi ) \rightarrow B
(\theta ; \psi )$ in the $L^2 (D \times D)$ norm as $j\rightarrow
\infty$. Therefore, as $j \rightarrow \infty$, $B^{(j)}
\rightarrow B$ in the operator norm. Due to the minimax principle,
the sequence $\{\lambda ^{(j)}_0 \}_{j \geq 1}$ monotonically
decreases, and since $\lambda ^{(j)}_0 \geq \widehat{ \lambda
_0}$, where $\widehat{ \lambda _0}> 0$ is the principal eigenvalue
of a smooth domain $\widehat{D} \supset D_j, \; \widehat{D} \neq
S$, there exists a limit $\lambda _0 = \lim _{j \rightarrow \infty
} \lambda ^{(j)}_0 > 0$. The same variational principle implies
that the limit does not depend on an approximating sequence of
domains $D_j$. We have to show yet that $\lambda _0$ is the first
eigenvalue of the operator $- \Delta ^* $ in $D$ and to construct
its corresponding eigenfunction $\varphi _0 (\theta )$.

Since $\varphi ^{(j)} _0 = \lambda ^{(j)}_0 B^{(j)} \varphi
^{(j)}_0$, we have
\begin{equation}
\varphi ^{(j)}_0 -\lambda _0 B \varphi ^{(j)}_0 = \lambda ^{(j)}_0
\left(B^{(j)} - B\right) \varphi ^{(j)}_0 + \left(\lambda ^{(j)}_0
- \lambda _0 \right) B \varphi ^{(j)}_0.
\end{equation}
The right-hand side of (B.3) tends to $0$ in $L^2 (D)-$norm as $j
\rightarrow \infty $, because $\parallel \varphi ^{(j)}_0
\parallel = 1$, and so that
\begin{equation}
\parallel \varphi ^{(j)}_0 - \lambda _0 B \varphi ^{(j)}_0 \parallel _{L^2 (D)}\;
\rightarrow 0 \mbox{ as } j \rightarrow \infty .
\end{equation}

Next, due to the compactness of the operator $B$ defined in (B.2),
the function family $\{B \varphi ^{(j)}_0 \}_{j \geq 1}$ is
compact and contains a fundamental subsequence $\{B \varphi
^{(j_s)}_0 \}_{s \geq 1 }$. So, due to (B.4), the sequence
$\{\varphi ^{(j_s)}_0 \}_{s \geq 1 }$ is itself fundamental and
converges to a function $\varphi _0 \in L^2 (D)$. Now the equation
(B.3) implies
\begin{equation}
\varphi _0 (\theta ) = \lambda _0 B \varphi _0(\theta ),
\end{equation}
therefore, we have constructed the first eigenfunction $\varphi _0
(\theta ) \in L^2 (D)$ of the operator $B$ in the domain $D$ and
proved that $\parallel \varphi ^{(j)}_0 - \varphi _0 \parallel
_{L^2 (D)} \rightarrow 0$ as $j\rightarrow \infty $. Moreover, it
is clear that $\lambda _0$ is simple and we can choose $\varphi _0
(\theta) \geq 0$.\\

We have proved the existence of $\varphi _0(\theta) \geq 0$. To
prove that it is strictly positive, we consider a function $u_0
(r,\theta)=r^{\mu ^+}\varphi_0(\theta )$ in the cone $K^D=D\times
(0, \infty)$, where $\mu ^+$ is the characteristic constant of the
domain $D$, that is, the positive (since $\lambda _0 >0$) root of
the quadratic equation $\mu (\mu + n - 2) = \lambda _0$. By the
construction, $u_0$ is harmonic in $K^D$ and vanishes at its
boundary. If there exists a point $\theta _0 \in D$ such that
$\varphi _0 (\theta _0) = 0$, then the nonnegative harmonic
function $u_0$ vanishes along the ray $(r, \theta _0), \; 0<r<
\infty , \; u_0 (r, \theta _0) = 0$. Therefore, $u_0 \equiv 0$ in
$K^D$. Thus, it would be $\varphi _0 (\theta ) \equiv 0$, despite
$\parallel \varphi ^{(j)} _0 \parallel _{L^2 (D)} = 1$. \\

To prove the last part of the theorem, we fix a compact set $Q
\subset D$ and consider $m-$fold iterated kernels
$\left(B^{(j)}\right)^{[m]}(\theta ; \psi )$ and $B^{[m]} (\theta
; \psi ), \; m =1,2,...$ (Cf. Lemma A.1). If $m$ is large enough,
these kernels are continuous in $D \times D$ by Corollary A.1;
moreover, if the domain $D$ is smooth, then they are even
continuous in $\overline{D} \times \overline{D} $. Therefore, the
functions $B^{[m]} \varphi ^{(j)}_0 (\theta ), \; j=1,2,...$, are
equicontinuous in $Q$. Iterating (B.5) and a similar equation for
$\varphi ^{(j)}_0$, we arrive at the equations
\[\varphi ^{(j)}_0 (\theta ) = (\lambda ^{(j)}_0)^m
\int _D \left(B^{(j)}\right)^{[m]} (\theta ; \psi ) \varphi
^{(j)}_0 (\psi ) d \sigma (\psi )\]
and
\[\varphi _0(\theta ) = \lambda _0 ^m \int _D B^{[m]}
(\theta ; \psi ) \varphi _0 (\psi ) d \sigma (\psi ).\] So, if $j$
is fixed, then all the iterated kernels
$\left(B^{(j)}\right)^{[m]},\; m =1,2,...$, have the same
principal eigenfunction $\varphi ^{(j)}_0$ corresponding to the
$m-$th powers of the same eigenvalue $\lambda ^{(j)}_0$, and all
$B^{[m]}$ have the same eigenfunction $\varphi _0$ corresponding
to the $m$-th powers of the same eigenvalue $\lambda _0$. We can
now rewrite (B.3)-(B.4) as
\begin{equation}
\varphi ^{(j)}_0 - \lambda ^m _0 B^{[m]} \varphi ^{(j)}_0
=(\lambda ^{(j)}_0 )^m \left(\left(B^{(j)}\right)^{[m]} - B^{[m]}
\right) \varphi ^{(j)}_0 + \left((\lambda ^{(j)}_0 )^m - \lambda
^m _0 \right) B^{[m]} \varphi ^{(j)}_0
\end{equation}
and
\[\parallel \varphi ^{(j)}_0 - \lambda ^m _0 B^{[m]}\varphi^{(j)}_0
\parallel _{L^2(D)} \; \rightarrow 0, \; j \rightarrow \infty . \]

By means of the compactness criterion of a sequence of continuous
functions, we can select a uniformly convergent in $Q$ subsequence
of functions $\{ \lambda _0 ^m B^{[m]} \varphi ^{(j_s)}_0 (\theta
) \}_{s \geq 1}$. In addition, the right-hand side of (B.6),
obviously, tends to zero uniformly in $\theta \in Q \subset D$,
and so that $\varphi ^{(j)}_0 (\theta ) \rightarrow \varphi _0
(\theta )$ uniformly in $Q$. The proof is complete. \hfill
\textbf{$\diamondsuit$}

\newpage

\section{Special solutions of the equation \\
$y''(r)+(n-1)r^{-1}y'(r)-\left(\lambda r^{-2}+q(r)\right)y(r)=0$}

The equation in title is precisely the equation (3.4). In this
section we establish necessary properties of two special linearly
independent positive solutions $V(r)$ and $W(r)$ of this equation
called, respectively, {\em upper} (increasing as $r\rightarrow
\infty$) solution and {\em lower} (decreasing as $r\rightarrow
\infty$) solution. If $q$ is continuous, these properties can be
found, for example, in \cite[Chap. XI, Sect. 6, in particular,
Corollary 6.5]{Hart}, however, they hold in the case of locally
integrable potentials as well. In this case we assume a standard
stipulation that solutions considered are absolutely continuous
along with their first derivatives and satisfy (3.4) almost
everywhere. In what follows, $\mu^{\pm}=(1/2)\left(2-n \pm
\sqrt{(n-2)^2+4 \lambda} \right)$ are the roots of the quadratic
equation $ \mu ( \mu + n - 2 ) = \lambda $. Denote also $\chi =
\mu^+ -\mu^- =\sqrt{(n-2)^2 +4\lambda}$. If it is necessary to
indicate dependence of the solutions on the parameter $\lambda $,
we write $V(r, \lambda)$ and $W(r, \lambda)$.

\begin{lem}\hspace{-.09in}\textbf{.} Let $\lambda \geq 0$ be a nonnegative
constant and $0 \leq q(r) \in L_{loc}(0, \; \infty )$, that is,
$q$ is integrable over any segment $ [a,b] \subset (0, \; \infty ) $. Then \\

$1^{o}$ The equation (3.4) has a fundamental set of positive
solutions $ \{V, \; W \} $ on $ (0, \; \infty ) $ such that
\[V(0^+)\geq 0\; \mbox{ and }\; \frac{d}{dr}\left(r^{-\mu ^{+}}V(r)\right)\geq 0
\mbox{ for } r >0,\]
\[ W(+\infty)=0\; \mbox{ and } \;\frac{d}{dr}\left(r^{-\mu ^{-}}W(r)\right)\leq 0
\mbox{ for } r>0. \] In particular, $V(r) \geq b \; r^{\mu ^{+}}$
for $1<r< \infty$ and does not decrease on $(0, \; \infty)$, and
$W(r) \leq b \; r^{ \mu ^{-}}$ for $0<r<1$ and steadily decreases
on $(0, \; \infty)$. If $\lambda > 0$, then $V(0^+)=0$ and
$V(r)$ steadily increases on $(0, \infty)$.\\

$2^{o}$
\begin{equation}
V(r)W(r) = \underline{\underline{O}}(r^{2-n}) \; \; as \; \; r
\rightarrow + \infty \; \; and \; \; as \; \; r \rightarrow 0^+.
\end{equation}\\

$3^{o}$ If $ \lambda \searrow \lambda _{0} \geq 0 $, then
\[V(r) \equiv V(r,\lambda)\rightarrow V(r,\lambda _{0})\]
and
\[W(r) \equiv W(r,\lambda)\rightarrow W(r,\lambda _{0})\]
uniformly with respect to $r$ on every compact set $[a,b]\subset
(0,\; \infty)$. \hfill \textbf{$\diamondsuit $}
\end{lem}

\noindent\textbf{Sketch of the proof.} If we assume, in addition,
that
\begin{equation}
\int _0 tq(t) dt < \infty ,
\end{equation}
then (3.4) is equivalent to a Volterra equation
\begin{equation}
y(r) = r^{\mu ^+ } + \int ^r _0 C(r,s) y(s) q(s) ds ,
\end{equation}
where $C(r,s) = (s / \chi ) \left( \left(\frac{r}{s}\right)^{\mu
^+ } - \left(\frac{r}{s}\right)^{\mu ^- } \right)$ is the Cauchy
function of (3.4). Iterating (C.3) as in Lemma A.1, we construct a
solution of (3.4)
\[ V(r) = r^{\mu ^+ } + \int ^r _0 K(r,s) s^{\mu ^+ } q(s) ds,\]
where $K(r,s) = \sum _{j=1} ^{\infty } C^{[j]}
(r,s)$ is the resolvent of (C.3) and $C^{[j]}$ are iterations of
the kernel $C(r,s)$, as in (A.2). Under the condition (C.2) this
series converges uniformly on $[0,b]$ for any $b>0$. A linearly
independent solution is given by the standard formula
\begin{equation}
W(r) = V(r) \int ^{\infty } _r t^{1-n} V^{-2} (t) dt.
\end{equation}

To remove the restriction (C.2), we consider a cut-off function

\[q_N (r) = \left\{\begin{array}{ll}
            q(r) \mbox{ for $\; 0<r \leq N$} \vspace{.3cm} \\
            0    \mbox{  for $\; \; N < r$}
                   \end{array}
            \right. \]
and its Kelvin transform $q_{1,N}(r) = r^{-4} q_N (1/r)$, which
obviously satisfies (C.2). Therefore, the equation
\begin{equation}
y^{''}(r)+(n-1)r^{-1}y^{'}(r)-\left(\lambda
r^{-2}+q_{1,N}(r)\right)y(r) = 0
\end{equation}
possesses a solution $V_{1,N}(r)$ with the properties needed.
Using (C.4) with $V_{1,N}$, we construct a linearly independent
solution $W_{1,N}(r)$ of (C.5). But now the function $V_N (r) =
r^{2-n} W_{1,N} (1/r)$ satisfies the equation
\begin{equation}
y^{''}(r)+(n-1)r^{-1}y^{'}(r)-\left(\lambda
r^{-2}+q_N(r)\right)y(r)=0,
\end{equation}
that is, $V_N (r) = V(r)$ for $0<r<N$. A linearly independent
solution of (C.6) is given by $W(r) = r^{2-n} V_1 (1/r)$, where
$V_1$ is the solution, just constructed, of the equation (3.4)
with $q_1 (r) = r^{-4} q(1/r)$. All the conclusions of parts
$1^o$ and $2^o$ follow straightforwardly from this construction.\\

To prove part $3^o$ of Lemma C.1, we use the same iteration method
and notations as before. The uniform convergence $V_{1,N}(r,
\lambda) \rightarrow V_{1,N} (r,\lambda _0), \; \lambda
\rightarrow \lambda _0$, on compact sets $[a,b]\subset (0,\infty
)$ follows by a direct estimation. Applying the inequality
$V_{1,N} (r,\lambda) \geq b r^{\mu ^+ }$ to estimate the tail of
the integral in the following representation of a linearly
independent solution,
\[W_{1,N} (r,\lambda) = V_{1,N} (r,\lambda) \int _r ^{\infty }
t^{1-n} V^{-2} _{1,N} (t,\lambda)dt ,\] we immediately see for
ourselves that $W_{1,N} (r,\lambda) \rightarrow W_{1,N} (r,\lambda
_0), \; \lambda \rightarrow \lambda _0$, uniformly on compact sets
$[a,b] \subset (0, \infty )$. Introducing $V(r,\lambda)= r^{2-n}
W_1 (1/r,\lambda)$, we prove that $V (r, \lambda) \rightarrow
V(r,\lambda _0)$ and $W (r,\lambda) \rightarrow W(r,\lambda _0)$
as $\lambda \rightarrow \lambda _0 $. \hfill
\textbf{$\diamondsuit $}\\

We need asymptotic formulas for these solutions. Denote $s(r)=r^2
q(r)$ and $\chi ^2 _k = (n-2)^2 + 4 ( \lambda + k)$. If there
exists the limit
\begin{equation}
\lim _{r \rightarrow \infty } s(r) = k < \infty ,
\end{equation}
then asymptotically, for any $\epsilon > 0$
\begin{equation}
b \; r^{ - \epsilon } < V(r) r ^{ - (2 - n + \chi _{k} )/2 } < b
\; r^{ \epsilon } .
\end{equation}
In general, the estimate (C.8) can not be improved. For example,
if
\[q(r) = \frac{a}{r^{2}(1+\ln r)}, \; a= const, \;\; r>1, \]
that is, (C.7) is valid with $k = 0$, then
\[V(r) = b \; r^{(2-n + \chi _{0})/2} \; (1 + \ln r)^{a \; / \; \chi _{0}}
\; (1+o(1)), \; r \rightarrow \infty .\]

To derive more precise asymptotic formulas, we have to restrict
the asymptotic behavior of a potential. Comparing solutions of the
equations (3.4) with different potentials $q$ \cite[Chap.XI]{Hart}
and using formulas from \cite[Theorem 2 and its Corollary]{VeM1},
we obtain the following conclusions. The case \textbf{(A)} is
considered in Lemma C.2 and the case \textbf{(B)} in Lemma C.3.
These classes are defined in
Section 3 in the paragraph before Theorem 3.2. \\

\begin{lem}\hspace{-.09in}\textbf{.} If (C.7) holds and in addition
$\int ^{\infty } t^{-1}\left(t^2 q(t) - k\right)^2 dt < \infty $,
then, as $r \rightarrow \infty$,
\begin{equation}
V(r)= b \; r^{(2-n+ \chi _{k}) \; /\; 2} \exp \left\{\frac{1}{\chi
_{k}} \int ^r _1 \left(t^2 q(t) - k\right) \frac{dt}{t} \right\}
(1+\bar{\bar{o}}(1))
\end{equation}
and
\begin{equation}
W(r)= b \; r^{(2-n - \chi _{k}) \; /\; 2} \exp \left\{-
\frac{1}{\chi _{k}} \int ^r _1 \left(t^2 q(t) - k\right)
\frac{dt}{t} \right\} (1+\bar{\bar{o}}(1)).
\end{equation}
\hfill \textbf{$\diamondsuit$}\\
\end{lem}

\begin{cor}\hspace{-.09in}\textbf{.} If, in addition,
$r^{-1}(r^2 q(r) - k) \in L(1, \infty )$, then
\[V(r)= b \; r^{(2-n+ \chi _{k}) \; /\; 2}
(1+\bar{\bar{o}}(1)), \; r \rightarrow \infty ,\] and
\[W(r)= b \; r^{(2-n - \chi _{k}) \; /\; 2}
(1+\bar{\bar{o}}(1)), \; r \rightarrow \infty .\]
\end{cor}

\pagebreak

\begin{lem}\hspace{-.09in}\textbf{.} If
$\lim _{r \rightarrow \infty } s(r) = \infty $, then $V$ grows
faster than any power of $r$, that is, $ \lim _{ r \rightarrow
\infty } r^{-\alpha} V(r) = \infty $ for every $\alpha >0$.
Moreover, under some regularity conditions such as, for instance,
the potential $q \in C^2 (r_0, \infty )$, $\int ^{\infty }
\frac{dr}{r^2 \sqrt{q(r)}} < \infty $, and
\[\int ^{\infty} \biggl|4 \widetilde{q}(r) \widetilde{q}''(r) -
5 \left(\widetilde{q}'(r)\right)^2 \biggl| \widetilde{q}^{-5/2}(r)
dr < \infty ,\] where $\widetilde{q}(r) = q(r)+ \left(\frac{n^2
-4n +3}{4} + \lambda \right) r^{-2}$, there hold the following
\emph{ JWKB}-asymptotic formulas \cite[p. 382, Exercise
9.6]{Hart}, as $r \rightarrow \infty$,
\begin{equation}
V(r)= b \; r^{(1-n)/2}q^{-1/4}(r) \exp \left\{\int
^{r}_{1}q^{1/2}(t)dt \right\} \; (1+\overline{\overline{o}}(1))
\end{equation}
\begin{equation}
W(r)= b \; r^{(1-n)/2}q^{-1/4}(r) \exp \left\{- \int
^{r}_{1}q^{1/2}(t)dt \right\} \; (1+\overline{\overline{o}}(1)).
\end{equation} \hfill \textbf{$\diamondsuit$}
\end{lem}

The conditions of Lemma C.3 can be replaced by various
combinations of simpler conditions - see, for example,
\cite[Theorem 3]{VeM1}. For example, the following three
conditions are sufficient for (C.11)-(C.12) to hold:
\[(\beta ') \; \; \; r^{-4} q^{-3/2} (r) \in L(1, \infty ) \vspace{.2cm} \]
\[(\beta '') \; \; \; q^{-3/2} (r) q'' (r) \in L(1, \infty ) \vspace{.2cm} \]
\[(\beta ''') \; \; \; q^{-3/2} (r) q'(r) \rightarrow 0 \mbox{ as }
\; r \rightarrow \infty . \] \hfill \textbf{$\diamondsuit$}

\end{document}